\documentclass[3p,11pt]{elsarticle}
\usepackage{amsmath}

\usepackage{amsfonts}
\usepackage{amssymb}
\usepackage{mathrsfs}
\usepackage{amsfonts}
\usepackage{graphicx}
\usepackage{color}
\usepackage{graphics,color,epsfig}
\usepackage{stmaryrd}
\usepackage{enumitem}
\usepackage{amsthm}

\newcommand{\px}[1][x]{\partial_{#1}}
\newcommand{\dx}[1][x]{\,{\rm d}#1}

\newtheorem{example}{Example}[section]
\newtheorem{remark}{Remark}[section]
\newcommand{\sjwchange}[1]{\textcolor{black}{#1}}

\journal{Journal of Computational Physics}
 \graphicspath{{figures/}}

\begin{document}

\begin{frontmatter}



\title{A discrete least squares collocation method for two-dimensional nonlinear
time-dependent partial differential equations}




\author{Fanhai Zeng\footnotemark[2], Ian Turner\footnotemark[2]${}^{,}$\footnotemark[3]${}^{,}$\footnotemark[1],
Kevin Burrage\footnotemark[2]${}^{,}$\footnotemark[4], Stephen J. Wright\footnotemark[5]}
\date{\today}

\renewcommand{\thefootnote}{\arabic{footnote}}
\renewcommand{\thefootnote}{\fnsymbol{footnote}}
\footnotetext[2]{School  of  Mathematical   Sciences,
Queensland   University  of  Technology,   Brisbane,   QLD 4001, Australia (fanhaiz@foxmail.com).}
\footnotetext[3]{Australian Research Council Centre of Excellence for Mathematical and Statistical Frontiers, Queensland University of
Technology, Brisbane, QLD 4001, Australia (i.turner@qut.edu.au).}
\footnotetext[4]{Visiting Professor, Department of  Computer  Science,     University  of
Oxford, OXI 3QD, UK (kevin.burrage@qut.edu.au).}
\footnotetext[5]{Computer Sciences Department, University of Wisconsin, Madison, WI 53706
(swright@cs.wisc.edu).}
\footnotetext[1]{Corresponding author.}

\begin{abstract}
In this paper, we develop regularized discrete least squares
collocation and finite volume methods for solving two-dimensional
nonlinear  {time-dependent} partial
differential equations on irregular domains. The solution is
approximated using tensor product cubic spline basis functions defined
on a background rectangular (interpolation) mesh, which leads to high
spatial accuracy and straightforward implementation, and establishes a
solid base for extending the computational framework to
three-dimensional problems. A semi-implicit time-stepping method is
employed to transform the nonlinear partial differential equation into
a linear boundary value problem. A key finding of our study is that
the newly proposed mesh-free finite volume method based on circular
control volumes reduces to the collocation method as the radius limits
to zero. Both methods produce a large constrained least-squares
problem that must be solved at each time step in the advancement of
the solution. We have found that regularization yields a relatively
well-conditioned system that can be solved accurately using QR
factorization. An extensive numerical investigation is performed to
illustrate the effectiveness of the present methods, including the
application of the new method to a coupled system of time-fractional
partial differential equations having different fractional indices in
different (irregularly shaped) regions of the solution domain.
\end{abstract}

\end{frontmatter}


\section{Introduction}\label{sec1}
\sjwchange{Meshfree and meshless methods have been
applied widely to solve partial differential equations (PDEs) on
irregular domains due to their flexibility in dealing with a wide
variety of geometries}, see, for
example, \cite{BELYTSCHKO1996,LiuGu2005,NGUYEN2008}.  In this work, we
take an alternative approach to deal with irregular domains and
develop  {efficient} regularized
discrete least squares (LS) collocation and finite volume methods
(FVMs) for solving
nonlinear  {time-dependent} partial
differential equations in  \sjwchange{irregular two-dimensional} domains.

Our work is motivated by the need to resolve multiscale transport
processes in heterogeneous porous media using upscaling methods. As an
example, the Extended Distributed Microstructure Model proposed
in \cite{CarrTP16} approximates the macroscopic flux as the average of
the microscopic fluxes computed over micro-cells having geometrical
features representative of the actual porous micro-structure. Such an
approach avoids the need for any effective parameters in the
formulation and accounts more accurately for a non-equilibrium field
evolution within the micro-cells. It is therefore important to have
fast and flexible computational methods for computing over complex,
irregularly shaped domains (micro-cells) that are often established
directly from images of the porous microstructure. Recent work
highlights that memory effects and non-Fickian behaviour are important
physical mechanisms evident for lignocellulosic materials such as wood
\cite{TurnerIP11}. The contribution of our work is to investigate new time-fractional models for use in our microscale formulation that can address these issues for regions of differing fractional properties within the micro-cell.

LS Galerkin and collocation methods have been widely applied to solve
PDEs due to their computational advantages and
simplifications, see for example \cite{BocGun98,DingShuYX14,Ernest76,GelbPR08,KeithPFD17,LarSchHer17,Daniele17}.
In the mesh-based (LS) Galerkin method for solving FDEs on irregular domains,
the domain is approximately divided into subdomains, on which
 basis functions are defined. The division of the irregular domain may cause errors that affect the accuracy of the method.
If moving boundary  conditions are considered in \sjwchange{this method, the
irregular domain must be redivided
and the mass and stiffness matrices  recomputed. These are  costly operations.}

The LS collocation (LSC) method with a completely orthogonal
computational mesh has been used to solve PDEs on irregular domains.
\sjwchange{Its advantages include simple implementation and high accuracy; see} \cite{BocGun98,Ernest76,Jiangbonan12,LaiblePinder89,LaiblePinder93,ZeiLaiPin95}
for further details.  However, these  methods may yield an
ill-conditioned LS system, which may cause computational difficulties
in real applications.  In this work, this disadvantage is overcome by
using a generalized regularization
technique \cite{GolubHansenO99,NocedalWright06}, which stabilizes the
solution procedure while preserving high accuracy.

{  The idea of developing the LSC method for solving
time-dependent PDEs has not been fully explored in previous work (see,
for example, \cite{Ernest76}). The main contribution of our work is to
combine this idea with the linearized time-stepping method to develop
regularized (penalized) LS methods for solving nonlinear
time-dependent PDEs on irregular domains, yielding linear systems that
can be solved efficiently.}

We first develop the regularized LSC method and least squares finite
volume method (LSFVM) for solving a two-dimensional linear
boundary-value problem on \sjwchange{an} irregular domain $\Omega$
(see \eqref{s31:eq-1}).
The key idea of the present LS method is
to \sjwchange{approximate the solution  by} the tensor product of cubic spline basis
functions defined on a rectangular mesh, for a division of a suitable rectangular
domain $\Omega_{\Box}\supset \Omega$,  leading to high spatial accuracy, straightforward implementation,
and easy extension to a three-dimensional framework. The use of regularization makes the penalized LS system relatively well-conditioned with  the smallest singular value of order
$O(\sqrt{\delta})$, where $\delta>0$ is
the regularization parameter that balances accuracy and well-conditioning
of the LS system. We use QR decomposition to solve the LS problem accurately and stably.
Numerical simulations show that a regularization parameter
$\delta=0.01$ yields satisfactory numerical solutions.


We then extend our LSC/LSFVM to solve a nonlinear time-dependent PDE.
\sjwchange{Our approach is to
apply} a semi-implicit time-stepping method (see \cite{ZengTBK2018})
to transform the nonlinear time-dependent PDE into a linear boundary
value problem, which is solved by the penalized LSC method or LSFVM
developed for the model problem.  {  One advantage of the
present LS method is that it  {requires solution of {\em
linear} LS systems}.}  \sjwchange{We present numerical simulations} to
verify the effectiveness of the present LS method, including a
fractional model with two fractional indices that models diffusion in
a composite medium.

A finite element method based on weighted extended B-splines on a
regular grid as basis functions was proposed to solve Dirichlet
problems on irregular domains in \cite{Hollig2001}, yielding a
well-conditioned linear system.  Recently, a least squares radial
basis function partition-of-unity method has been developed
in \cite{LarSchHer17}, which can deal with irregular domains easily.
The newly developed meshfree LSFVM in \cite{FOY2017} is another
alternative to deal with irregular domains and achieves high accuracy
in solving multi-phase porous media models. A key finding of the work
presented here is that as the radius of the circular finite volumes
approaches zero, the LSC method is recovered.

This paper is organized as follows. Section~\ref{sec-2} presents the
penalized LSC method and LSFVM for a two-dimensional linear model
problem.  We also investigate the conditioning of the LS system and
seek techniques for solving it efficiently.  In
Section~\ref{sec3-fpde}, we show how to extend the LS method for the
linear model problem to nonlinear time-dependent PDEs.  Two numerical
examples are given in Section~\ref{sec:numerical} to show the
effectiveness of the LSC method for solving nonlinear generalized PDEs
and a coupled system of time-fractional PDEs. We make some closing
remarks in Section~\ref{concl}.

\section{Discrete least squares for a model problem}\label{sec-2}

\begin{figure}[!h]
\begin{center}
\begin{minipage}{0.47\textwidth}\centering
\epsfig{figure=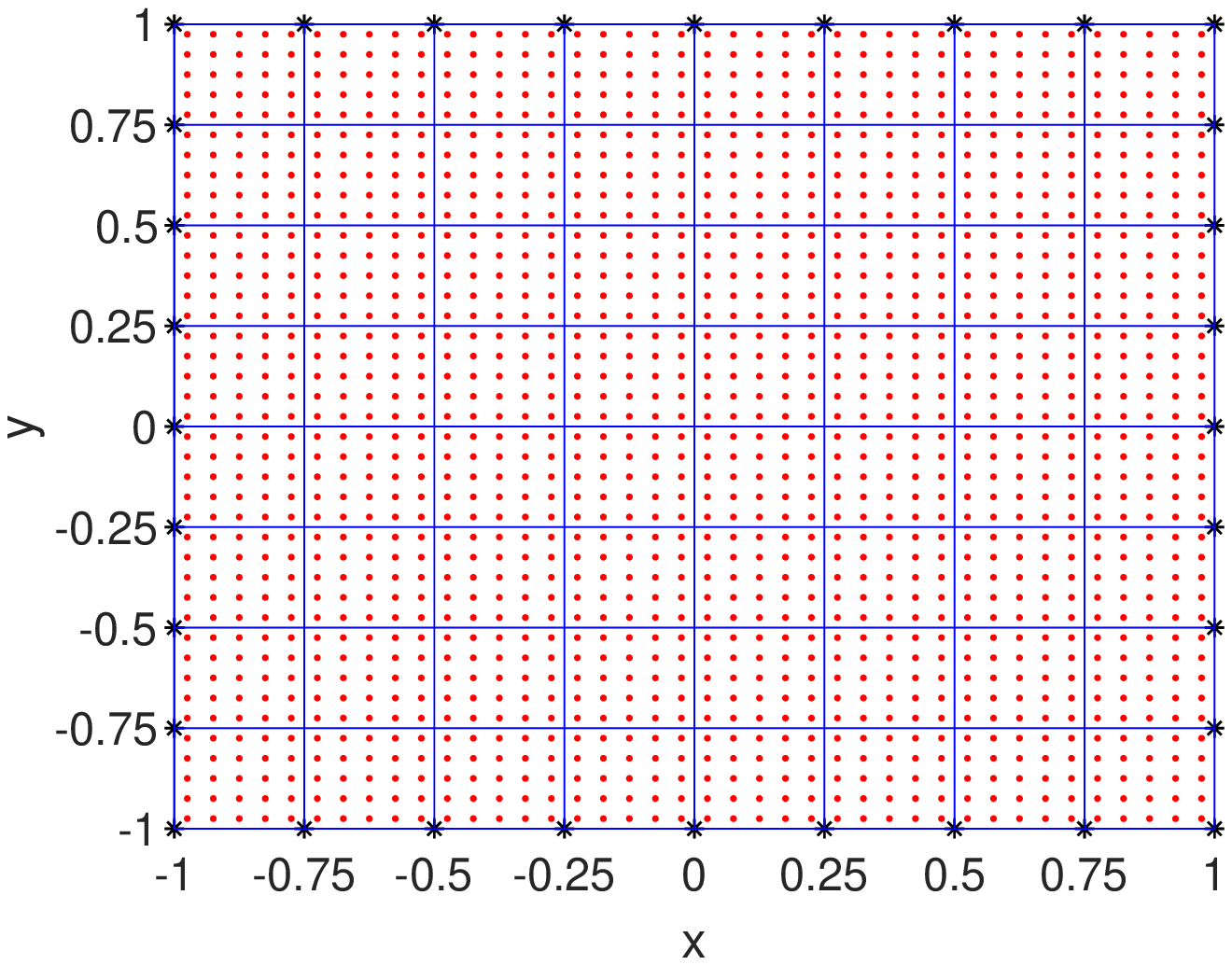,width=6cm} \par {(a) }
\end{minipage}
\begin{minipage}{0.47\textwidth}\centering
\epsfig{figure=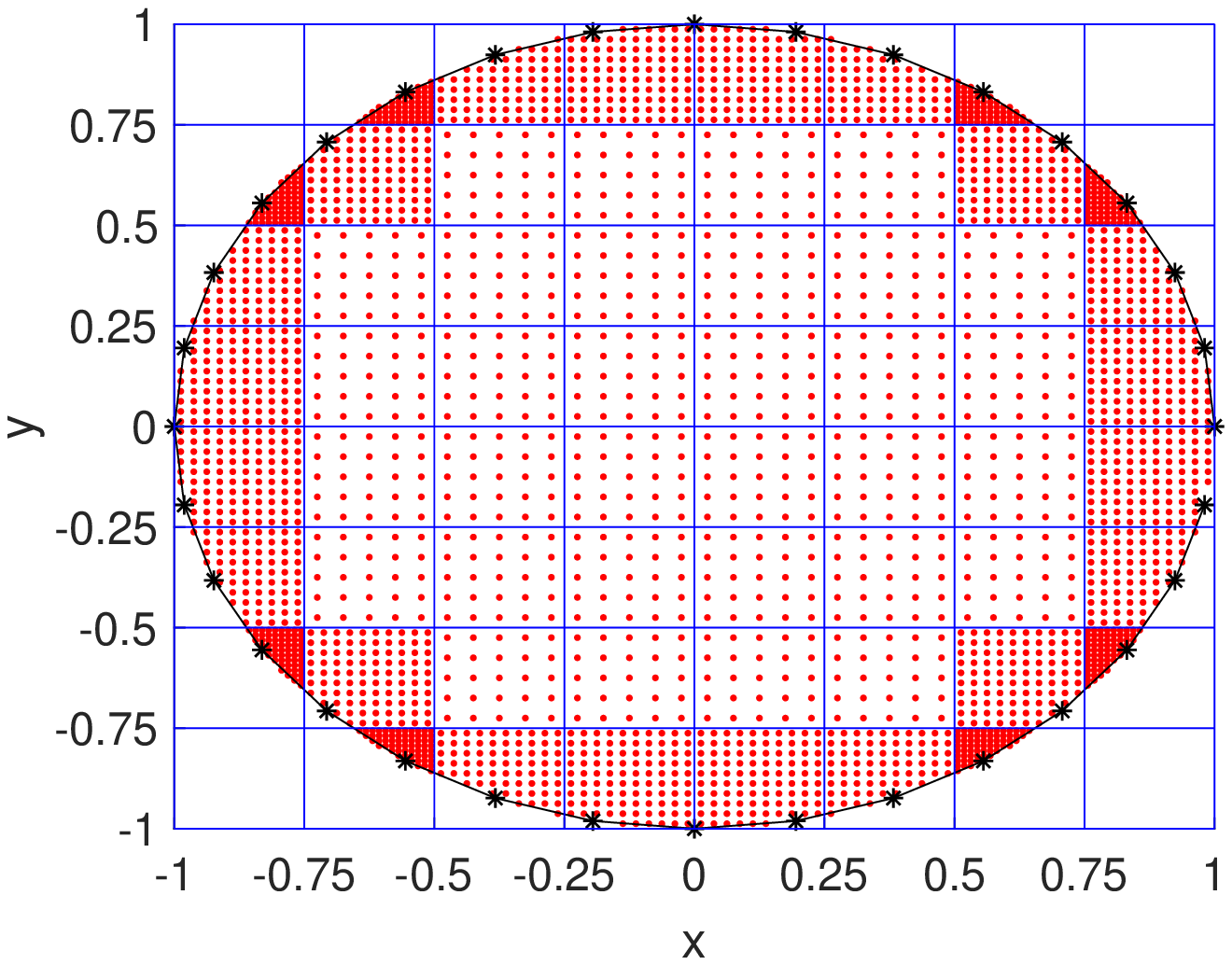,width=6cm} \par {(b)   }
\end{minipage}
\begin{minipage}{0.47\textwidth}\centering
\epsfig{figure=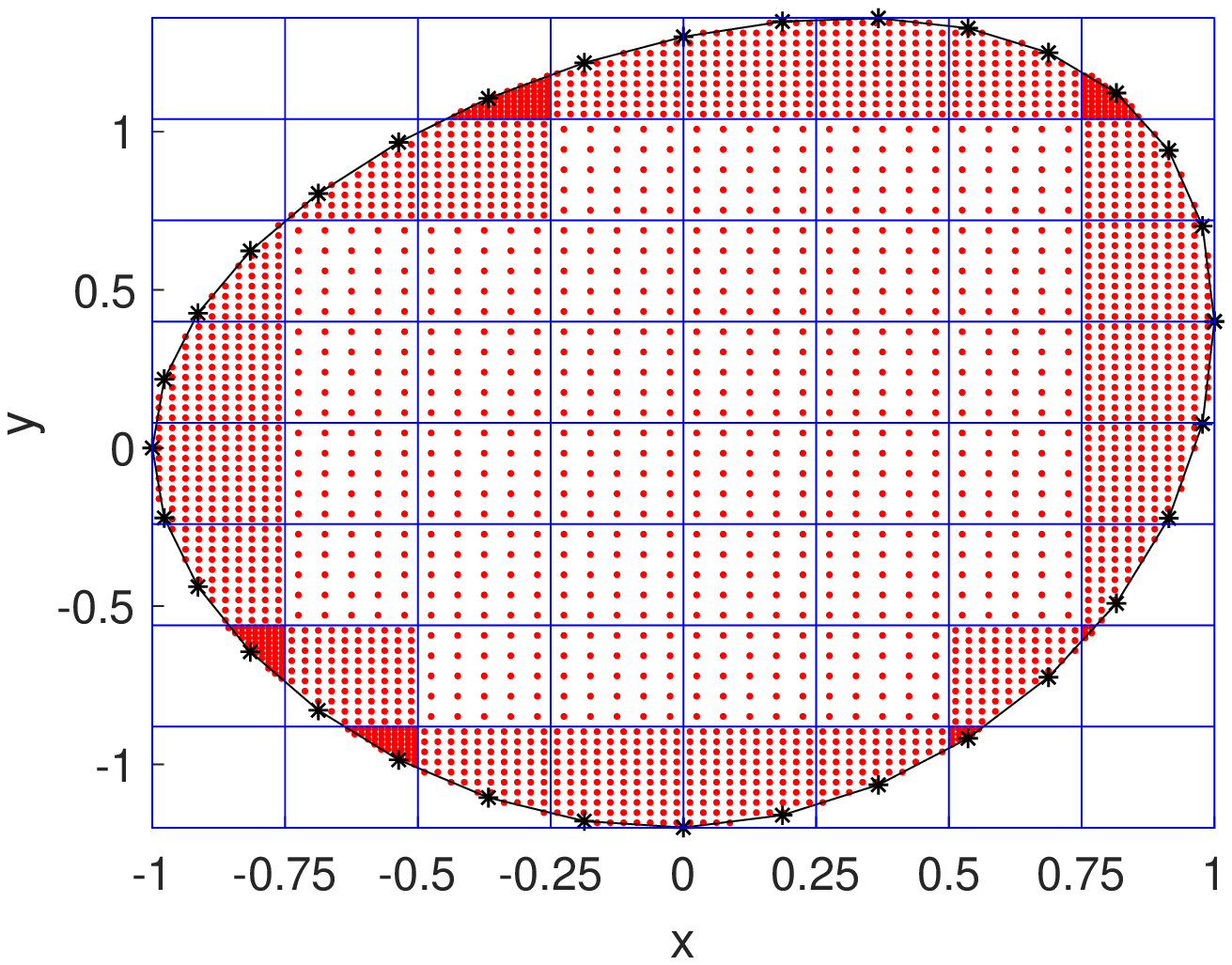,width=6cm} \par {(c)   }
\end{minipage}
\begin{minipage}{0.47\textwidth}\centering
\epsfig{figure=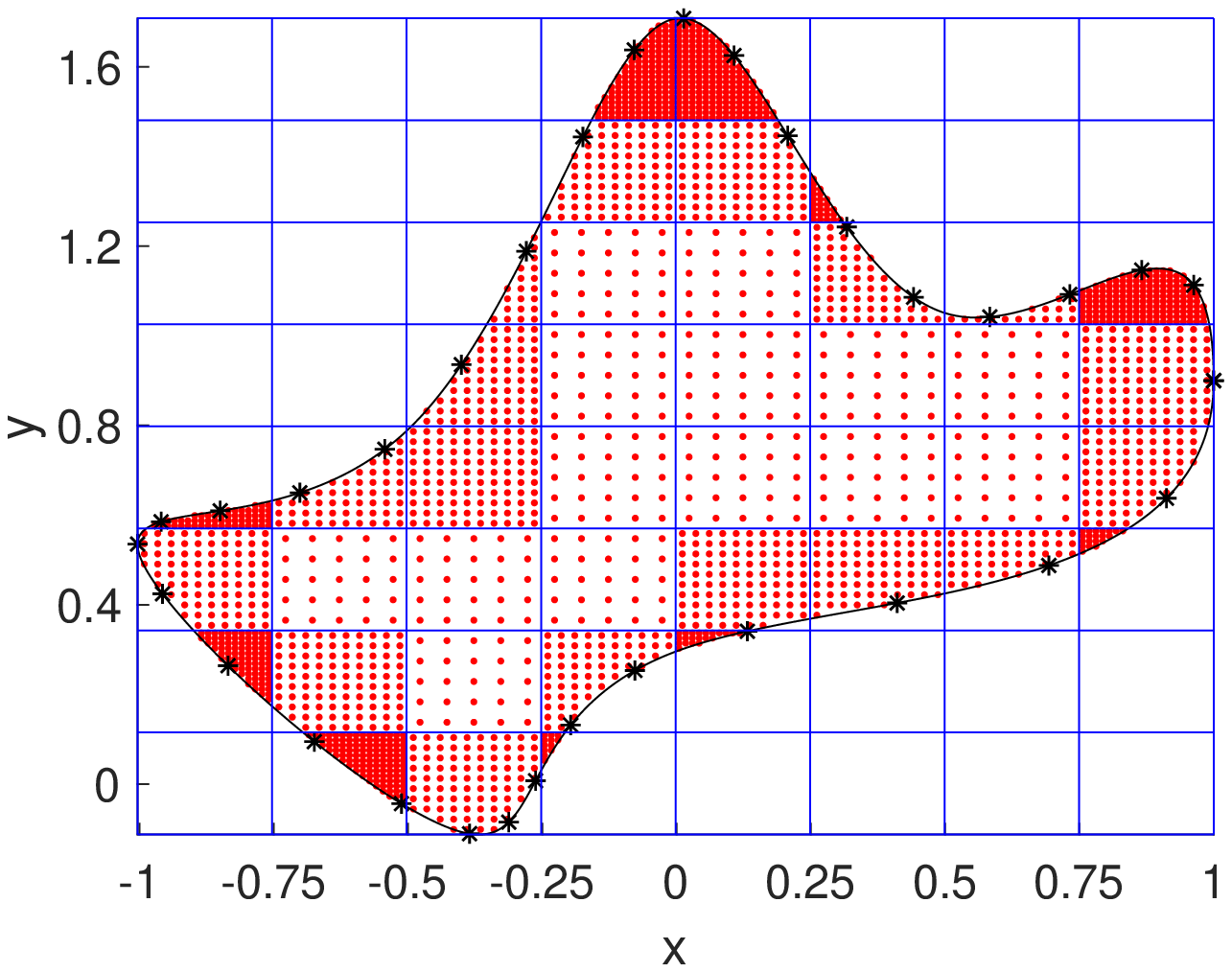,width=6cm} \par {(d)   }
\end{minipage}
\end{center}
\caption{The division of the different domains
(the area surrounded by the black curve for (b), (c), and (d)) with
the corresponding boundary points (stars $*$) and collocation points
(red dots { $\bullet$}).\label{fig:iregular_domain}}
\end{figure}

We develop the regularized LSC method and LSFVM for the following
model problem
\begin{equation}\label{s31:eq-1}
u(x,y) - \nu\Delta u(x,y) = f(x,y),{\quad}(x,y)\in\Omega,
\end{equation}
subject to the following generalized boundary condition
\begin{equation}\label{s31:bcs}
\mathcal{B}(u)= u_b \quad \text{on}\quad  \partial\Omega,
\end{equation}
where $\Delta $ is the Laplacian in two-dimensional space; $\nu$ is a
positive \sjwchange{scalar} parameter; $\mathcal{B}$ represents
Dirichlet, Neumann, or mixed boundary conditions; and
$\Omega$  {is an irregular domain with a
piecewise smooth boundary;} Figure~\ref{fig:iregular_domain} shows
several irregular domains that will be used for the numerical
\sjwchange{computations} performed in this work to verify our discretisation
methods.

\subsection{Notation}

A background (interpolation) mesh \sjwchange{and a corresponding set
of} basis functions are needed for the LS method.  For any finite
domain $\Omega$, we can find a suitable rectangular domain
$\Omega_{\Box}$ such that $\Omega \subset \Omega_{\Box}=(a,b)\times
(c,d)$; see Figure \ref{fig:iregular_domain}. The default values of
$a,b,c,$ and $d$ are given by $a=\inf\{x|(x,y)\in \Omega\}$,
$b=\sup\{x|(x,y)\in \Omega\}$, $c=\inf\{y|(x,y)\in \Omega\}$, and
$d=\sup\{y|(x,y)\in \Omega\}$, respectively.

Let $\delta_x=\{x_i\}$ (or $\delta_y=\{y_j\}$) be a division of $I^x$
(or $I^y$), satisfying $a=x_0<x_1<\cdots<x_{N_x}=b$ (or
$c=y_0<y_1<\cdots<y_{N_y}=d$).  Denote $I_{i,j}={I}^x_i \times
{I}^y_j$, where $I^x_i=(x_{i-1},x_{i})$ and $I^y_j=(y_{j-1},y_{j})$.
Then the domain ${\Omega}_{\Box}$ is divided into non-overlapping
rectangular subdomains that satisfy
${\bar{\Omega}_{\Box}}=\cup_{i,j} \bar{I}_{i,j}$.

The $C^1$ continuous basis function space defined on the mesh
$\delta_x\times\delta_y$ is given by
\begin{equation}
\mathcal{M}_{\Omega}(\delta_x\times\delta_y)=\left\{v \, : \, v\in C^1(\bar{\Omega}), \;
v|_{\bar{I}_{i,j}} \in P_3(\bar{I}^x_i)\otimes P_3(\bar{I}^y_j), \; I_{i,j} \cap \Omega \neq \emptyset \right\},
\end{equation}
where $P_k(\bar{I})$ denotes the polynomial space defined on the
interval $\bar{I}$ with degree no greater than $k$.  In this work, we
apply the tensor product cubic spline basis functions defined on the
mesh grids $\delta_x\times\delta_y$ to approximate the solution of the
differential equation.  Therefore, any function
$v\in\mathcal{M}_{\Omega}(\delta_x\times\delta_y)$ can be expressed by
$v|_{\bar{I}_{i,j}}
= \sum_{m=1}^4\sum_{n=1}^4c^{(i,j)}_{m,n}\phi^{(i)}_m(x)\psi^{(j)}_n(y)$,
where the cubic spline basis functions $\phi^{(i)}_m(x)$ in the $x$
direction are given by (see \cite{WWSun1999}):
\begin{equation}
\begin{aligned}
\phi^{(i)}_1(x) &=  2\hat{x}^3-3\hat{x}^2+1,
{\qquad}\phi^{(i)}_2(x) =  (x_i-x_{i-1})(\hat{x}^3-2\hat{x}^2+\hat{x}), \\
\phi^{(i)}_3(x) &=  -2\hat{x}^3+3\hat{x}^2,  {\qquad\quad}
\phi^{(i)}_4(x) =  (x_i-x_{i-1})(\hat{x}^3-\hat{x}^2),\\
 \hat{x} &= (x-x_{i-1})/(x_{i}-x_{i-1}).
\end{aligned}\end{equation}
The basis functions $\psi^{(j)}_n(y)(1\leq n \leq 4)$ in the $y$ direction
are defined similarly.
%

Before presenting the LS method, we introduce further notation. Let
$$\{\Phi_1(x,y),\Phi_2(x,y),\dotsc,\Phi_M(x,y)\}$$ be a
basis for $\mathcal{M}_{\Omega}(\delta_x\times\delta_y)$.
Denote
\begin{align}
\Phi(x,y)&=(\Phi_1(x,y),\Phi_2(x,y),\dotsc,\Phi_M(x,y))^T,\\
\Delta\Phi(x,y) &= (\Delta\Phi_1(x,y),\Delta\Phi_2(x,y),\dotsc,\Delta\Phi_M(x,y))^T,
\end{align}
with  $\Phi^T(x,y)=(\Phi(x,y))^T$ and $\Delta\Phi^T(x,y)=(\Delta\Phi(x,y))^T$.
For any $\mathbf{c},\mathbf{c}^n\in \mathbb{R}^M$, denote $\mathbf{c}=(c_1,c_2,\dotsc,c_M)^T$
and $\mathbf{c}^n=(c_1^n,c_2^n,\dotsc,c_M^n)^T$.
Then any $U_h\in \mathcal{M}_{\Omega}(\delta_x\times\delta_y)$ can be expressed by
  \begin{equation}\label{s31:Uh}\begin{aligned}
U_h(x,y) = \Phi^T(x,y)\mathbf{c}= (\Phi^T\mathbf{c})(x,y),{\quad}\mathbf{c}\in\mathbb{R}^M.
\end{aligned}\end{equation}

\sjwchange{The set $\mathcal{S}^{(p,q)}_{\Omega}$  of collocation points  in the domain $\Omega$ is defined by the following four steps.}
\begin{itemize}
  \item Step 1) Let $p$ and $q$ be two positive integers.  Define
  $\widehat{\mathcal{S}}_{p,q}$ as
      $$\widehat{\mathcal{S}}_{p,q}
      =\left\{(\hat{x}_k,\hat{y}_{\ell}) \, : \, 1\leq k\leq p,\; 1\leq \ell \leq q\right\},$$
 where $(\hat{x}_k,\hat{y}_{\ell})=(\frac{2k-1}{2p},\frac{2\ell-1}{2q})$\footnote{\sjwchange{We can make other choices for these points}, provided that  $0<\hat{x}_1<\hat{x}_2< \cdots< \hat{x}_p<1$ and $0<\hat{y}_1<\hat{y}_2< \cdots< \hat{y}_q<1$.}.
Obviously,  $\widehat{\mathcal{S}}_{p,q}$ is a set of points uniformly distributed
in the reference domain $(0,1)\times (0,1)$.
  \item  Step 2) For any $1\leq i\leq N_x$ and $1\leq j\leq N_j$, define the mapping
  $\mathcal{P}_{i,j}$ as follows:
      \begin{equation}\label{Pij}
\mathcal{P}_{i,j}(\widehat{\mathcal{S}}_{p,q})
=\{(\widetilde{x}_k,\widetilde{y}_{\ell}) \, : \, \widetilde{x}_k={h_i^x\hat{x}_k}+x_{i-1}, \;
\widetilde{y_\ell}={h_j^y\hat{y}_{\ell}}+y_{j-1}, \; (\hat{x}_k,\hat{y}_{\ell})\in \widehat{\mathcal{S}}_{p,q}\},
\end{equation}
where $h_i^x=x_i-x_{i-1}$ and $h_j^y=y_j-y_{j-1}$.

\item   Step 3)  Define  $ {\mathcal{S}}^{(i,j)}_{p,q}$   as follows:
    \begin{equation}\label{Sijpq}
    {\mathcal{S}}^{(i,j)}_{p,q} =\left\{\begin{aligned}
&\mathcal{P}_{i,j}(\widehat{\mathcal{S}}_{5,5}),&&{\quad} I_{i,j} \subset \Omega,\\
&\mathcal{P}_{i,j}(\widehat{\mathcal{S}}_{p,q}), &&{\quad}  \text{two or three vertices of the rectangle } I_{i,j} \text{ are in } \Omega,\\
&\mathcal{P}_{i,j}(\widehat{\mathcal{S}}_{2p,2q}),&&{\quad}  \text{only one vertex of the rectangle } I_{i,j} \text{ is in } \Omega,\\
& \emptyset,&&{\quad} \text{otherwise}.
\end{aligned}\right.\end{equation}
The set ${\mathcal{S}}^{(i,j)}_{10,10}$ is applied if we do not specify $p$ and $q$.

  \item Step 4) The set $\mathcal{S}^{(p,q)}_{\Omega}$  of collocation points is given by
      \begin{equation}\label{Spq}
    \mathcal{S}^{(p,q)}_{\Omega}
    = \Big(\bigcup_{i,j} {\mathcal{S}}^{(i,j)}_{p,q}\Big)\bigcap \Omega.
    \end{equation}
\end{itemize}
\sjwchange{The set} $\mathcal{S}^{(p,q)}_{\Omega}$ \sjwchange{has} more collocation points in the
cells $I_{i,j}$ near the boundary of $\Omega$, which makes the
coefficient matrix from the LS method more well-conditioned than the
matrix based on the uniform distribution of collocation point set
$\mathcal{\overline{S}}^{(p,q)}_{\Omega}$  \sjwchange{defined by}
\begin{equation}\label{Spq2}
\mathcal{\overline{S}}^{(p,q)}_{\Omega}
 = \Big(\bigcup_{i,j} \mathcal{P}_{i,j}(\widehat{\mathcal{S}}_{p,q})\Big)\bigcap \Omega.
\end{equation}

\begin{figure}[!h]
\begin{center}
\begin{minipage}{0.47\textwidth}\centering
\epsfig{figure=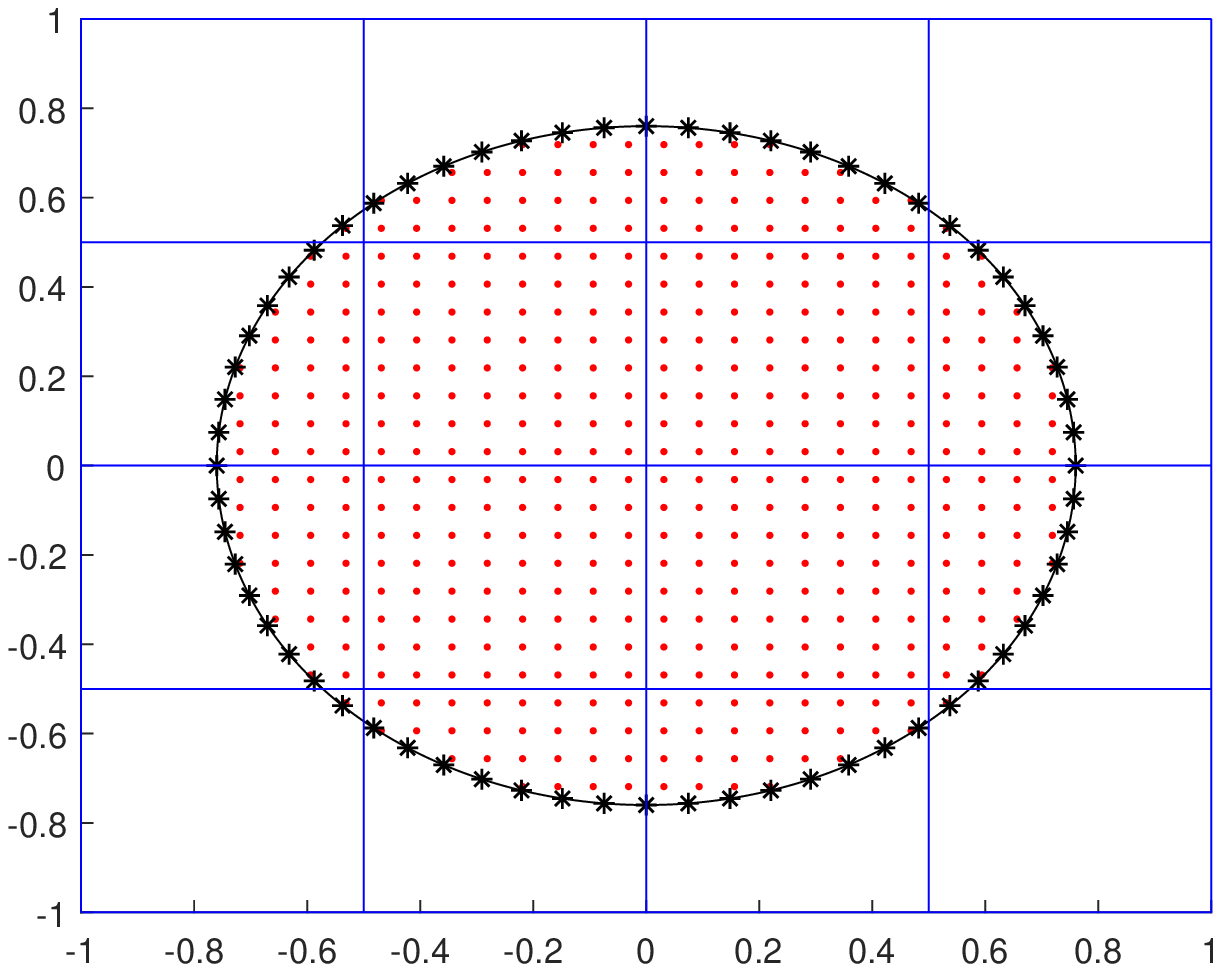,width=7cm}
\par {(a) Distribution of  $\mathcal{\overline{S}}^{(8,8)}_{\Omega}$. }
\end{minipage}
\begin{minipage}{0.47\textwidth}\centering
\epsfig{figure=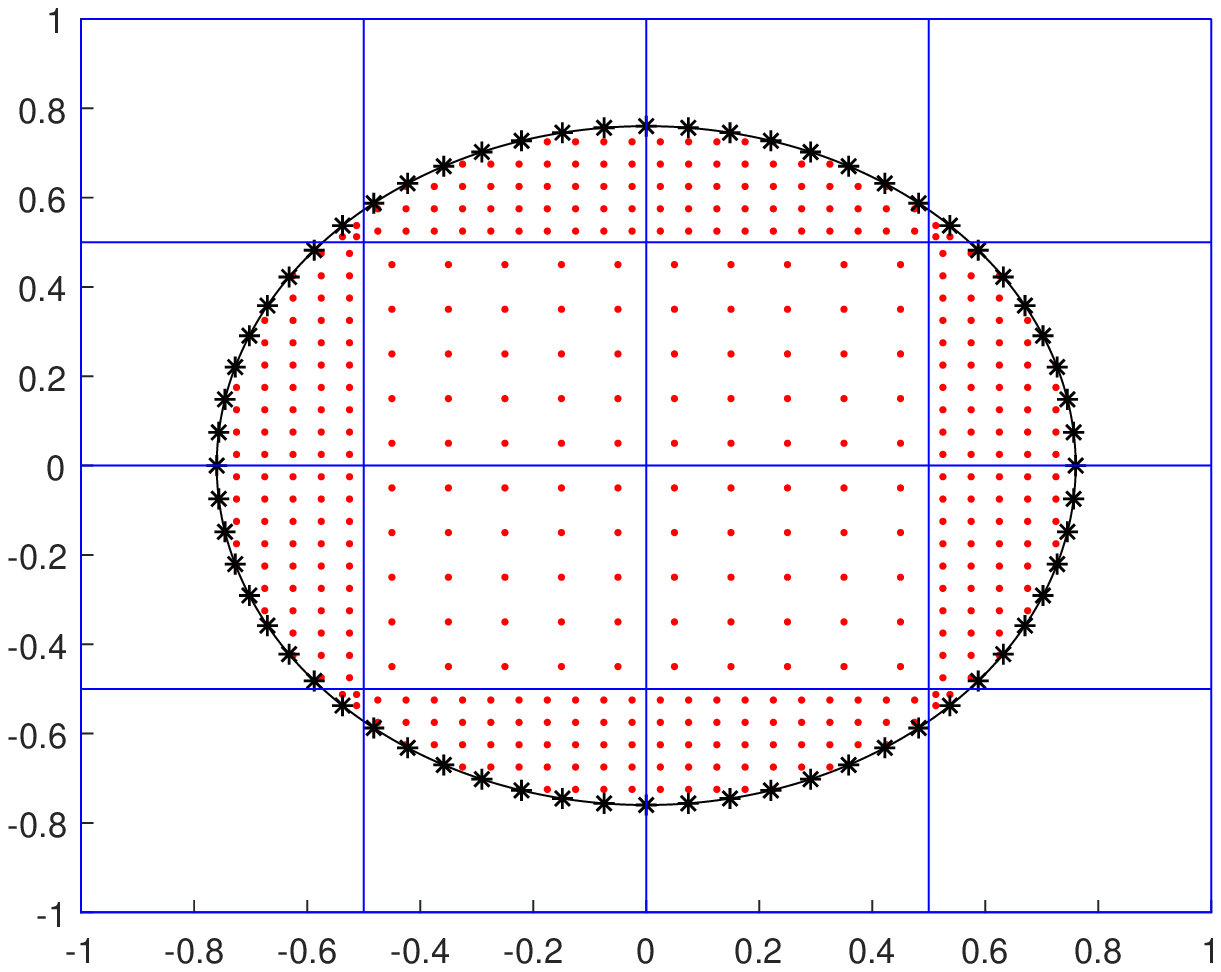,width=7cm}
\par {(b) Distribution of $\mathcal{S}^{(8,8)}_{\Omega}$. }
\end{minipage}
\end{center}
\caption{Distribution of collocation points of  $\mathcal{\overline{S}}^{(8,8)}_{\Omega}$  and
$\mathcal{S}^{(8,8)}_{\Omega}$,
where $\Omega =\{(x,y)|x^2+y^2\leq 0.76^2\}$,
$\Omega_{\Box} = (-1,1)^2$, and $N_x=N_y=4$.\label{fig2-0-0}}
\end{figure}

For simplicity, we denote
$\xi^{(in)}_i=(x^{(in)}_i,y^{(in)}_i)\in \mathcal{S}^{(p,q)}_{\Omega}$,
$1\leq i\leq N_{in}$ .  The set of collocation points on
$\partial\Omega$ is denoted by $\mathcal{S}^{(p,q)}_{\partial\Omega}$
with
$\xi^{(b)}_i=(x^{(b)}_i,y^{(b)}_i)\in\mathcal{S}^{(p,q)}_{\partial\Omega}$,
$1\leq i \leq N_b$.  Figure~\ref{fig:iregular_domain} shows the
distribution of collocation and boundary points in different domains.
{ Figure \ref{fig2-0-0} contrasts
$\mathcal{\overline{S}}^{(p,q)}_{\Omega}$ and
$\mathcal{S}^{(p,q)}_{\Omega}$ for $p=q=8$.}

\subsection{LSC method and LSFVM}

\sjwchange{We now} present the penalized LSC method and  LSFVM
for solving the model problem \eqref{s31:eq-1}, yielding a
better-conditioned linear system that can be solved by many known
methods.

\subsubsection{LSC   method}
Replacing $u(x,y)$ in \eqref{s31:eq-1} with $U_h(x,y)
= \Phi^T(x,y)\mathbf{d}$ defined by \eqref{s31:Uh}, and letting
$(x,y)=\xi^{(in)}_i\in\mathcal{S}^{(p,q)}_{\Omega}$, we obtain
\begin{equation}\label{s31:eq-10}
\left(\Phi^T(\xi^{(in)}_i)  - \nu\Delta \Phi^T(\xi^{(in)}_i)\right)\mathbf{d}
= f(\xi^{(in)}_i),{\quad}1\leq i \leq N_{in},
\end{equation}
subject to \sjwchange{boundary conditions}
\begin{equation}\label{s31:eq-11}
\mathcal{B}(U_h(\xi^{(b)}_i))=\mathcal{B}(\Phi^T(\xi^{(in)}_i))\mathbf{d}
=u_b(\xi^{(b)}_i),{\quad}1\leq i \leq N_{b}.
\end{equation}
Clearly, the system \eqref{s31:eq-10}-\eqref{s31:eq-11} is
overdetermined. If we enforce boundary conditions exactly and other
conditions in a LS sense, we obtain the following constrained LS
problem
\begin{equation}\label{s31:LeastSquare}\begin{aligned}
\mathbf{c} &= \arg\min_{\mathbf{d}\in  \mathbb{R}^{M}} \,
\sum_{i=1}^{N_{in}}\left((  \Phi^T(\xi^{(in)}_i)
-\nu \Delta\Phi^T(\xi^{(in)}_i))\mathbf{d}  -f(\xi^{(in)}_i) \right)^2\\
& = \arg\min_{\mathbf{d}\in  \mathbb{R}^{M} } \,
\|(\mathbf{A}_{in}-\nu\mathbf{S}_{in})\mathbf{d}-\mathbf{f}_{in}\|^2, \\
& \mbox{subject to} {\quad} \mathbf{A}_b\mathbf{d}=\mathbf{u}_b
=\left(u_b(\xi^{(b)}_1),u_b(\xi^{(b)}_2),
\dotsc,u_b(\xi^{(b)}_{N_{b}})\right)^T,
\end{aligned}\end{equation}
where  $\mathbf{f}_{in}=(f(\xi^{(in)}_1),f(\xi^{(in)}_2),
\dotsc,f(\xi^{(in)}_{N_{in}}))^T$, \sjwchange{and the matrices}
{ $\mathbf{A}_{in},\mathbf{S}_{in} \in R^{N_{in}\times M}$ and
$\mathbf{A}_{b}\in R^{N_{b}\times M}$ are given by}
\begin{equation}\label{s31:A}
\mathbf{A}_{in}=\left[
 \begin{array}{cccc}
\Phi^T(\xi^{(in)}_1)  \\
\Phi^T(\xi^{(in)}_2)  \\
                 \vdots \\
\Phi^T(\xi^{(in)}_{N_{in}})  \\
\end{array}
\right],{\quad}
\mathbf{S}_{in}=\left[
 \begin{array}{cccc}
\Delta\Phi^T(\xi^{(in)}_1)  \\
\Delta\Phi^T(\xi^{(in)}_2)  \\
                 \vdots \\
\Delta\Phi^T(\xi^{(in)}_{N_{in}})  \\
\end{array}
\right], {\quad}
\mathbf{A}_{b}=\left[
 \begin{array}{cccc}
\mathcal{B}(\Phi^T(\xi^{(b)}_1))  \\
\mathcal{B}(\Phi^T(\xi^{(b)}_2))  \\
                 \vdots \\
\mathcal{B}(\Phi^T(\xi^{(b)}_{N_{b}}))  \\
\end{array}
\right].
\end{equation}

One way to enforce boundary conditions approximately is to include
them directly in the overdetermined system using a weight $\lambda$
that can be varied according to how accurately we want the conditions
to be satisfied
(see \cite{Loan85}). This \sjwchange{quadratic-penalty} approach
yields the following unconstrained LS formulation:
\begin{equation}\label{s31:LeastSquare-2}\begin{aligned}
\mathbf{c}_{\lambda}&= \arg\min_{\mathbf{d}\in  \mathbb{R}^{M} }
\left\{\|(\mathbf{A}_{in}-\nu\mathbf{S}_{in})\mathbf{d}-\mathbf{f}_{in}\|^2
+ \lambda^2\|\mathbf{A}_b\mathbf{d}-\mathbf{u}_b\|^2\right\}\\
&=\arg\min_{\mathbf{d}\in  \mathbb{R}^{M} }\left\|\left(
           \begin{array}{c}
             \mathbf{A}_{in}-\nu\mathbf{S}_{in} \\
             \lambda\mathbf{A}_b \\
           \end{array}
         \right)\mathbf{d}
         -\left(
           \begin{array}{c}
             \mathbf{f}_{in} \\
             \lambda\mathbf{u}_b  \\
           \end{array}
         \right)
\right\|^2.
\end{aligned}\end{equation}
\sjwchange{Large values of $\lambda$ increase the importance of satisfying the
 boundary residuals relative to the interior
residuals},  {that is, the boundary
conditions will be more accurately satisfied as $\lambda$ becomes
larger;} see \cite{Ernest76,Loan85}.

There is computational difficulty in solving
both \eqref{s31:LeastSquare} and \eqref{s31:LeastSquare-2} due to the
ill-conditioning (even rank deficiency) of the matrices
$\mathbf{A}_{in}-\nu\mathbf{S}_{in}$ and
$\binom{\mathbf{A}_{in}-\nu\mathbf{S}_{in}}{\lambda\mathbf{A}_b}$.
\sjwchange{This is caused} by a possibly very small support of some cells
$I_{i,j}$ in $\Omega$, which means that
$|I_{i,j}\cap\Omega|/|I_{i,j}|$ is too small to contain enough
collocation points, where $|I_{i,j}|$ denotes the volume of the domain
$I_{i,j}$; see \cite{Hollig2001}. We tackle such ill-conditioning by
formulating the following penalized LS problem
\begin{equation}\label{s31:LeastSquare2}\begin{aligned}
\mathbf{c}_{\delta} &= \arg\min_{\mathbf{d}\in  \mathbb{R}^{M} }\left\{
\|(\mathbf{A}_{in}-\nu\mathbf{S}_{in})\mathbf{d}-\mathbf{f}_{in}\|^2
+ \delta\|{\mathbf{M}}(\mathbf{d}- \mathbf{d_*})\|^2\right\}, \\
&\mbox{subject to} {\quad} \mathbf{A}_b\mathbf{d}=\mathbf{u}_b,
\end{aligned}\end{equation}
where $\delta\geq 0$, $\mathbf{M}\in \mathbb{R}^{M\times M}$, and
$\mathbf{d_*}$ is a reference solution (\sjwchange{a prior estimate of
$\mathbf{d}$}), so that $\|{\mathbf{M}}(\mathbf{d}-\mathbf{d_*})\|$ is
small.  The penalty term
$\delta\|{\mathbf{M}}(\mathbf{d}- \mathbf{d_*})\|^2$ balances
well-conditioning and accuracy in the method \eqref{s31:LeastSquare2}
for
solving \eqref{s31:eq-1}.  \sjwchange{Formulation \eqref{s31:LeastSquare2}
can be viewed as a type of Tikhonov
regularization} \cite{GolubHansenO99,Neumaier98}.  In this work, we
use $\mathbf{M}=\mathbf{I}$ (the identity matrix), which appears to
work well. \sjwchange{Other possible choices for $\mathbf{M}$ are
described in}   \cite{FOY2017,Neumaier98}.

The regularized analog of \eqref{s31:LeastSquare-2} is as follows
\begin{equation}\label{s31:LeastSquare-22}\begin{aligned}
\mathbf{c}_{\lambda,\delta}&= \arg\min_{\mathbf{d}\in  \mathbb{R}^{M} }
\left\{\|(\mathbf{A}_{in}-\nu\mathbf{S}_{in})\mathbf{d}-\mathbf{f}_{in}\|^2
+ \lambda^2\|\mathbf{A}_b\mathbf{d}-\mathbf{u}_b\|^2
+ \delta\|\mathbf{d}-\mathbf{d_*}\|^2\right\}\\
&=\arg\min_{\mathbf{d}\in  \mathbb{R}^{M} }
\left\|\left(
           \begin{array}{c}
             \mathbf{A}_{in}-\nu\mathbf{S}_{in} \\
             \lambda\mathbf{A}_b \\
             \sqrt{\delta}\mathbf{I}\\
           \end{array}
         \right)\mathbf{d} -\left(
           \begin{array}{c}
             \mathbf{f}_{in} \\
             \lambda\mathbf{u}_b  \\
             \sqrt{\delta}\mathbf{d_*}\\
           \end{array}
         \right)\right\|^2.
\end{aligned}\end{equation}
\sjwchange{Approaches for solving \eqref{s31:LeastSquare2}
and \eqref{s31:LeastSquare-22} are discussed in} Section~\ref{sec223}.
We will focus on applying \eqref{s31:LeastSquare-22} to the model
problem \eqref{s31:eq-1}, then extend it to solve nonlinear
time-dependent PDEs in Section~\ref{sec3-fpde}.

\subsubsection{LSFVM}\label{sec2-2-2}

Let $V_i$ be a control volume centered at an interior point
$\xi_i^{(in)}=(x^{(in)}_i,y^{(in)}_i) \in \mathcal{S}^{(p,q)}_{\Omega}$
\sjwchange{with  radius} $\rho>0$. Suppose that
$U_h\in\mathcal{M}_{\Omega}(\delta_x\times\delta_y)$ is an approximate
solution of $u(x,y)$ to \eqref{s31:eq-1}. \sjwchange{Gauss's
Divergence Theorem yields}
\begin{equation}\label{s5:eq-2}\begin{aligned}
\int_{V_i}U_h\dx[V_i]
&=\nu \int_{\px[]V_i}\nabla U_h\cdot \mathbf{\hat{n}}\dx[s_i]+\int_{V_i}f(x,y)\dx[V_i],
\end{aligned}\end{equation}
where $\mathbf{\hat{n}}$ is the unit normal on the boundary
$\px[]V_i$.  Let $x=x_i^{(in)}+\rho\cos\theta$ and
$y=y_i^{(in)}+\rho\sin\theta$, so that $\dx[s_i]=\rho\dx[\theta]$.
\sjwchange{By substituting from \eqref{s31:Uh} for $U_h$ in $\int_{\px[]V_i}\nabla
U_h\cdot \mathbf{\hat{n}}\dx[s_i]$, we rewrite \eqref{s5:eq-2} as
follows:}
\begin{equation}\label{s5:eq-3}\begin{aligned}
 \int_{V_i}U_h\dx[V_i]
=\nu\rho\int_{0}^{2\pi}J^{(i)}(\theta)\dx[\theta]+\int_{V_i}f(x,y)\dx[V_i],
\end{aligned}\end{equation}
where $J^{(i)}(\theta)=J^{(i)}_x(\theta)+J^{(i)}_y(\theta)$ and
\begin{align*}
J^{(i)}_x(\theta)&= \px[x]U_h(x_i^{(in)}+\rho\cos\theta,y_i^{(in)}+\rho\sin\theta)\cos\theta,\\
J^{(i)}_y(\theta)&=\px[y]U_h(x_i^{(in)}+\rho\cos\theta,y_i^{(in)}+\rho\sin\theta)\sin\theta.
\end{align*}
By Taylor expansion at $\xi_i^{(in)}$, the integral
$\int_{V_i}U_h\dx[V_i]$ in \eqref{s5:eq-3} can be calculated exactly
by
\begin{equation}\label{s5:eq-4}
\int_{V_i}U_h\dx[V_i]= \pi\rho^2U_h(\xi_i^{(in)})+\frac{\pi\rho^4}{8}{\Delta}U_h(\xi_i^{(in)}),
\end{equation}
provided that $\rho$ is suitably small.
\sjwchange{Since
$J^{(i)}(\theta)$ is a trigonometric polynomial of degree at most six,}
 the second integral in \eqref{s5:eq-3}
can also be evaluated exactly using
the trapezoidal formula with at least six quadrature points, that is,
\begin{equation}\label{s5:eq-5}
\int_{0}^{2\pi}J^{(i)}(\theta)\dx[\theta]
=\frac{2\pi}{K}\sum_{r=1}^K J^{(i)}(\theta_r),
{\quad}\theta_r=\frac{2r\pi}{K},K\geq 6.
\end{equation}
By combining \eqref{s5:eq-3}, \eqref{s5:eq-4}, and \eqref{s5:eq-5}, we have
\begin{equation}\label{s5:eq-6}
\left(1+\frac{\rho^2}{8}{\Delta}\right)U_h(\xi_i^{(in)})
 =\frac{2\nu}{K\rho} \sum_{r=1}^KJ^{(i)}(\theta_r)
 +  f(\xi_i^{(in)})+O(\rho^2),\quad 1\leq i \leq N_{in},
\end{equation}
where we used $\int_{V_i}f\dx[V_i]= \pi\rho^2f(\xi_i^{(in)})+O(\rho^4).$
\sjwchange{By omitting the $O(\rho^2)$ term in this} equation, we  obtain
\begin{equation}\label{s5:eq-7}
\mathbf{\widetilde{A}}_{in}  \mathbf{c}
= \nu\mathbf{\widetilde{S}}_{in}\mathbf{c} + \mathbf{f}_{in},
\end{equation}
where
$\mathbf{\widetilde{A}}_{in}=(\mathbf{A}_{in}+\frac{\rho^2}{8}\mathbf{S}_{in})$,
and the matrix $\mathbf{\widetilde{S}}_{in}$ depends on $\rho$. If
$\rho$ is suitably small, then $\mathbf{\widetilde{S}}_{in}$ can also
be expressed by
\begin{equation*}
\mathbf{\widetilde{S}}_{in}=\mathbf{S}_{in}+\frac{\rho^2}{4}\left[
 \begin{array}{cccc}
\px^2\px[y]^2\Phi^T(\xi^{(in)}_1)  \\
\px^2\px[y]^2\Phi^T(\xi^{(in)}_2)  \\
                 \vdots \\
\px^2\px[y]^2\Phi^T(\xi^{(in)}_{N_{in}})  \\
\end{array}
\right].
\end{equation*}

The regularized LSFVM for \eqref{s31:eq-1}-\eqref{s31:bcs} is given by
\begin{equation}\label{s31:LeastSquare-33}\begin{aligned}
\mathbf{c}_{\delta} &= \arg\min_{\mathbf{d}\in  \mathbb{R}^{M} }\left\{
\|(\mathbf{\widetilde{A}}_{in}-\nu\mathbf{\widetilde{S}}_{in})\mathbf{d}-\mathbf{f}_{in}\|^2
+ \delta\|\mathbf{d}-\mathbf{d_*}\|^2\right\}, \\
&\mbox{subject to} {\quad} \mathbf{A}_b\mathbf{d}=\mathbf{u}_b,
\end{aligned}\end{equation}
where  $\mathbf{A}_{b}$, $\mathbf{f}_{in}$, and $\mathbf{u}_{b}$
are given in \eqref{s31:LeastSquare-2}, and $\mathbf{\widetilde{A}}_{in}$
and $\mathbf{\widetilde{S}}_{in}$ are defined in \eqref{s5:eq-7}.

By penalizing the constraint in \eqref{s31:LeastSquare-33}   as in \eqref{s31:LeastSquare-22}, we obtain
\begin{equation}\label{s31:LeastSquare-3}\begin{aligned}
\mathbf{c}_{\lambda,\delta}
&=\arg\min_{\mathbf{d}\in  \mathbb{R}^{M} }\left\|\left(
           \begin{array}{c}
             \mathbf{\widetilde{A}}_{in} -\nu\mathbf{\widetilde{S}}_{in} \\
             \lambda\mathbf{A}_b \\
             \sqrt{\delta}\mathbf{I}
           \end{array}
         \right)\mathbf{d}
         -\left(
           \begin{array}{c}
             \mathbf{f}_{in} \\
             \lambda\mathbf{u}_b  \\
             \sqrt{\delta}\mathbf{d_*}
           \end{array}
         \right)
\right\|.
\end{aligned}\end{equation}

\subsubsection{Solving the LS problems} \label{sec223}
In this subsection, we discuss how to solve the LS formulations
\eqref{s31:LeastSquare2} and  \eqref{s31:LeastSquare-22}. The LSFVMs
\eqref{s31:LeastSquare-33} and \eqref{s31:LeastSquare-3}
can be solved similarly.

\begin{figure}[!h]
\begin{center}
\begin{minipage}{0.45\textwidth}\centering
\epsfig{figure=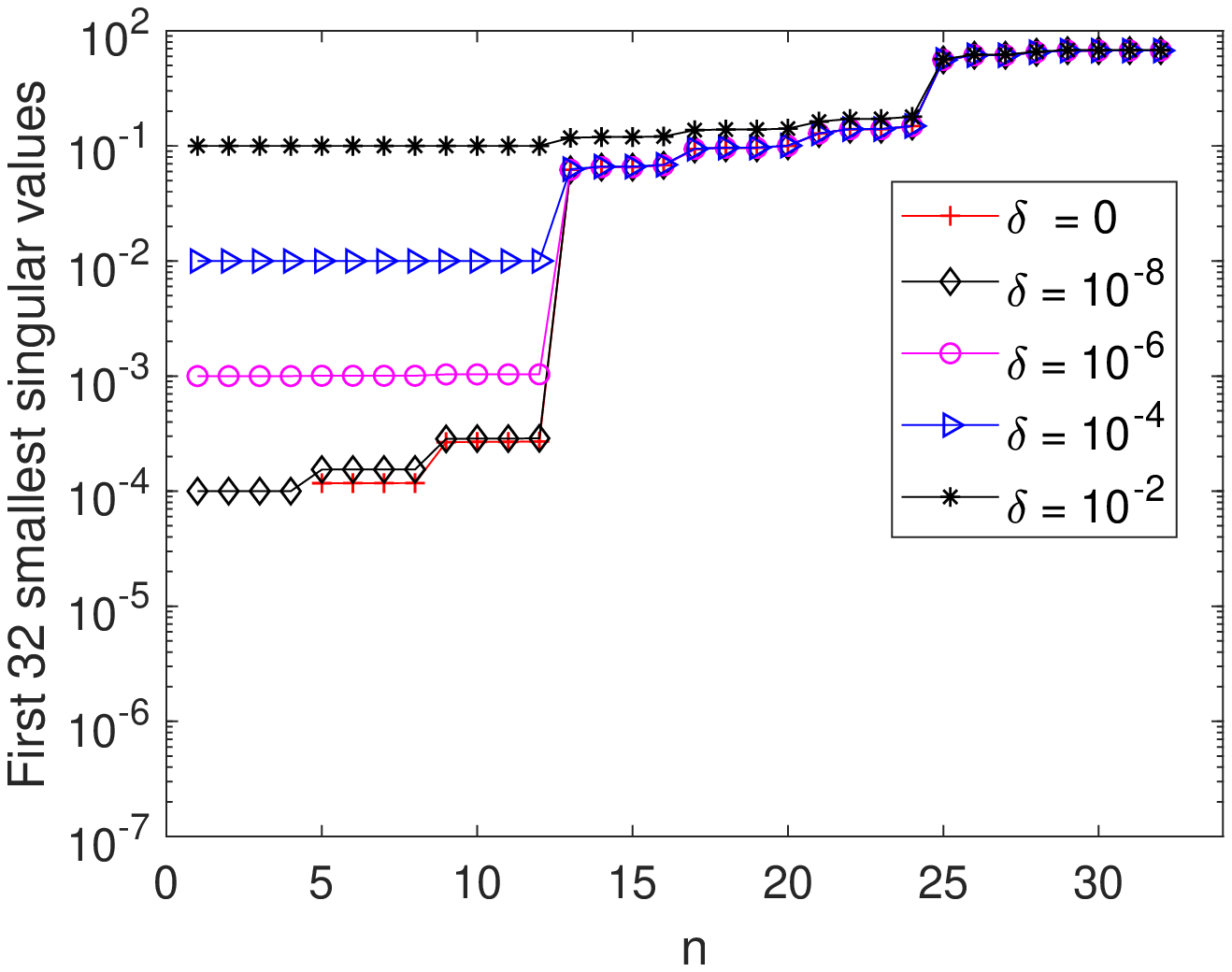,width=7cm}
\par {(a) Singular values based on $\mathcal{\overline{S}}^{(15,15)}_{\Omega}$.}
\end{minipage}
\begin{minipage}{0.45\textwidth}\centering
\epsfig{figure=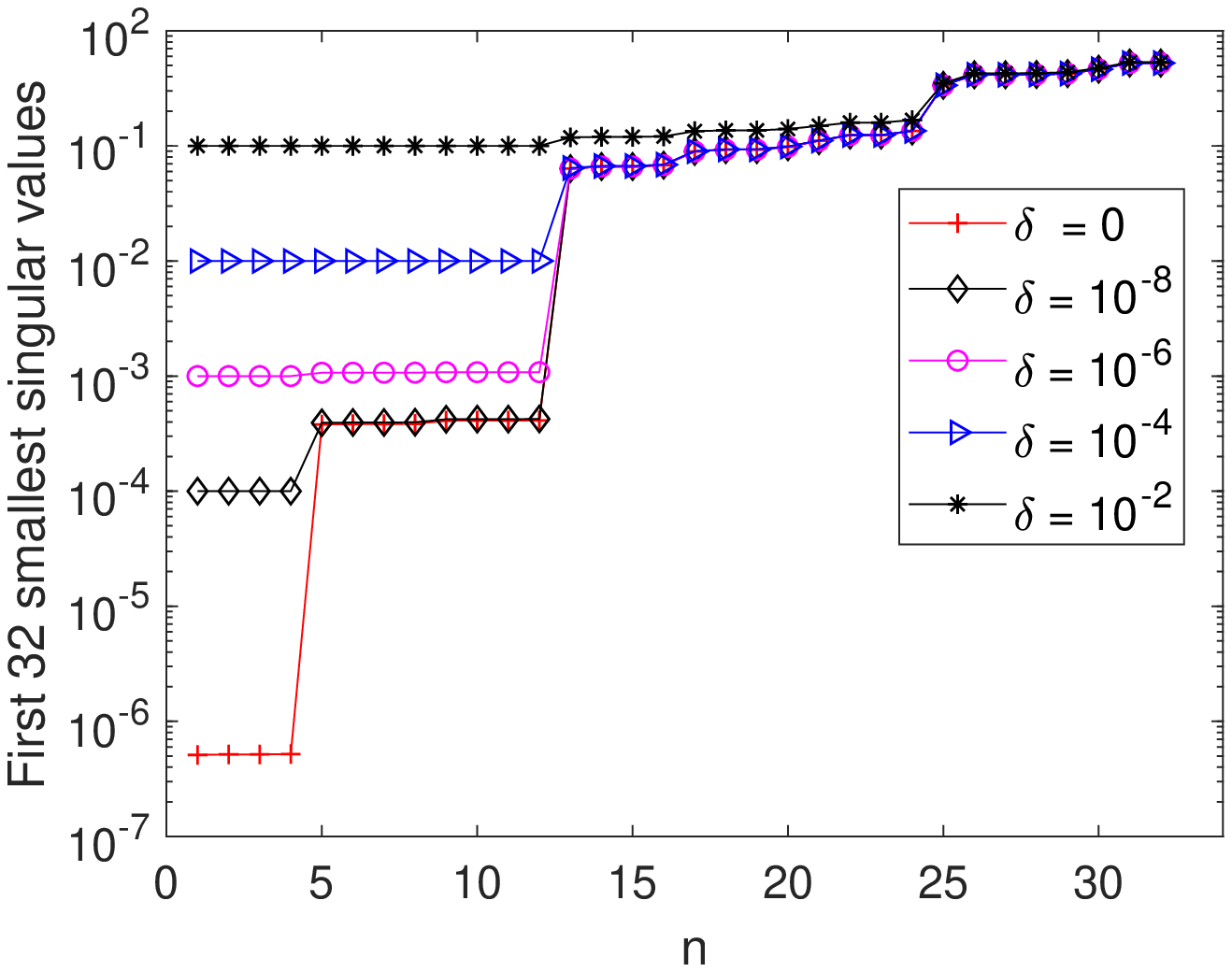,width=7cm}
\par {(b) Singular values based on $\mathcal{S}^{(15,15)}_{\Omega}$.}
\end{minipage}
\end{center}
\caption{The first 32 singular values   of the matrix \eqref{matrix-3} based on different densities of the collocation point sets  $\mathcal{\overline{S}}^{(p,q)}_{\Omega}$ (see \eqref{Spq}) and
$\mathcal{S}^{(p,q)}_{\Omega}$ (see \eqref{Spq2}), where $\Omega
=\{(x,y)|x^2+y^2\leq 0.76^2\}$, $\Omega_{\Box} = (-1,1)^2$,
$N_x=N_y=4$, and $\lambda=10^{4}$.\label{fig2-0}}
\end{figure}

{  The condition number of the linear LS system may depend
on both largest and smallest singular values of the coefficient matrix
as well as the solution and the \sjwchange{optimal residual} of the LS
system \cite{Grcar03}.  For example, the spectral norm absolute or
relative condition number of the LS system \eqref{s31:LeastSquare-22}
can be bounded by
\begin{equation}\label{cond-1-1}
\frac{\left(\|\mathbf{r}_0\|^2+\|\mathbf{d}_0\|^2(\sigma^2_{\min}+\delta)\right)^{1/2}}
{\sigma^2_{\min}+\delta}
\end{equation}
and
\begin{equation}\label{cond-1-2}
\left(\frac{\|\mathbf{r}_0\|^2}{\|\mathbf{d}_0\|^2 \sqrt{\sigma^2_{\min}+\delta}}+1\right)
\sqrt{\frac{\sigma^2_{\max}+\delta}{\sigma^2_{\min}+\delta}},
\end{equation}
\sjwchange{respectively (see \cite[Theorem~5.1]{Grcar03}).}
Here, $\sigma_{\max}$ and $\sigma_{\min}$ are the largest and smallest
nonzero singular values of
${\mathbf{A}_{in}-\nu\mathbf{S}_{in}}\choose{\lambda\mathbf{A}_b}$,
respectively; $\mathbf{d}_0$ is the solution of the LS
system \eqref{s31:LeastSquare-22}; and $\mathbf{r}_0$ is
the \sjwchange{optimal residual}.
\sjwchange{We see from \eqref{cond-1-1} and \eqref{cond-1-2}
that the smallest nonzero singular value of the coefficient matrix
in \eqref{s31:LeastSquare-22}, which is
\begin{equation}\label{matrix-3}\begin{aligned}
\left(
       \begin{array}{c}
             \mathbf{A}_{in}-\nu\mathbf{S}_{in} \\
             \lambda\mathbf{A}_b \\
             \sqrt{\delta}\mathbf{I}\\
           \end{array}
         \right),
\end{aligned}\end{equation}
plays an important role on the ill-conditioning of the LS
system \eqref{s31:LeastSquare-22}.}

The regularization parameter $\delta$ helps to improve the condition
number of system \eqref{s31:LeastSquare-22}.
Larger values of $\delta$ improve the conditioning
of \eqref{s31:LeastSquare-22}, but may lead to inaccurate numerical
solutions.  A choice of \sjwchange{reference solution $\mathbf{d_*}$}
that is not too far from the solution allows us to use a larger
\sjwchange{value of $\delta$ in \eqref{s31:LeastSquare-22}  and still} maintain accuracy.
\sjwchange{We could obtain $\mathbf{d_*}$ by using alternative methods to
solve \eqref{s31:LeastSquare}, or derive it by using our method on a
coarser grid.} We will show how different choices of $\delta$ and
$\mathbf{d_*}$ affect the accuracy of the
methods \eqref{s31:LeastSquare2} and \eqref{s31:LeastSquare-22} in
Section~\ref{sec2-4}.

{  If a division $\delta_x\times\delta_y$ of $\Omega_{\Box}$
has a cell $I_{i,j}$ with very small support in $\Omega$, then the
smallest singular value of \eqref{matrix-3} may be very small, leading
to a very large condition number of
\eqref{matrix-3}. \sjwchange{Consider the plots in  Figure~\ref{fig2-0}
with $\delta=0$. The} coefficient matrix \eqref{matrix-3} based on the
uniformly distributed collocation set
$\mathcal{\overline{S}}^{(15,15)}_{\Omega}$ (see \eqref{Spq2}) is
singular (see Figure \ref{fig2-0}(a)), while the coefficient matrix
based on the nonuniformly distributed collocation set
$\mathcal{S}^{(15,15)}_{\Omega}$ (see \eqref{Spq}) is nonsingular with
smallest singular value of magnitude of $O(10^{-7})$ (see
Figure \ref{fig2-0}(b)).  \sjwchange{The use of regularization
($\delta>0$)} increases the smallest singular value
$\sqrt{\sigma^2_{\min}+\delta}$ of \eqref{matrix-3}. The effects on
small singular values of the matrix \eqref{matrix-3} are apparent from
Figure \ref{fig2-0}.}  }

There are a number of methods documented in the literature for solving
the LS problem
\eqref{s31:LeastSquare2},
even if \eqref{s31:LeastSquare2} is
ill-conditioned \cite{Neumaier98,ScottTuma17}.  The constrained
formulation \eqref{s31:LeastSquare2} can be solved by quadratic
programming techniques  \sjwchange{or
by writing the optimality conditions as a linear system with symmetric
indefinite matrix, as follows:}
\begin{equation} \label{eq:kkt}
\left[ \begin{array}{ccc}
 {\mathbf{I}} & {-(\mathbf{A}_{in}-\nu\mathbf{S}_{in})} & {\mathbf{0}} \\
 {-(\mathbf{A}_{in}-\nu\mathbf{S}_{in})^{T}} & {-\delta \mathbf{I}} & {(\mathbf{A}_b)^{T}} \\
 {\mathbf{0}} & {\mathbf{A}_b} & {\mathbf{0}}\end{array}\right]
 \left[ \begin{array}{l}{\mathbf{r}} \\
 {\mathbf{d}} \\
 {\mathbf{\gamma}}\end{array}\right]
 =\left[ \begin{array}{c}{-\mathbf{f}_{in}} \\
 {-\delta \mathbf{d}^{*}} \\ \mathbf{u}_b\end{array}\right],
\end{equation}
\sjwchange{where $\gamma$ is a vector of Lagrange multipliers for
the constraints in \eqref{s31:LeastSquare2}.}
In the numerical simulations \sjwchange{below}, we
show that the constrained formulation \eqref{s31:LeastSquare2} and
penalized formulation \eqref{s31:LeastSquare-22} yield numerical
solutions of similar accuracy.
In the following, we discuss how to solve the unconstrained LS
problem \eqref{s31:LeastSquare-22}, which is extended to solve
time-fractional PDEs in Section \ref{sec3-fpde}.

\sjwchange{The solution of \eqref{s31:LeastSquare-22} can be obtained by forming
and solving its normal equations, which are as follows:}
\begin{equation}\label{s223-2}\begin{aligned}
((\mathbf{A}_{in}-\nu\mathbf{S}_{in})^T(\mathbf{A}_{in}
&-\nu\mathbf{S}_{in})+\lambda^2\mathbf{A}_b^T\mathbf{A}_b
+ \delta \mathbf{I})\mathbf{c} \\
&=(\mathbf{A}_{in}-\nu\mathbf{S}_{in})^T\mathbf{f}_{in}
+ \lambda^2\mathbf{A}_b^T\mathbf{u}_b + \delta\mathbf{d_*}.
\end{aligned}\end{equation}
\sjwchange{The condition number of the coefficient matrix in this system is
approximately} ${\sigma^2_{\max}}/{\delta}$.
\sjwchange{Iterative methods such as conjugate gradients can be used to
solve \eqref{s223-2}, or it can be solved directly using a Cholesky
factorization.}  Even when the linear system \eqref{s223-2} is
ill-conditioned, it is still possible to obtain accurate solutions
with an appropriate choice of
algorithms \cite{CarsonHigham18,EldenSimoncini12,Neumaier98,ScottTuma17}. However,
our computational experiments show that those methods for
solving \eqref{s223-2} are slower than \sjwchange{methods that apply a QR
factorization directly to the matrix in
\eqref{s31:LeastSquare-22}; see Table~\ref{tb2-1}.
The QR-based approach can be implemented compactly by means of the
backslash operator in Matlab.}


\subsection{Error analysis}
Let $\mathbf{c}$, $\mathbf{c}_{\lambda}$, $\mathbf{c}_{\delta}$, and $\mathbf{c}_{\lambda,\delta}$ be the solutions of \eqref{s31:LeastSquare},
\eqref{s31:LeastSquare-2}, \eqref{s31:LeastSquare2}, and
\eqref{s31:LeastSquare-22}, respectively.
\sjwchange{Our goal in this section is to estimate the difference norm
$\|\mathbf{c}-\mathbf{c}_{\lambda,\delta}\|$ between the regularized
solution $\mathbf{c}_{\lambda,\delta}$ and the true solution
$\mathbf{c}$.}

Assume that the rank of the matrix
$\mathbf{A}_{in}-\nu\mathbf{S}_{in}$ is $M$.  Then the error bounds
$\|\mathbf{c}-\mathbf{c}_{\lambda}\|$ and
$\|\mathbf{c}_{\delta}-\mathbf{c}_{\lambda,\delta}\|$ can be derived
directly \sjwchange{from \cite[equation~(2.20)]{Loan85}.  We have}
\begin{equation}\label{sec2-3:eq-2}
\|\mathbf{c}-\mathbf{c}_{\lambda}\|
\leq \frac{1}{2\lambda\beta_{\min}}\sqrt{\|\mathbf{r}_{\lambda}\|^2-\|\mathbf{r}\|^2},
\end{equation}
where $\beta_{\min}$ is the smallest singular value of $\mathbf{A}_b$ and
\begin{equation*}\label{sec2-3:eq-3}
\mathbf{r}= \left(\mathbf{A}_{in}-\nu\mathbf{S}_{in}\right)\mathbf{c} - \mathbf{f}_{in} ,{\qquad}
\mathbf{r}_{\lambda}=\mathbf{A}_{\lambda}\mathbf{c}_{\lambda} - \mathbf{b}_{\lambda},
\end{equation*}
with
$\mathbf{A}_{\lambda}=\left(
\begin{array}{c}
\mathbf{A}_{in}-\nu\mathbf{S}_{in} \\
\lambda\mathbf{A}_b \\
\end{array}
 \right)$ and
 $\mathbf{b}_{\lambda}=\left(
           \begin{array}{c}
             \mathbf{f}_{in} \\
             \lambda\mathbf{u}_b  \\
           \end{array}
         \right)$.
Numerical simulations show  that $\mathbf{c}_{\lambda}$
is a good approximation of $\mathbf{c}$ for a sufficiently large value $\lambda>0$, that is,
 $\mathbf{c}_{\lambda}$ is bounded if   $\mathbf{c}$ is bounded; see \cite{Loan85}.

Next, we estimate $\mathbf{c}_{\lambda}-\mathbf{c}_{\lambda,\delta}$.
\sjwchange{We denote by $\hat{\sigma}_i$ the  singular values of
$\mathbf{A}_{\lambda}$, with
$0<\hat{\sigma}_{\min}=\hat{\sigma}_M\leq\hat{\sigma}_{M-1}\leq\cdots\leq\hat{\sigma}_1$.}
Then there exist orthogonal matrices $\mathbf{U}$ and $\mathbf{V}$
such that $\mathbf{A}_{\lambda}=\mathbf{U}\Sigma \mathbf{V}^T$, where
$\Sigma=\mathrm{diag}(\hat{\sigma}_1,\hat{\sigma}_2,\dotsc,\hat{\sigma}_M)$,
$\mathbf{U}=(\mathbf{u}_1,\mathbf{u}_2,\dotsc,\mathbf{u}_{M})$,
$\mathbf{V}=(\mathbf{v}_1,\mathbf{v}_2,\dotsc,\mathbf{v}_M)$,
satisfying $\mathbf{U}^T\mathbf{U}=\mathbf{V}^T\mathbf{V}=\mathbf{I}$.
The LS solutions of \eqref{s31:LeastSquare-2} and
\eqref{s31:LeastSquare-22} can be expressed by
\begin{align*}
\mathbf{c}_{\lambda}&=(\mathbf{A}_{\lambda}^T\mathbf{A}_{\lambda})^{-1}
\mathbf{A}_{\lambda}^T\mathbf{b}_{\lambda},\\
\mathbf{c}_{\lambda,\delta}&=(\mathbf{A}_{\lambda}^T\mathbf{A}_{\lambda}+\delta\mathbf{I})^{-1}
\left(\mathbf{A}_{\lambda}^T\mathbf{b}_{\lambda}+\delta\mathbf{d}_{*}\right).
\end{align*}
From the above two equations and $\mathbf{A}_{\lambda}=\mathbf{U}\Sigma \mathbf{V}^T$, we have
\begin{equation}\label{sec2-3:eq-4}\begin{aligned}
\|\mathbf{c}_{\lambda}-\mathbf{c}_{\lambda,\delta}\|
&=\|(\mathbf{A}_{\lambda}^T\mathbf{A}_{\lambda})^{-1}
\mathbf{A}_{\lambda}^T\mathbf{b}_{\lambda}-(\mathbf{A}_{\lambda}^T\mathbf{A}_{\lambda}+\delta\mathbf{I})^{-1}
\left(\mathbf{A}_{\lambda}^T\mathbf{b}_{\lambda}+\delta\mathbf{d}_{*}\right)\|\\
&=\delta\|(\mathbf{A}_{\lambda}^T\mathbf{A}_{\lambda})^{-1}
(\mathbf{A}_{\lambda}^T\mathbf{A}_{\lambda}+\delta\mathbf{I})^{-1}
\mathbf{A}_{\lambda}^T\mathbf{b}_{\lambda}
- (\mathbf{A}_{\lambda}^T\mathbf{A}_{\lambda}+\delta\mathbf{I})^{-1}\mathbf{d}_{*}\|\\
&=\delta\|\Sigma^{-1}(\Sigma^2 + \delta \mathbf{I})^{-1}\mathbf{U}^T\mathbf{b}_{\lambda}
- (\Sigma^2 + \delta \mathbf{I})^{-1} \mathbf{V}^T\mathbf{d}_{*}\|\\
& \leq \delta\sum_{i=1}^M
\left(\frac{|\mathbf{u}_i^T\mathbf{b}_{\lambda}|}{\hat{\sigma}_i(\hat{\sigma}_i^2+\delta)}
+\frac{|\mathbf{v}_i^T\mathbf{d}_{*}|}{(\hat{\sigma}_i^2+\delta)}\right).
\end{aligned}\end{equation}
This bound shows that $\mathbf{c}_{\lambda,\delta}$ converges to
$\mathbf{c}_{\lambda}$ as $\delta\to0$, but it does not show the
effectiveness of $\mathbf{d}_{*}$ in improving the accuracy of the
numerical solutions. However, by writing $\mathbf{c}_{\lambda}$ as
$\mathbf{c}_{\lambda}=(\mathbf{A}_{\lambda}^T\mathbf{A}_{\lambda}
+\delta\mathbf{I})^{-1}
\left(\mathbf{A}_{\lambda}^T\mathbf{b}+\delta\mathbf{c}_{\lambda}\right)$, we obtain
\begin{equation}\label{sec2-3:eq-6}
\|\mathbf{c}_{\lambda}-\mathbf{c}_{\lambda,\delta}\|
=\delta\|(\mathbf{A}_{\lambda}^T\mathbf{A}_{\lambda}+\delta\mathbf{I})^{-1}
\left(\mathbf{c}_{\lambda}-\mathbf{d}_{*}\right)\|
\leq\frac{\delta}{\delta+\hat{\sigma}^2_{\min}}
\|\mathbf{c}_{\lambda}-\mathbf{d}_{*}\|.
\end{equation}
By combing \eqref{sec2-3:eq-2} and \eqref{sec2-3:eq-6}, we obtain
\begin{equation}\label{sec2-3:eq-7}\begin{aligned}
\|\mathbf{c}-\mathbf{c}_{\lambda,\delta}\|
\leq\frac{1}{2\lambda\beta_{\min}}\sqrt{\|\mathbf{r}_{\lambda}\|^2-\|\mathbf{r}\|^2}
+\frac{\delta}{\delta+\hat{\sigma}^2_{\min}}
\|\mathbf{c}_{\lambda}-\mathbf{d}_{*}\|.
\end{aligned}\end{equation}
This bound shows that a relatively large value of $\delta$ can be
chosen to balance the accuracy and well conditioning of the LS
system \eqref{s31:LeastSquare-22} when $\mathbf{d}_{*}$ is close to
$\mathbf{c}_{\lambda}$.  Even if the matrix $\mathbf{A}_{\lambda}$ is
singular, highly accurate numerical solutions can still be obtained,
as we show in Sections \ref{sec2-4} and \ref{sec:numerical}.

\subsection{Examples}\label{sec2-4}
We present an example to show the effectiveness of \sjwchange{methods
based on the
formulations} \eqref{s31:LeastSquare2}, \eqref{s31:LeastSquare-22}, \eqref{s31:LeastSquare-33},
and \eqref{s31:LeastSquare-3} for solving a model problem on three
irregular domains, including a nonconvex domain.  \sjwchange{We compare QR
decomposition for solving the unconstrained LS formulation with
solution of the KT system \eqref{eq:kkt} that arises from the
constrained LS formulation. Our results show that both methods achieve
a similar level of accuracy, while the method based on QR
factorization is slightly more efficient.}
\begin{example}\label{eg2-1}
Solve \eqref{s31:eq-1} subject to suitable Dirichlet boundary
conditions and a source term $f(x,y)$, such that the solution
of \eqref{s31:eq-1} is
$$u(x,y)=\exp(x+y).$$
\end{example}
This problem is solved on the  domains defined by the following three cases.
\begin{itemize}[leftmargin=*]
\item Case I: the circular domain $\Omega=\{(x,y)|x^2+y^2<1\}$, see
  Figure~\ref{fig:iregular_domain}(b).
\item Case II: the irregular
  domain $\Omega$ defined as shown in
  Figure~\ref{fig:iregular_domain}(c) with the boundary defined by the
  B-spline interpolation using the Matlab function
  \texttt{spline}$(\mathbf{x,y})$, where \begin{equation*}
\begin{aligned}
&\mathbf{x} = \left[
              \begin{array}{cccccccccc}
               0  &\pi/2   &  \pi & 3\pi/2 &  2\pi
              \end{array}
            \right],\\
&\mathbf{y}=\left[
              \begin{array}{cccccccccc}
               0  &1   &  0 & -1 &  0  & 1  & 0   \\
              1.7 &0.9 & 1.8& 0.5& -0.7& 0.9& 1.7 \\
              \end{array}
            \right].
\end{aligned}\end{equation*}
\item Case III : the irregular domain $\Omega$ defined as shown in
  Figure \ref{fig:iregular_domain}(d) with the boundary defined by the
  Matlab function \texttt{spline}$(\mathbf{x,y})$, where
  \begin{equation*}
\begin{aligned}
&\mathbf{x} = \left[
              \begin{array}{cccccccccc}
               0  &2\pi/7   &  4\pi/7 & 6\pi/7 &  8\pi/7 &10\pi/7 &12\pi/7 &2\pi
              \end{array}
            \right],\\
&\mathbf{y}=\left[
              \begin{array}{cccccccccc}
                0   & 1   & 0.5 &      0& -0.5& -1 &-0.4 &      0 & 1 & 0   \\
                1.7 & 0.9 & 1.051 & 1.708&  0.791&  0.511 &-0.107 & 0.296 & 0.9 & 1.7 \\
              \end{array}
            \right].
\end{aligned}\end{equation*}
\item Case IV : The irregular domain with oscillatory boundaries defined by
the Matlab function \texttt{interp1q}$(\mathbf{x,y,x_i})$, where
\begin{equation*}
\begin{aligned}
&\mathbf{x} = \left[
              \begin{array}{cccccccccccccccc}
               1 &&-1 &&-1 &&1  &&1/2  &&1  &&0  &&1
              \end{array}
            \right],\\
&\mathbf{y}=\left[
              \begin{array}{cccccccccccccccc}
                1 &&1  &&-1 &&-1 &&-1/2 &&0&&   \Lambda &&1 \\
              \end{array}
            \right];
\end{aligned}\end{equation*}

see Figure~\ref{fig:iregular_domain-2}(a) for $\Lambda=1/2$ and
Figure~\ref{fig:iregular_domain-2}(b) for $\Lambda=3/4$.

\begin{figure}[!h]
\begin{center}
\begin{minipage}{0.47\textwidth}\centering
\epsfig{figure=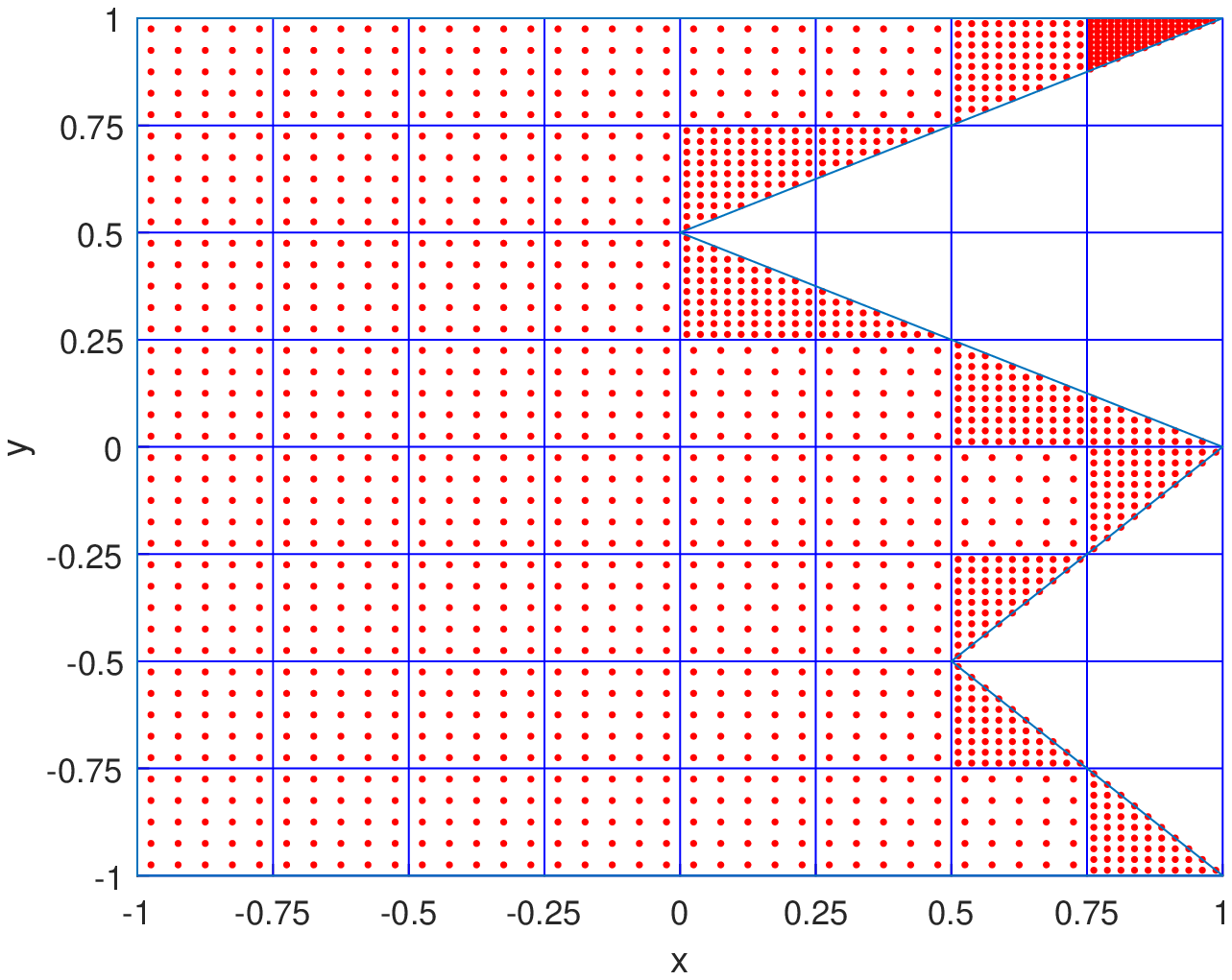,width=6cm} \par {(a) $\Lambda=1/2$.}
\end{minipage}
\begin{minipage}{0.47\textwidth}\centering
\epsfig{figure=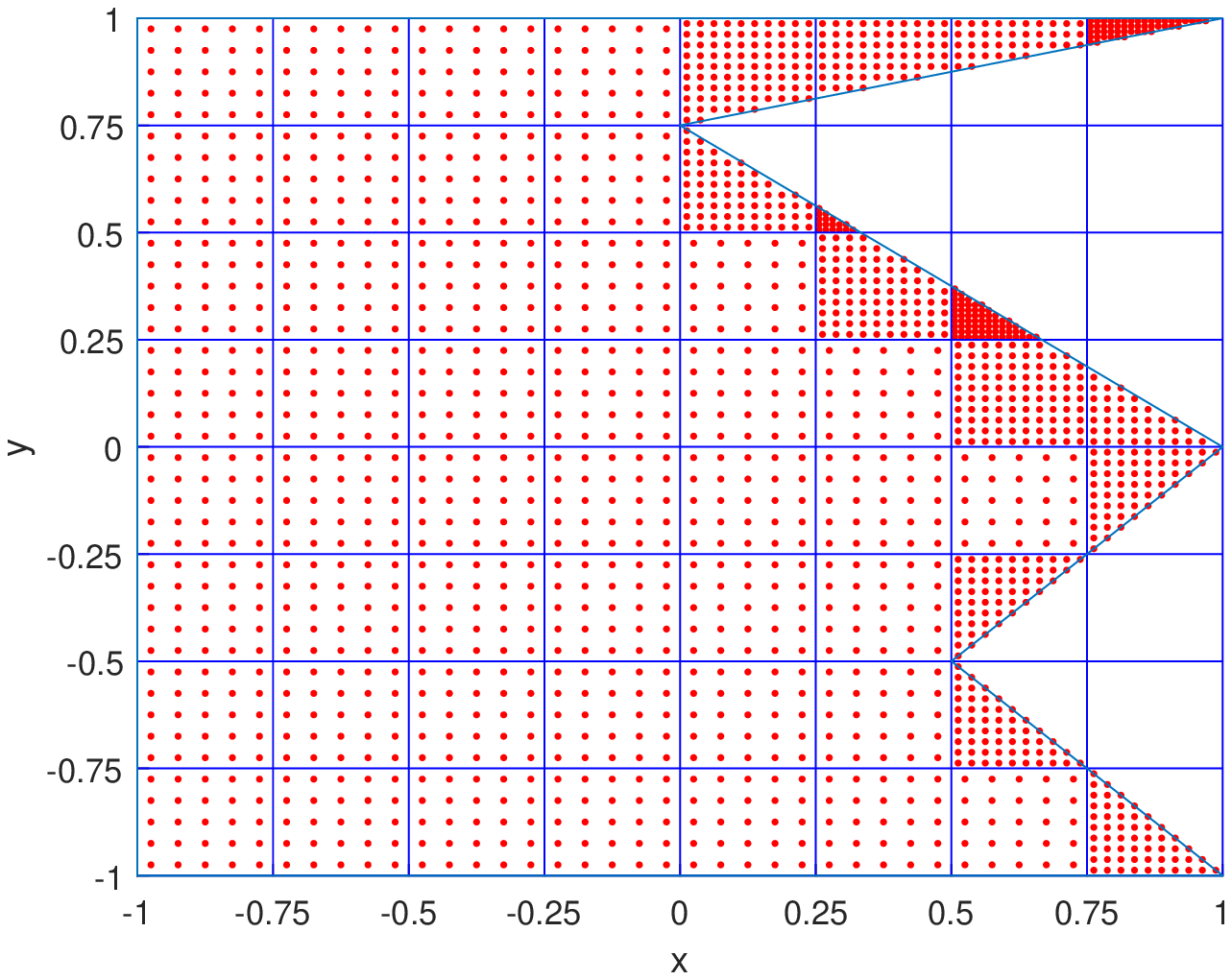,width=6cm} \par {(b)  $\Lambda=3/4$. }
\end{minipage}
\end{center}
\caption{The irregular domain with oscillatory boundaries.}\label{fig:iregular_domain-2}
\end{figure}

\end{itemize}


The $L^2$ error is measured by
$$\|e\|=\left(h_xh_y\sum_{i=0}^{200}\sum_{j=0}^{200}|e_{i,j}|^2\right)^{1/2},$$
where $e_{i,j}=u(\hat{x}_i,\hat{y}_j)-U_h(\hat{x}_i,\hat{y}_j)$ for
$(\hat{x}_i,\hat{y}_j)\in\bar{\Omega}$ and $e_{i,j}=0$ for
$(\hat{x}_i,\hat{y}_j)\notin\Omega$, $\hat{x}_i=a+i\hat{h}_x$,
$\hat{y}_j=-c+j\hat{h}_y$, $\hat{h}_x=(b-a)/200$, and
$\hat{h}_y=(d-c)/200$.

The default solver for the unconstrained LS
problems \eqref{s31:LeastSquare-22} and \eqref{s31:LeastSquare-3} is
the QR decomposition, \sjwchange{implemented via the backslash
operator ``$\setminus$'' in Matlab 2018a}) with $\lambda=10^5$.  We
take $\nu=0.1$ and $N_x=N_y=N$ in the numerical simulations. The
nonuniformly distributed collocation point set
$\mathcal{S}^{(10,10)}_{\Omega}$ is used in the four
methods \eqref{s31:LeastSquare2}, \eqref{s31:LeastSquare-33},
\eqref{s31:LeastSquare-22}, and \eqref{s31:LeastSquare-3}, \sjwchange{unless
specified otherwise.}

 The choice of $\mathbf{d_*}$ in the
formulations \eqref{s31:LeastSquare2}, \eqref{s31:LeastSquare-22}, \eqref{s31:LeastSquare-33},
and \eqref{s31:LeastSquare-3} is as follows. { First, we let
$\delta=\delta_0=10^{-6}$ and $\mathbf{d}_*=\mathbf{0}$ in
\eqref{s31:LeastSquare-22} and solve it by   QR factorization. The solution is denoted as $\mathbf{d}^*$.}
Then, we take
$\mathbf{d_*}=\mathbf{d}^*\cdot(1+\epsilon \mathbf{\tilde{d}})$, where
$\epsilon$ is a nonnegative number and $\mathbf{\tilde{d}}$ is a
random vector generated by Matlab function \texttt{rand}.
(Call \texttt{rng}$(1)$ before calling
$\mathbf{\tilde{d}}=\texttt{rand}(M,1)$
\sjwchange{to derive the same value of $\mathbf{\tilde{d}}$ each time the code is run}.)

\begin{figure}[!h]
\begin{center}
\begin{minipage}{0.47\textwidth}\centering
\epsfig{figure=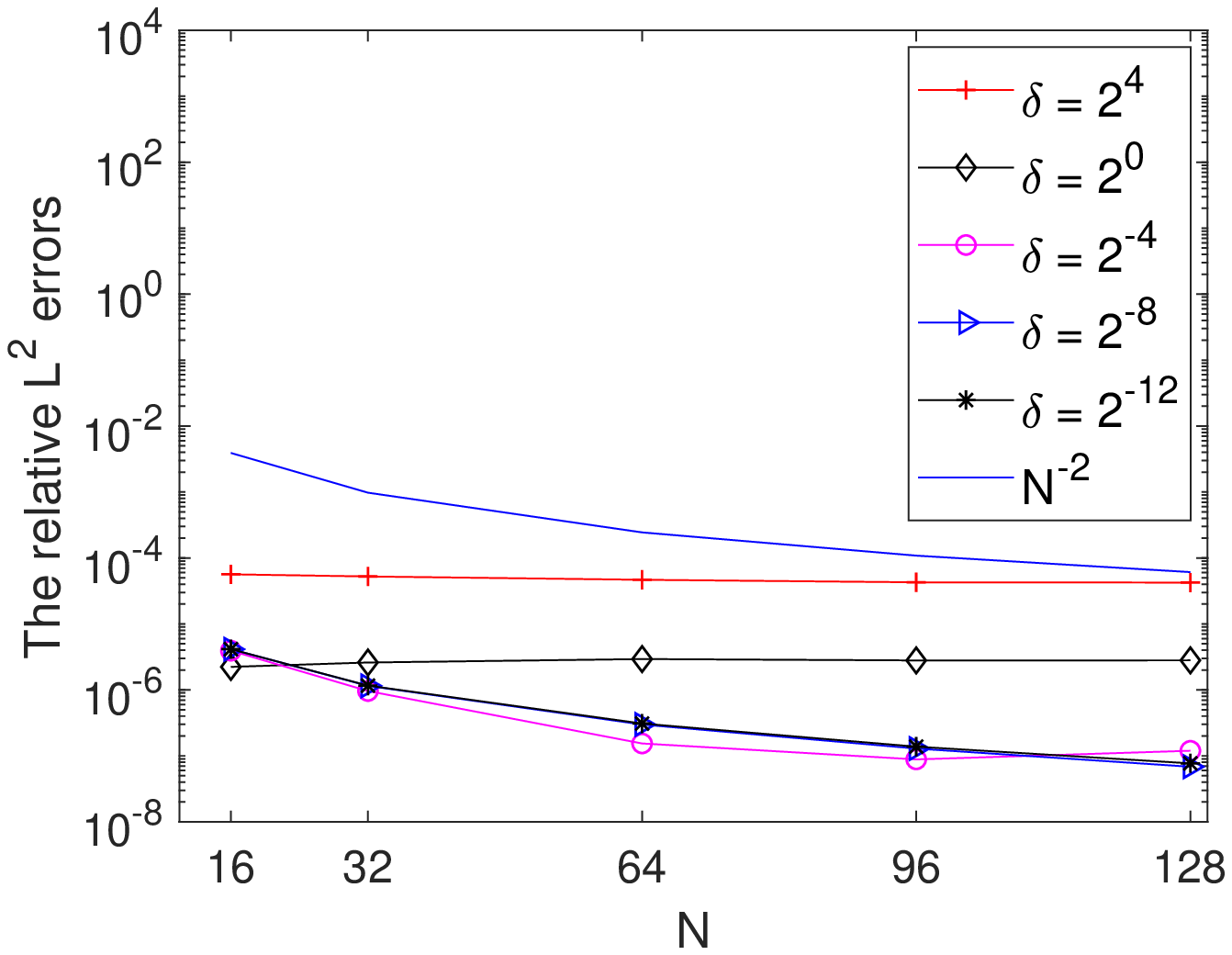,width=7cm} \par {(a) $\epsilon=0.01$ }
\end{minipage}
\begin{minipage}{0.47\textwidth}\centering
\epsfig{figure=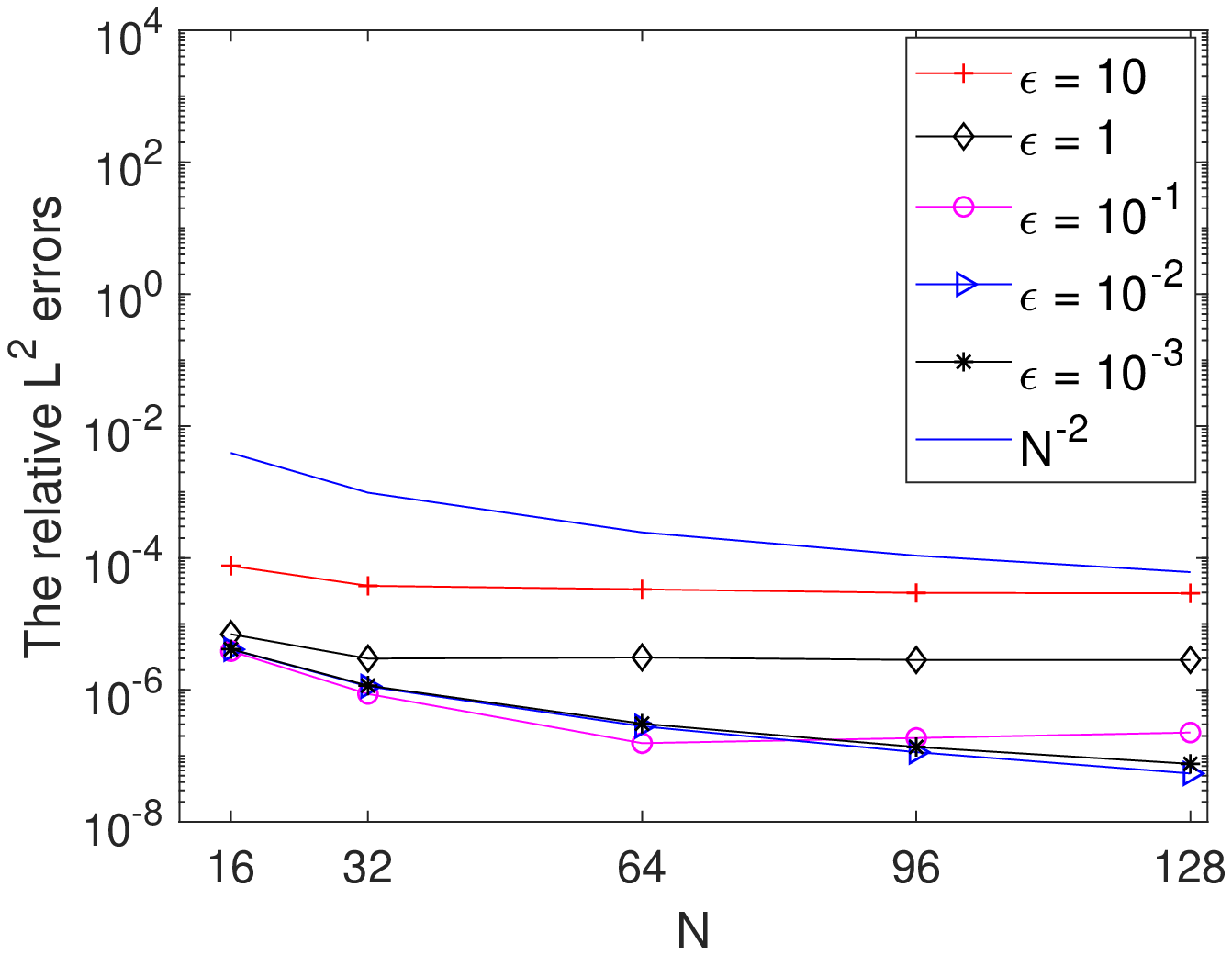,width=7cm} \par {(b) $\delta = 0.01$ }
\end{minipage}
\end{center}
\caption{The relative $L^2$ error $\|e\|/\|u\|$ of
 the LSC method \eqref{s31:LeastSquare-22},  Example \ref{eg2-1}, Case I.\label{fig2-1}}
\end{figure}

We show first how $\delta$ and $\epsilon$ affect the accuracy of the
LSC method \eqref{s31:LeastSquare-22} for Case I. In
Figure~\ref{fig2-1}(a), we fix $\epsilon=0.01$ and see that the
accuracy increases as $\delta$ decreases.  Second-order accuracy is
observed for smaller $\delta=2^{-4},2^{-8}$, and $2^{-12}$. In
Figure~\ref{fig2-1}(b), we fix $\delta=0.01$ and see that the error
decreases as $\epsilon$ decreases. Second-order accuracy is observed
for $\epsilon=10^{-2}$ and $10^{-3}$.

\begin{figure}[!h]
\begin{center}
\begin{minipage}{0.47\textwidth}\centering
\epsfig{figure=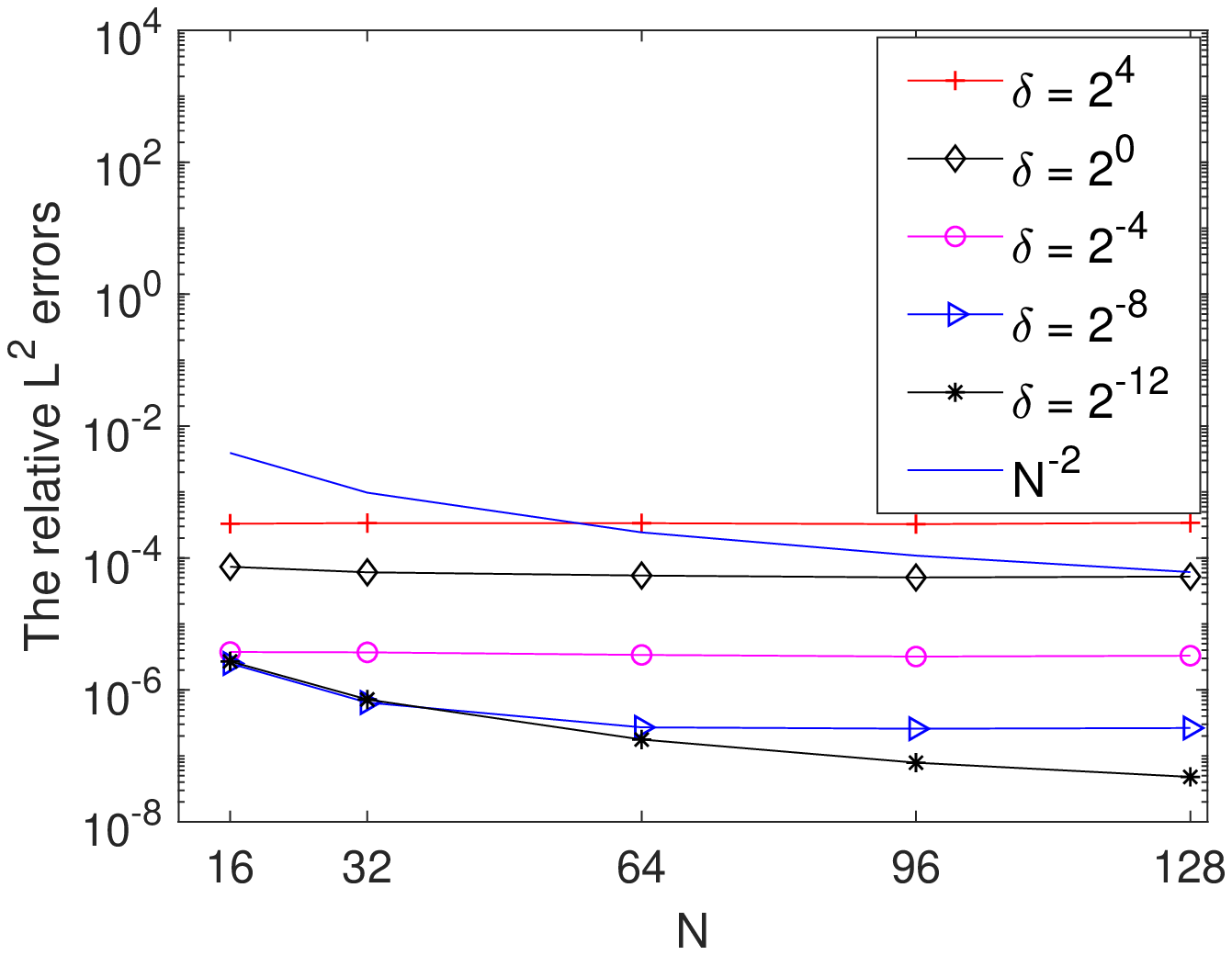,width=7cm} \par {(a) $\epsilon=0.01$ }
\end{minipage}
\begin{minipage}{0.47\textwidth}\centering
\epsfig{figure=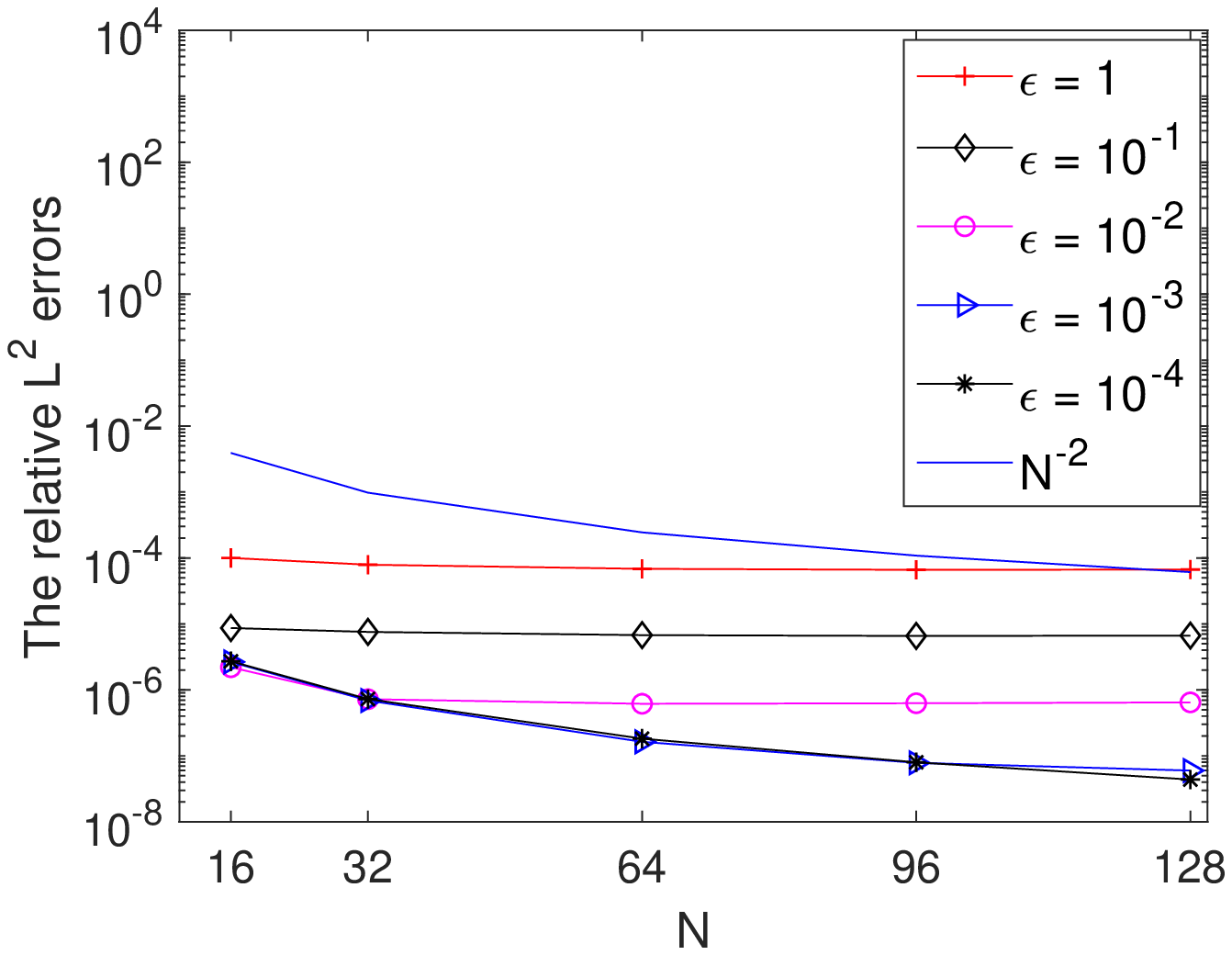,width=7cm} \par {(b) $\delta = 0.01$ }
\end{minipage}
\end{center}
\caption{The relative $L^2$ error $\|e\|/\|u\|$ of the LSC method \eqref{s31:LeastSquare-22},
Example \ref{eg2-1},
Case II.\label{fig2-2}}
\end{figure}

Figure~\ref{fig2-2} shows the relative $L^2$
errors \sjwchange{obtained from the LSC
formulation \eqref{s31:LeastSquare-22}} for Case II. We can see that
decreasing $\delta$ or $\epsilon$ improves the accuracy, and
second-order is observed for \sjwchange{smaller values of these
parameters}.

\begin{figure}[!h]
\begin{center}
\begin{minipage}{0.47\textwidth}\centering
\epsfig{figure=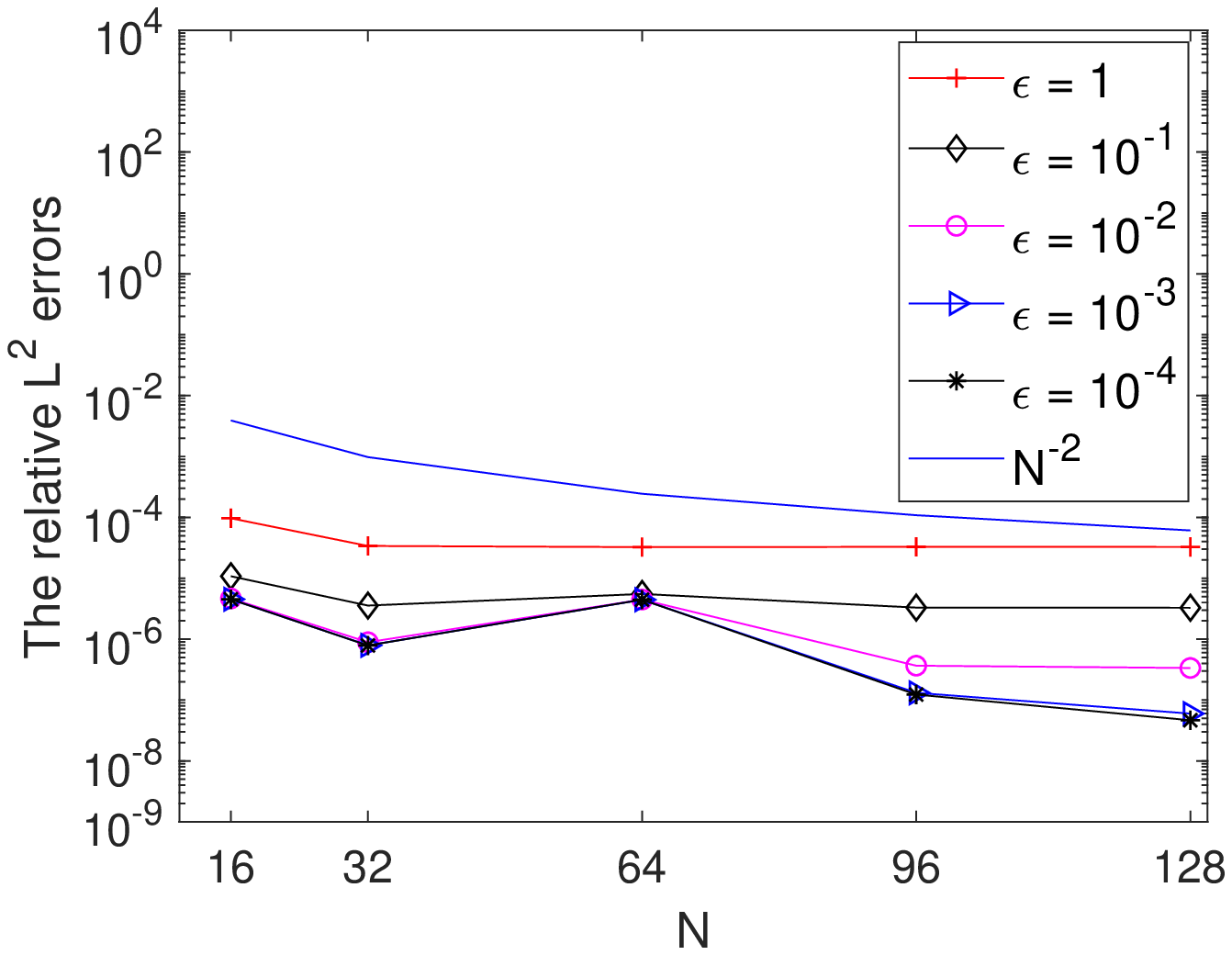,width=6.8cm} \par {(a1)   $p=q=10$. }
\end{minipage}
\begin{minipage}{0.47\textwidth}\centering
\epsfig{figure=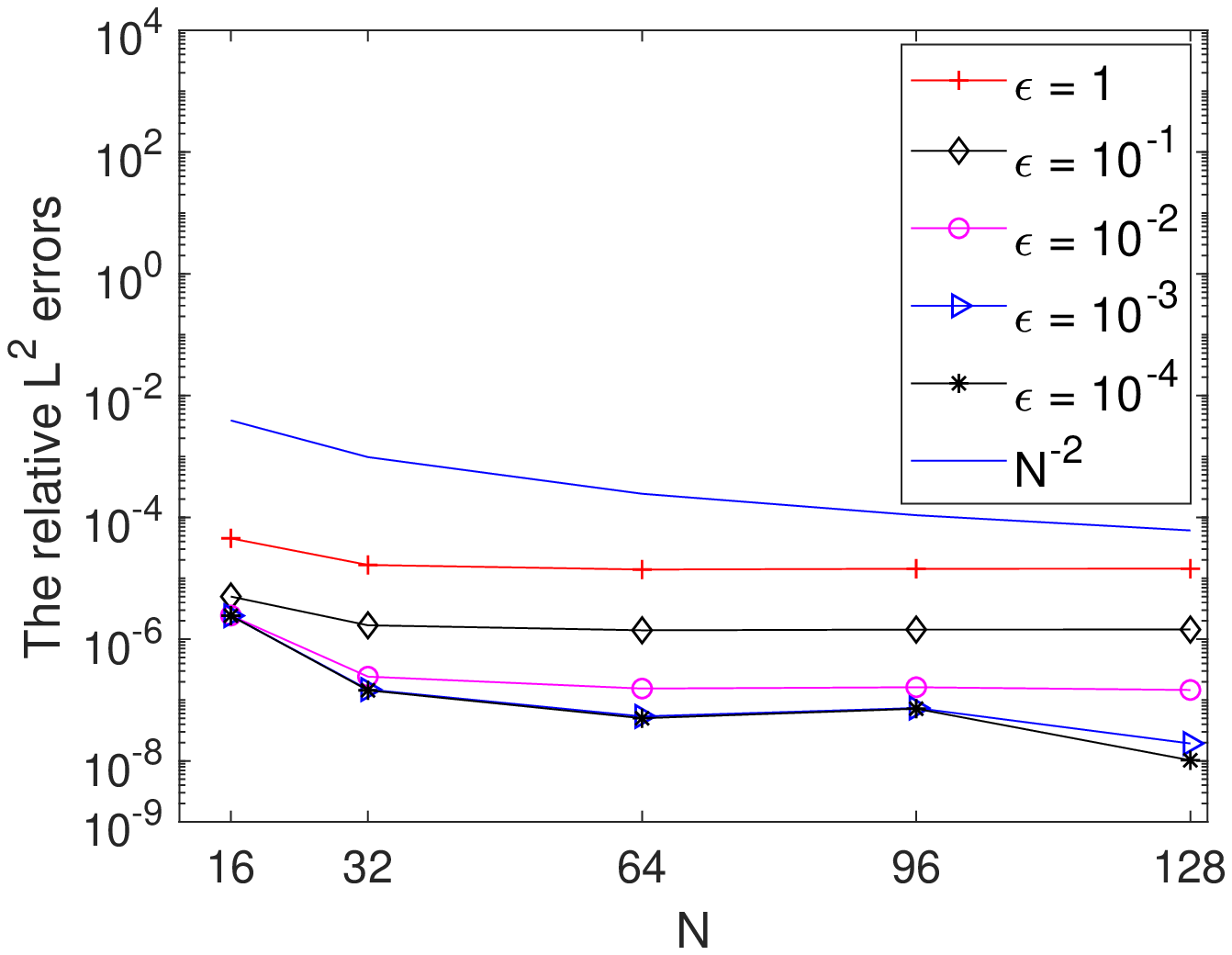,width=6.8cm} \par {(a2)  $p=q=16$. }
\end{minipage}
\begin{minipage}{0.47\textwidth}\centering
\epsfig{figure=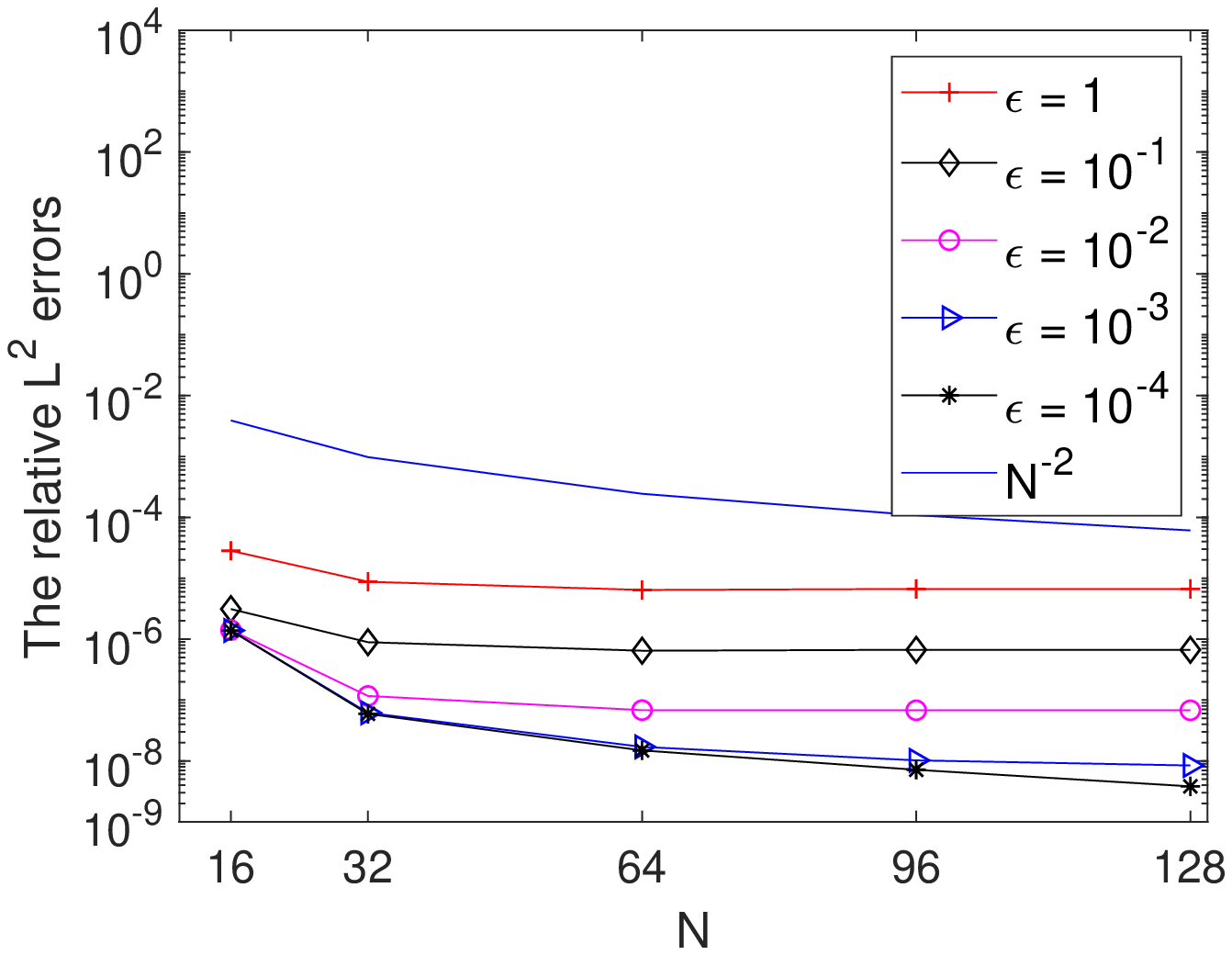,width=6.8cm} \par {(a3) $p=q=24$. }
\end{minipage}
\begin{minipage}{0.47\textwidth}\centering
\epsfig{figure=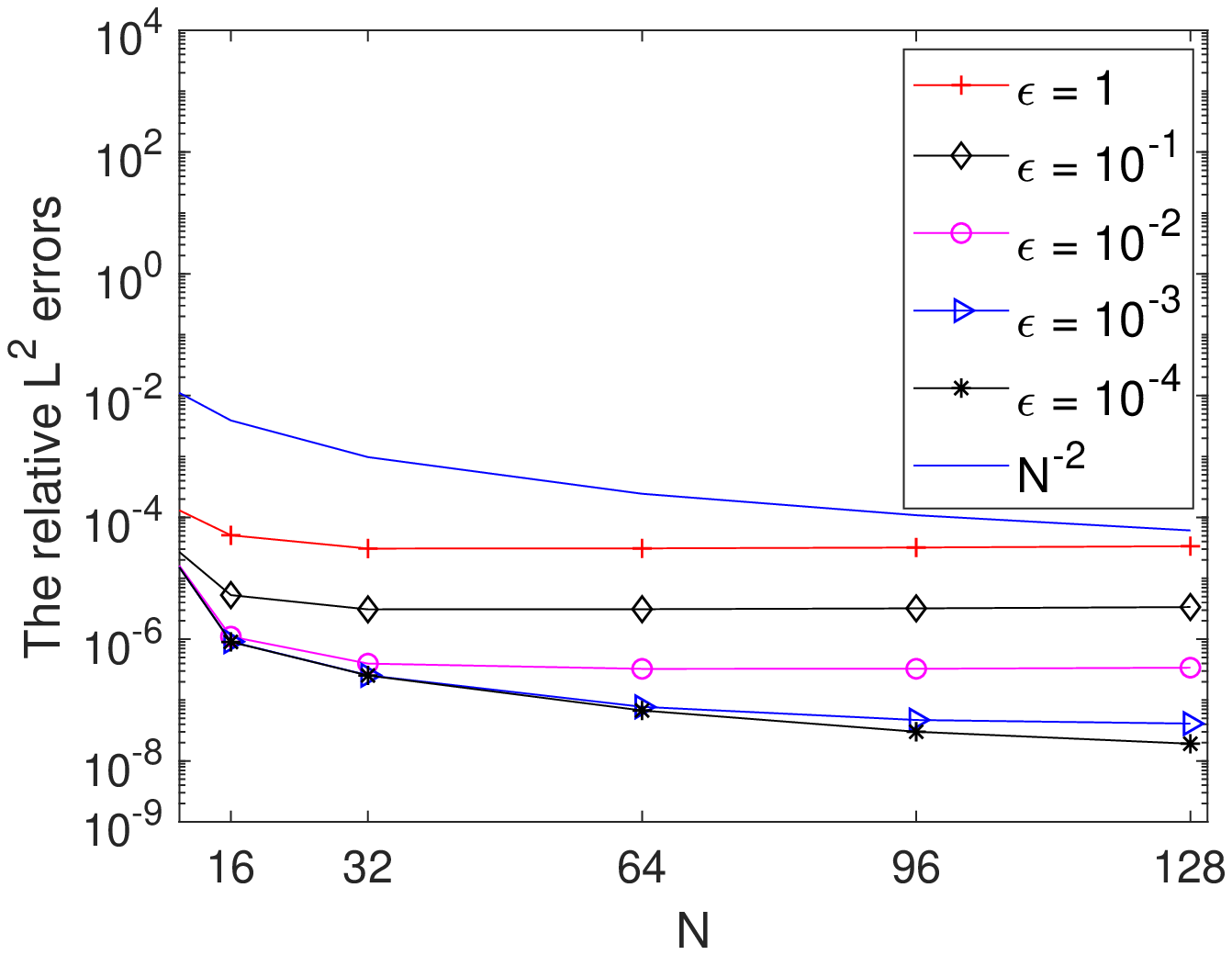,width=6.8cm}
\par {(b)   }
\end{minipage}
\end{center}
\caption{The relative $L^2$ errors $\|e\|/\|u\|$ of the LSC method
\eqref{s31:LeastSquare-22}, Example \ref{eg2-1},  Case III, $\delta=0.01$.
(a1), (a2), and (a3) :  the LS method
based on the uniformly distributed collocation point set $\mathcal{\overline{S}}^{(p,q)}_{\Omega}$; (b)
 the LS method \eqref{s31:LeastSquare-22}
 based on non-uniformly distributed collocation point set
 $\mathcal{S}^{(10,10)}_{\Omega}$.\label{fig2-2-3}}
\end{figure}

\sjwchange{On the nonconvex domain of Case III, we still}
observe similar results as for Cases I and II; see
Figure~\ref{fig2-2-3}(b).  In Figures~\ref{fig2-2-3}(a1)-(a3), we also
show the errors of the LSC method \eqref{s31:LeastSquare-22} based on
uniformly distributed collocation points of different densities.
\sjwchange{The LSC method works less well when there are
insufficient} collocation points in the boundary cells; see
Figure~\ref{fig2-2-3}(a1).  As the number of collocation points in the
boundary cells increases, better results are obtained; see
Figures~\ref{fig2-2-3}(a2)-(a3).
\sjwchange{We conclude from these results that}
sufficient collocation points in the boundary cells are needed to
achieve second-order accuracy, and that adding more collocation points
to the inner cells of $\Omega$ does not necessarily improve the
accuracy of the method.  In the following, we will use the LSC method
or LSFVM based on nonuniformly distributed collocation points to solve
the considered differential equations.

\begin{figure}[!h]
\begin{center}
\begin{minipage}{0.47\textwidth}\centering
\epsfig{figure=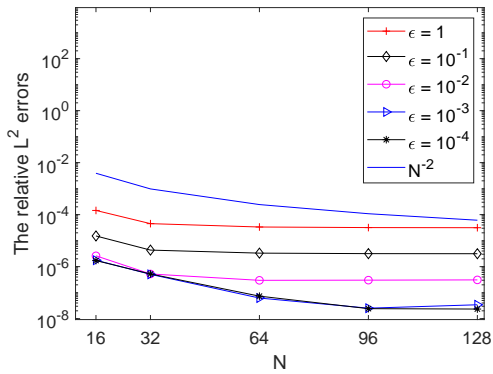,width=7cm} \par {(a) $\Lambda=1/2$. }
\end{minipage}
\begin{minipage}{0.47\textwidth}\centering
\epsfig{figure=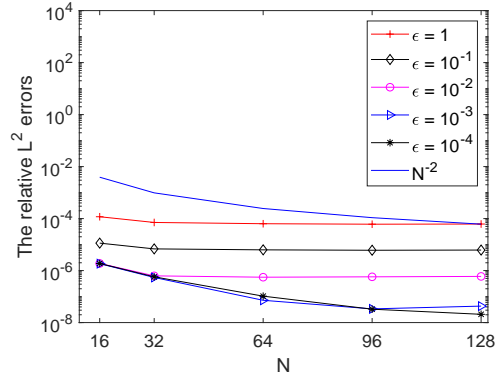,width=7cm} \par {(b)  $\Lambda=3/4$ }
\end{minipage}
\end{center}
\caption{The relative $L^2$ errors $\|e\|/\|u\|$ of the LSC method
\eqref{s31:LeastSquare-22} based on non-uniformly distributed collocation point set
 $\mathcal{S}^{(10,10)}_{\Omega}$, Example \ref{eg2-1},  Case IV, $\delta=0.01$.\label{fig2-2-4}}
\end{figure}

{  For Case IV, we solve the problem on the domain with
oscillatory boundaries.  Figure~\ref{fig2-2-4} shows that decreasing
$\epsilon$ increases accuracy, and second-order accuracy is observed
for $\epsilon=10^{-4}$, which is similar to that exhibited in
Figure \ref{fig2-2} (b).  }

\sjwchange{Table~\ref{tb2-1} shows the relative $L^2$ errors and the
corresponding computational time of two approaches for Case III: QR
decomposition applied to the penalized linear-least-squares
formulation \eqref{s31:LeastSquare-22}, and direct solution of the KKT
conditions \eqref{eq:kkt} for the constrained
formulation \eqref{s31:LeastSquare2}.  Both solvers achieve results of
identical accuracy, but the method based on QR decomposition appears
to be slightly more efficient, especially when the matrix is large.}



Table~\ref{tb2-3} compares the LSC method \eqref{s31:LeastSquare-22}
with the LSFVM \eqref{s31:LeastSquare-3}. \sjwchange{The two
approaches attain similar levels} of accuracy for $\rho=10^{-4}$,
which is in line with the theoretical analysis in
Section~\ref{sec2-2-2}. The results in Table~\ref{tb2-3}
also \sjwchange{suggest} that the LSFVM reduces to the LSC method when
$\rho\to 0$.

{ The constrained
formulation \eqref{s31:LeastSquare2} (or \eqref{s31:LeastSquare-33}) yielded similar results  as those from \eqref{s31:LeastSquare-22} (or \eqref{s31:LeastSquare-3}) for $\lambda=10^5$, these results are not shown here.}



\begin{table}[!ht]
\caption{{The relative errors $\|e\|/\|u\|$   and computation times for a QR-based solver applied to the formulation \eqref{s31:LeastSquare-22} and direct solution of the KKT system \eqref{eq:kkt} on  Example \ref{eg2-1},
Case III, $\delta=\epsilon=10^{-4}$. }}\label{tb2-1}
\centering\footnotesize
\begin{tabular}{|c|cc|cc|cc|c|c|c|c|c|c|}
\hline
 Method & $N=32$& Time (s)& $N=64$  & Time (s)& $N=128$& Time (s)\\
 \hline
QR,  \eqref{s31:LeastSquare-22}  &2.5601e-7 & 1.5560e-2 & 6.8002e-8 & 5.0821e-2 & 1.6083e-8 & 2.0375e-1 \\
KKT, \eqref{eq:kkt} &2.5601e-7&2.3683e-2&6.8002e-8&9.2296e-2&1.6083e-8&4.5509e-1 \\
\hline\end{tabular}
\end{table}

\begin{table}[!ht]
\caption{{Comparison of the relative $L^2$ \sjwchange{difference} $\|e\|/\|u\|$  between
 the LCM method \eqref{s31:LeastSquare-22} and the LSFVM \eqref{s31:LeastSquare-3}, Example \ref{eg2-1},
 Case II, $\delta=0.01$, and $N=64$.}}\label{tb2-3}
\centering\footnotesize
\begin{tabular}{|c|c|c|c|c|c|c|c|c|c|c|c|c|}
\hline
 & & $\epsilon=1$& $\epsilon=10^{-2}$ &  $\epsilon=10^{-4}$  \\
 \hline
LSC  &$-$              &6.8916e-5&6.2460e-7&1.8022e-7\\    \hline
&$\rho=10^{-4}$        &6.8916e-5&6.2462e-7&1.8011e-7\\
&$\rho=10^{-3}$        &6.8879e-5&6.2656e-7&1.6937e-7\\
LSFVM&$\rho=10^{-2}$   &6.1956e-5&1.3875e-6&9.0710e-7\\
&$\rho=2\times10^{-2}$ &7.3184e-5&4.5220e-6&4.1602e-6\\
&$\rho=10^{-1}$        &1.5751e-4&1.0895e-4&1.0873e-4\\
\hline
\end{tabular}
\end{table}

\sjwchange{One further example shows} that the Neumann boundary conditions can be handled easily in the present framework.
{
\begin{example}\label{eg2-2}
Consider \eqref{s31:eq-1} subject to the following Neumann boundary
conditions
\begin{equation}\label{s24:eq-NBC}
\mathcal{B}(u)=(\mathbf{\hat{n}}\cdot \nabla u)(x,y)=u_b(x,y),
\quad (x,y)\in\partial \Omega,
\end{equation}
where $\mathbf{\hat{n}}=(\hat{n}_x,\hat{n}_y)^T$  denotes the unit normal
to the outer boundary of $\Omega$.
\end{example}
}
\sjwchange{In a similar fashion to} \eqref{s31:eq-10}, \eqref{s31:eq-11},  we  obtain the following overdetermined system:
\begin{equation}\label{s24:eq-1}\left\{\begin{aligned}
&\left(\Phi^T\mathbf{d}
- \nu\Delta \Phi^T\mathbf{d} \right)(\xi^{(in)}_i) -
f(\xi^{(in)}_i)=0,{\quad}1\leq i \leq N_{in},\\ &\left(\Phi^T\mathbf{d}
- \nu\Delta \Phi^T\mathbf{d} \right)(\xi^{(b)}_i)-f(\xi^{(b)}_i)
-\widetilde{\lambda} \left(\mathcal{B}(\Phi^T\mathbf{d})-u_b\right) (\xi^{(b)}_i)
=0 , {\quad}1\leq i \leq N_{b},
\end{aligned}\right.\end{equation}
where the Neumann boundary conditions (the second equation
in \eqref{s24:eq-1}) are imposed according to
\cite[equation (28)]{LorenzisEHR15}. \sjwchange{The parameter $\widetilde{\lambda}$} is a positive constant
dependent on the mesh, which is chosen as
$\widetilde{\lambda}=4\max\{(b-a)N_x,(d-c)N_y\}$ in the numerical simulations in
this work; see \cite[Section~4.1.4]{LorenzisEHR15}.

The penalized LSC method   for \eqref{s31:eq-1} subject to the
Neumann boundary condition \eqref{s24:eq-NBC} is given by
\begin{equation}\label{s31:LSC-nbc}\begin{aligned}
\mathbf{c}&=\arg\min_{\mathbf{d}\in  \mathbb{R}^{M} }
\left\|\left(
           \begin{array}{c}
             \mathbf{A}_{in}-\nu\mathbf{S}_{in} \\
             \mathbf{A}_{b}-\nu\mathbf{S}_{b} -\widetilde{\lambda} \mathbf{B} \\
             \sqrt{\delta}\mathbf{I}\\
           \end{array}
         \right)\mathbf{d} -\left(
           \begin{array}{c}
             \mathbf{f}_{in} \\
             \mathbf{f}_{b}-\widetilde{\lambda}\mathbf{u}_b  \\
             \sqrt{\delta}\mathbf{d_*}\\
           \end{array}
         \right)\right\|^2,
\end{aligned}\end{equation}
where  $\mathbf{B}$ is related to $(\mathbf{\hat{n}}\cdot \nabla U_h)(x,y)$
(We avoid an explicit specification here).

\begin{figure}[!h]
\begin{center}
\begin{minipage}{0.47\textwidth}\centering
\epsfig{figure=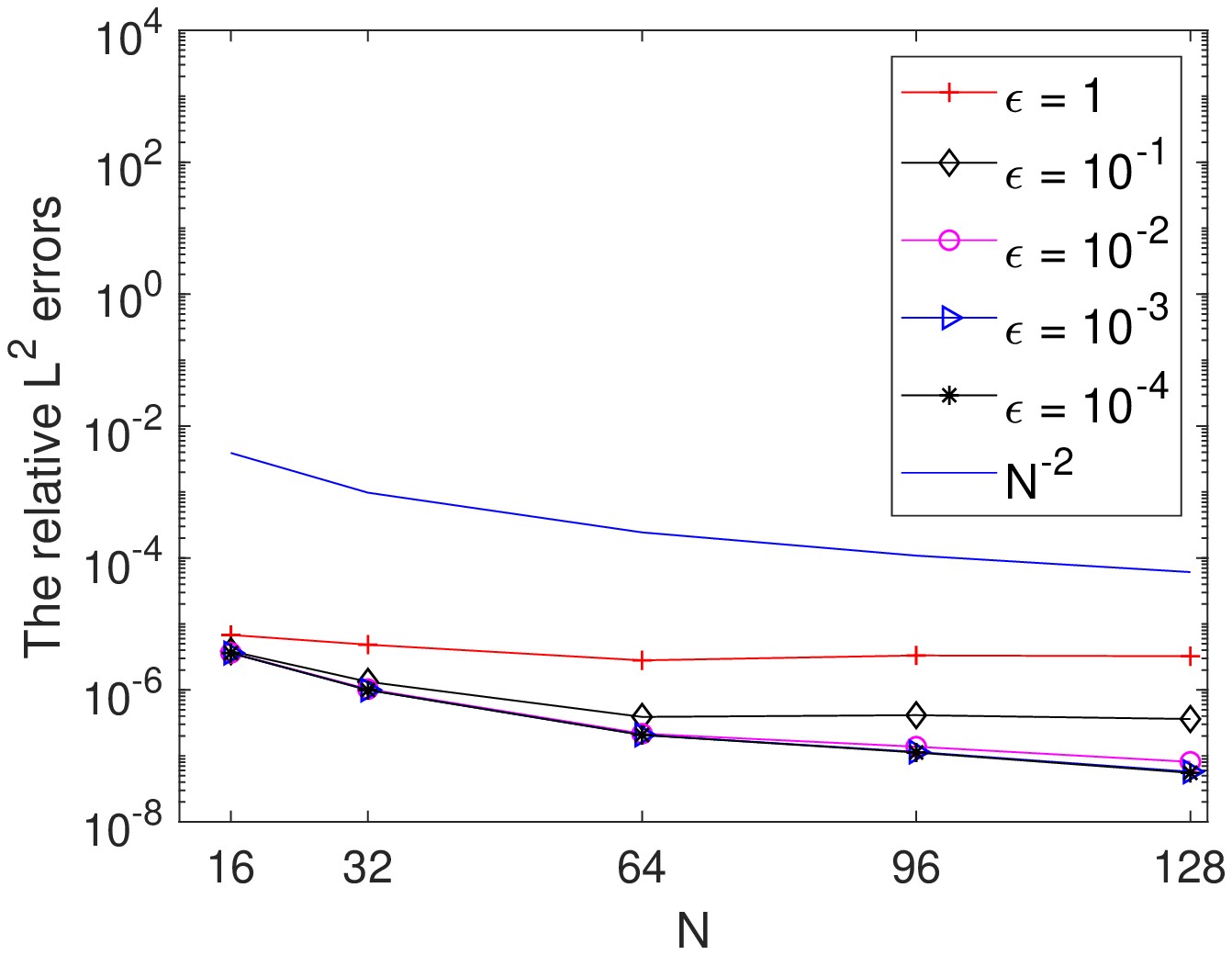,width=7cm} \par {(a) The circular domain in Case I }
\end{minipage}
\begin{minipage}{0.47\textwidth}\centering
\epsfig{figure=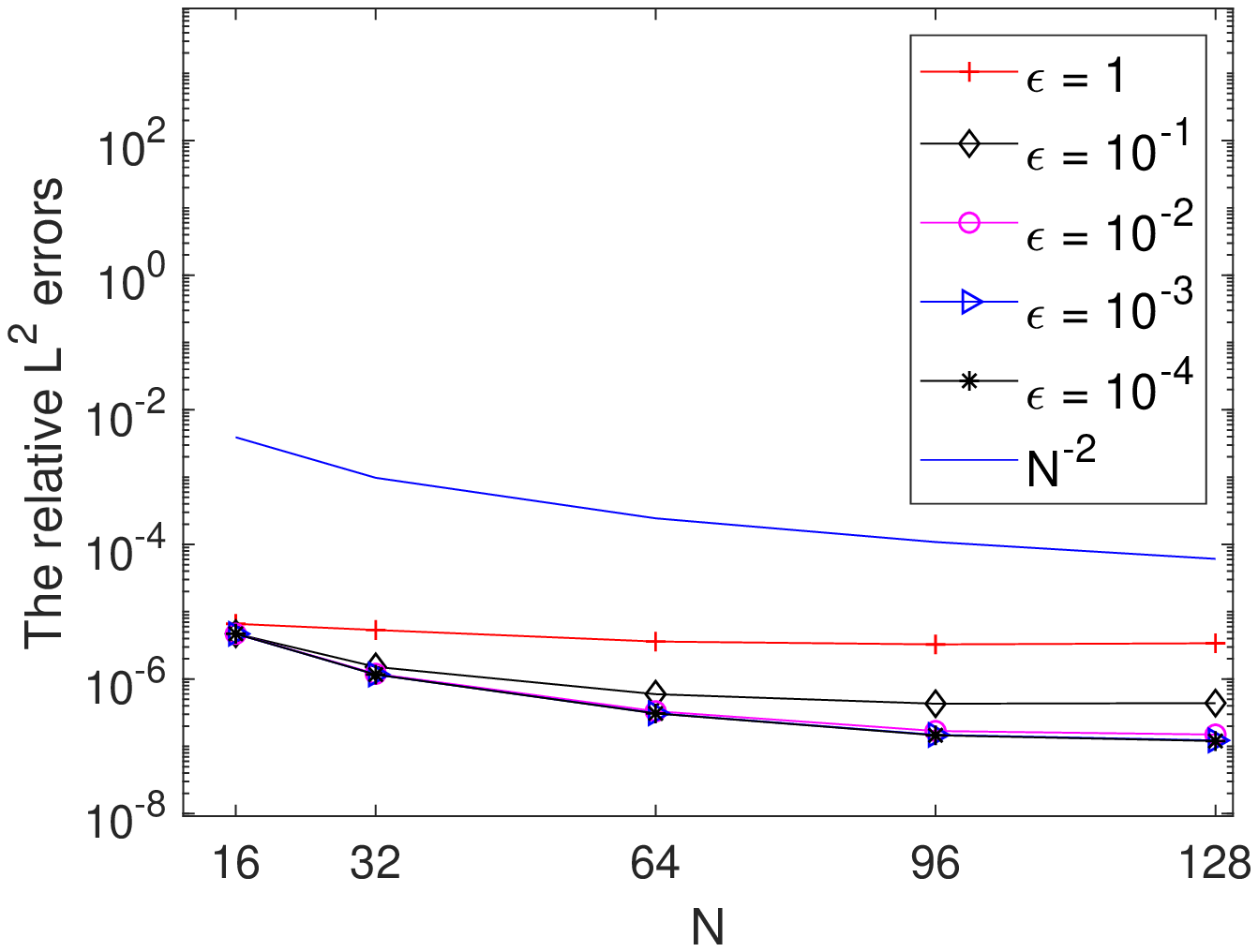,width=7cm} \par {(b) The irregular domain in Case II }
\end{minipage}
\end{center}
\caption{Relative errors  of the LSC method \eqref{s31:LSC-nbc} subject to the Neumann boundary
conditions \eqref{s24:eq-NBC}, compared to the the exact solution $u=\exp(x+y)$, Example~\ref{eg2-2}, $\nu=0.1$, $\delta=0.0001$, $\widetilde{\lambda}=4\max\{(b-a)N,(d-c)N\}$.\label{fig2-5}}
\end{figure}

We apply \eqref{s31:LSC-nbc} to solve \eqref{s31:eq-1} subject to the
Neumann boundary condition \eqref{s24:eq-NBC} on the circular domain
(see Figure~\ref{fig:iregular_domain}(b) and Case I in
Example~\ref{eg2-1}) and the irregular domain (see
Figure~\ref{fig:iregular_domain}(c) and Case II in
Example~\ref{eg2-1}).  Relative errors are shown in
Figure~\ref{fig2-5}. Note that second-order accuracy is still observed
\sjwchange{for small $\epsilon$.}

\section{Application to nonlinear time-dependent models}\label{sec3-fpde}
In this section, we extend the discrete LSC method and LSFVM
to solve the following nonlinear time-fractional diffusion equation
\begin{equation}\label{s5:eq-1}\left\{\begin{aligned}
&{}_{C}D^{\alpha}_{0,t}u(x,y,t) =  \nu \Delta u(x,y,t) + f(u,x,y,t),
{\quad}&&(x,y,t)\in\Omega\times(0, T],\\
&u(x,y,0)=\phi(x,y),{\quad}&&(x,y)\in\Omega,\\
&u(x,y,t)=u_b(x,y,t),{\quad}&&(x,y,t)\in\px[]\Omega\times[0,T].
\end{aligned}\right.\end{equation}
\sjwchange{We apply} a stable semi-implicit time-stepping method
to discretize the time direction of \eqref{s5:eq-1} to derive a linear
boundary value problem. The boundary value problem is solved by the
\sjwchange{LSC or  LSFVM approaches} developed in the previous section.

\subsection{Semi-implicit time-stepping method}\label{sec-3-1}

We recall a semi-implicit method proposed in \cite{ZengTBK2018} for
the fractional initial value problem
\begin{equation}\label{s4:fiv}
{}_{C}D^{\alpha}_{0,t}u(t) = \nu u(t) +  f(u(t),t),{\quad}u(0)=u^0,{\quad}t\in(0, T],
\end{equation}
where $0<\alpha\leq 1$, $\nu<0$, $f(u,t)$ is a nonlinear function with
respect to $u$, and the Caputo fractional operator
${}_{C}D^{\alpha}_{0,t}$ is defined by
\begin{equation}\label{caputo}
{}_{C}D^{\alpha}_{0,t}u(t) =\frac{1}{\Gamma(1-\alpha)}\int_0^t(t-s)^{-\alpha}u'(s)\dx[s].
\end{equation}
We assume that the solution $u(t)$ to \eqref{s4:fiv} satisfies
\begin{equation}\label{s4:solu}
u(t)-u(0) = \sum_{k=1}^{m}c_kt^{\gamma_k}
+ t^{\gamma_{m+1}}\tilde{u}(t),{\quad}0<\gamma_k<\gamma_{k+1},
\end{equation}
where $\tilde{u}(t)$ is uniformly bounded for $t\in[0,T]$.  This
assumption holds in real applications; see for example,
\cite{Diethelm-B10,FordMorReb13,Luchko11,Pod-B99,BurrageBTZ18}.
If $f(u(t),t)$ is smooth for $t\in[0,T]$, then
$\gamma_k\in\{i+j\alpha,i=0,1,\dotsc;j=0,1,\dotsc\}$;
see \cite{Diethelm-B10}.

Combining \eqref{s4:fiv} and \eqref{s4:solu}, we obtain
\begin{equation}\label{s4:fut}
f(u(t),t) -f(u(0),0)= {}_{C}D^{\alpha}_{0,t}u(t)- \nu u(t)
= \sum_{k=1}^{m_f}\hat{c}_kt^{\hat{\gamma}_k} +
t^{\hat{\gamma}_{k+1}}v(t),
\end{equation}
where $0<\hat{\gamma}_k<\hat{\gamma}_{k+1}$,
$\hat{\gamma}_k\in\{\gamma_k\}\cup\{\gamma_k-\alpha\}$, and $v(t)$ is
bounded for $t\in[0,T]$.

Denote by $t_j=j\tau\,(j\geq 0)$ the grid points, where $\tau=T/n_T$
is the stepsize and $n_T$ is a positive integer. Let $u^n=u(t_n)$ and
denote
\begin{equation}\begin{aligned}\label{s2:Dalf}
D_{\tau}^{(\alpha,n,m,\gamma)}u
=\frac{1}{\tau^{\alpha}}\sum_{j=0}^{n}\omega^{(\alpha)}_{n-j}(u^j-u^0)
+\frac{1}{\tau^{\alpha}}\sum_{j=1}^{m}w^{(\alpha)}_{n,j}(u^j-u^0),
\end{aligned}\end{equation}
where $m\geq 0$ is the number of correction terms and the  quadrature weights $\omega^{(\alpha)}_{j}$
satisfy
\begin{equation}
\omega^{(\alpha)}(z)=(1-z)^{\alpha}\left(1+\frac{\alpha}{2}-\frac{\alpha}{2}z\right)
=\sum_{n=0}^{\infty}\omega^{(\alpha)}_{n}z^n.
\end{equation}
 The starting weights $\{w^{(\alpha)}_{n,j}\}$ are chosen such that
$$\sum_{j=0}^{n}\omega^{(\alpha)}_{n-j}u^j
+ \sum_{j=1}^{m}w^{(\alpha)}_{n,j}u^j
=\frac{\Gamma(\gamma_k+1)}{\Gamma(\gamma_k+1-\alpha)}n^{\gamma_k-\alpha}$$
for some $u(t)=t^{\gamma_k}\,(0<\gamma_{k}<\gamma_{k+1})$.  (We refer
readers to \cite{DieFord06,Lub86,ZengZK17} for \sjwchange{further
information on calculation and properties of the starting weights.})

We also introduce the notation $E^{(n,m,\gamma)}_2$ defined by
\begin{equation}\label{s4:Ep-m}
E^{(n,m,\sigma)}_2(u)=u^n-2u^{n-1}
+u^{n-2}- \sum_{j=1}^{m}w^{(u)}_{n,j}(u^j - u^0),
\end{equation}
where  $\{w^{(u)}_{n,j}\}$ are chosen such that
$E^{(n,m,\gamma)}_2(u)=0$ for $u=t^{\gamma_r}$ $(0<\gamma_{r-1}<\gamma_r)$.


By applying the second-order semi-implicit time-stepping method
proposed in \sjwchange{\cite[equation~(17)]{ZengTBK2018}} to
discretize \eqref{s4:fiv}, we obtain
\begin{equation}\label{s4:eq-1}
D_{\tau}^{(\alpha,n,m,\gamma)}u
= \nu u^n +  f^n - E^{(n,m_f,\hat{\gamma})}_2(f)-\kappa  E^{(n,m_u,\gamma)}_2(u) + R^n,
\end{equation}
where $f^n = f(u^n,t_n)$, $\kappa$ is a constant that balances the
stability and accuracy of method \eqref{s4:eq-1}, and $R^n$ is the
time discretization error satisfying
\begin{equation}\begin{aligned}\label{s4:Rn}
R^n=O(\tau^2t_n^{\gamma_{m+1}-2-\alpha})
+O(\tau^2t_n^{\hat{\gamma}_{m_f+1}-2})+O(\tau^2t_n^{\gamma_{m_u+1}-2}).
\end{aligned}\end{equation}

{
\begin{remark}
Generally speaking, the solution of \eqref{s4:fiv} has singularity at
$t=0$ due to the singularity of the kernel in the fractional
derivative operator \eqref{caputo}, which may lead to low accuracy of
some numerical methods \cite{Diethelm-B10}.  In this work, we use the correction
${\tau^{-\alpha}}\sum_{j=1}^{m}w^{(\alpha)}_{n,j}(u^j-u^0)$ to
preserve high accuracy of the time-stepping method \eqref{s2:Dalf} for
the approximation of ${}_{C}D^{\alpha}_{0,t}u(t)$.
If   $u(t)$ satisfies \eqref{s4:fut}
and $\gamma_1<1$, then
$E^{(n,0,\sigma)}_2(u)=u^n-2u^{n-1} +u^{n-2}=O(\tau^2t_n^{\gamma_1-2})=O(\tau^{\gamma_1})$
for a small $n$, which yields large  discretization error of
\eqref{s4:eq-1}. We also use the correction terms to achieve
$E^{(n,m_u,\sigma)}_2(u)=O(\tau^2)$ for all $n>0$ when
$\gamma_{m_u+1}\geq 2$ (see \eqref{s4:Rn}).
See \cite{DieFord06,Lub86,ZengZK17,ZengTBK2018} for more details.
\end{remark}
}

\subsection{The fully discrete scheme}
We apply the time-stepping method \eqref{s4:eq-1} to discretize the
time direction of \eqref{s5:eq-1} and obtain the \sjwchange{following
semi-discrete scheme: For all $n>\max\{m,m_u,m_f\}$,} find
$U^n=U^n(x,y)$ such that
\begin{equation} \label{IMEX-FPDE}\begin{aligned}
D_{\tau}^{(\alpha,n,m,\gamma)}U
= \nu  \Delta U^n +  F^n - E^{(n,m_f,\hat{\gamma})}_2(F)-\kappa  E^{(n,m_u,\gamma)}_2(U),
\end{aligned}\end{equation}
where $F^n =f(x,y,U^n,t_n)$, $\kappa$ is a nonnegative constant.
\sjwchange{(Nonnegative $\kappa$ helps to enhance the stability
of \eqref{IMEX-FPDE}.)}
 {The method \eqref{IMEX-FPDE} is unconditionally stable for   $\tau>0$ if $\kappa > 0.75\max|\px[u]f|$; see \cite[Theorem~3]{ZengTBK2018} and its numerical verification. In real applications, $\kappa$ can be
estimated. For example, we can apply the fully implicit method with
coarse grid to solve \eqref{s4:fiv} to get an approximation of $u$, \sjwchange{and use this approximation to estimate $\max|\px[u]f|$.}}

Rewriting \eqref{IMEX-FPDE}  as the following linear boundary value problem
\begin{equation} \label{IMEX-FPDE-2}\left\{\begin{aligned}
&\left(\omega_0^{(\alpha)}+\kappa\tau^{\alpha}\right) U^n
- \nu \tau^{\alpha} \Delta U^n= RHS^{n-1},&(x,y)\in\Omega, \\
&U^n(x,y)=u_b(x,y,t_n),&(x,y)\in\px[]\Omega,
\end{aligned}\right.\end{equation}
 where
\begin{equation} \label{IMEX-FPDE-3}\begin{aligned}
RHS^{n-1}  = &RHS^{n-1}(x,y)=\tau^{\alpha}D_{\tau}^{(\alpha,n,m,\gamma)}U-\omega_0^{(\alpha)}U^n\\
&+F^n - E^{(n,m_f,\hat{\gamma})}_2(F)+\kappa (U^n-  E^{(n,m_u,\gamma)}_2(U)).
\end{aligned}\end{equation}

Let
$U_h^n=\Phi^T(x,y) \mathbf{c}^n\in\mathcal{M}_{\Omega}(\delta_x\times\delta_y)
$ be the approximate solution of \eqref{IMEX-FPDE-2}.  By inserting $U_h^n$
into the first equation of \eqref{IMEX-FPDE-2} and choosing
$(x,y)=\xi_i^{(in)}\in \mathcal{S}^{(p,q)}_{\Omega}$, we obtain
\begin{equation} \label{IMEX-FPDE-4}\begin{aligned}
\left(\omega_0^{(\alpha)}+\kappa\tau^{\alpha}\right) U_h^n(\xi_i^{(in)})
- \nu \tau^{\alpha} \Delta U^n_h(\xi_i^{(in)}) = RHS^{n-1}(\xi_i^{(in)}),
\quad 1\leq i \leq N_{in},
\end{aligned}\end{equation}
subject to the constraints
\begin{equation}\label{CS-1}\begin{aligned}
U_h^n(\xi^{(b)}_i) = \Phi^T(\xi^{(b)}_i)\mathbf{c}^n
= u_b(\xi^{(b)}_i,t_n), {\quad}1\leq i \leq N_b.
\end{aligned}\end{equation}
As in \eqref{s31:LeastSquare-22},
the discrete LS solution of \eqref{IMEX-FPDE-4}--\eqref{CS-1} is approximated by
\begin{equation}\label{LeastSquare}\begin{aligned}
\mathbf{c}^n
&=\arg\min_{\mathbf{d}\in  \mathbb{R}^{M} }\left\|\left(
           \begin{array}{c}
             (\omega_0^{(\alpha)}+\kappa\tau^{\alpha})\mathbf{A}_{in}-\nu\tau^{\alpha}\mathbf{S}_{in} \\
             \lambda\mathbf{A}_b \\
             \sqrt{\delta}\mathbf{I}\\
           \end{array}
         \right)\mathbf{d}
         -\left(
           \begin{array}{c}
             \mathbf{RHS}^{n-1}_{in} \\
             \lambda\mathbf{u}_b^n  \\
             \sqrt{\delta}\mathbf{d}^n_{*,r}\\
           \end{array}
         \right)
\right\|,
\end{aligned}\end{equation}
where  $\mathbf{A}_{in}$, $\mathbf{A}_{b}$, and $\mathbf{S}_{in}$  are defined in \eqref{s31:A},
$\delta\geq 0$,  $\mathbf{RHS}^{n-1}_{in}$ is  a
vector whose $i$th component is ${RHS}^{n-1}
(\xi^{(in)}_i) $;
 $\mathbf{d}^n_{*,r}$ is \sjwchange{a reference solution,} chosen such
that $\|\mathbf{c}^n-\mathbf{d}^n_{*,r}\|$ is suitably
small; \sjwchange{and $\lambda$ is the positive weight parameter that
controls fidelity to the boundary conditions}.  We choose the
reference solution $\mathbf{d}^n_{*,r}$ to be
\begin{equation}\label{LeastSquare-01}
\mathbf{d}^n_{*,r}=\left\{\begin{aligned}
&\mathbf{0},{\qquad\qquad} r = 0,\\
&\mathbf{c}^{n-1},{\qquad\,\,\,}r=1.
\end{aligned}\right.\end{equation}
The choice of $\mathbf{d}^n_{*,r}$ is inspired by the relation
$\|\mathbf{c}^{n}-\mathbf{d}^n_{*,r}\|=O(\tau^r)$ observed from the
time-stepping methods for solving time-dependent problems, which is
numerically verified in our computational results of the next section.


The LSFVM for \eqref{IMEX-FPDE-4} can be obtained similarly, where we
replace the matrices $\mathbf{A}_{in} $ and $\mathbf{S}_{in} $
in \eqref{LeastSquare} with $\mathbf{\widetilde{A}}_{in} $ and
$\mathbf{\widetilde{S}}_{in} $, respectively, where
$\mathbf{\widetilde{A}}_{in}$ and $\mathbf{\widetilde{S}}_{in}$ are
defined in \eqref{s5:eq-7}.



\section{Numerical examples and applications}\label{sec:numerical}

We present Examples \ref{s5-eg-3} and \ref{s5-eg-4}, to display the
efficiency of the LSC method \eqref{LeastSquare}.  For each example,
we consider two different cases based on a known solution with a
source term and an unknown solution with zero source but oscillating
initial conditions. \sjwchange{Since LSFVM gives similar results to
the LSC when $\rho$ is suitably small, we do not discuss the results
of LSFVM further.}

The following four different domains are considered for
Example \ref{s5-eg-3}.
\begin{itemize}
  \item[(\textbf{i})] the rectangular domain $\Omega=(-1,1)\times (-1,1)$, see Figure~\ref{fig:iregular_domain}(a);
  \item[(\textbf{ii})] the circular domain $\Omega=\{(x,y)|x^2+y^2<1\}$, see Figure~\ref{fig:iregular_domain}(b);
  \item[(\textbf{iii})] the irregular domain $\Omega$ defined by Figure~\ref{fig:iregular_domain}(c), see  Case II in Example~\ref{eg2-1};
   \item[(\textbf{iv})] the irregular domain $\Omega$ defined by Figure~\ref{fig:iregular_domain}(d), see  Case III in Example~\ref{eg2-1}.
\end{itemize}
We consider Example~\ref{s5-eg-4} on the irregular domain
(\textbf{iv}).

\begin{example}\label{s5-eg-3}
Consider the following nonlinear equation
\begin{equation}\label{s5eg1:eq-1}
{}_{C}D^{\alpha}_{0,t}u(x,y,t) = \Delta u(x,y,t) + u(1-u^2)+g(x,y,t),
{\quad}(x,y,t)\in\Omega\times(0, T],
\end{equation}
subject to  suitable initial and Dirichlet boundary conditions.
\end{example}

The following two cases are considered in this example.
\begin{itemize}
\item Case I: Choose initial and boundary conditions, and source
  term $g(x,y,t)$, such that the analytical solution
  of \eqref{s5eg1:eq-1} is \begin{equation}\label{eq:s5-eg-2} u(x,y,t)
  = E_{\alpha}(-t^{\alpha})\sin(x)\sin(y).  \end{equation}
\item Case  II: Choose initial value $u(x,y,0)=\sin(\pi x)\sin(\pi y)$, source
  term $g=0$, and homogenous boundary conditions.
\end{itemize}

We always choose $m_f=m_u=m$ with correction indices
$\gamma_k=k\alpha$ (see \cite{ZengTBK2018,ZengZK17} and \eqref{s4:fut}
for choosing $\hat{\gamma}_k$), $\lambda=10^5$, and $N_x=N_y=N$
when \eqref{s5eg1:eq-1} is solved
by \eqref{LeastSquare}. If \eqref{s5eg1:eq-1} is solved on the
rectangular domain, \sjwchange{then \eqref{LeastSquare} is relatively
well-conditioned}, so we take $\delta=0$ in the computation.  When
$\alpha=1$, this model reverts to a standard diffusion equation.

We verify the accuracy of the LSC method \eqref{LeastSquare} in Case
I, where the analytical solution is given explicitly on four different
domains.  Tables~\ref{eg1:tb1} shows \sjwchange{relative $L^2$ errors
for this approach} at $t=2$ on the rectangular domain.  Second-order
accuracy is observed once again.


\begin{table}[!ht]
\caption{{The  $L^2$ error $\|e^n\|$ of the LSC method  \eqref{LeastSquare} at $t=2$ on
the rectangular domain (i), Example \ref{s5-eg-3}, Case I, with $\kappa=2$, $\gamma_k=k\alpha$,   $\tau=2^{-10}$.}}\label{eg1:tb1}
\centering\footnotesize
\begin{tabular}{|c|c|c|c|c|c|c|c|c|c|c|c|c|}
\hline
 $N$ & $\alpha=0.1,m=3$ & Order& $\alpha=0.5,m=3$ & Order & $\alpha=0.8,m=1$ & Order & ${ \alpha=1,m=1}$ &{ Order}\\
 \hline
$8 $&3.4971e-6&      &2.5365e-6&      &1.7604e-6&      &{ 1.1468e-6}&                  \\
$16$&6.1913e-7&2.4979&4.5642e-7&2.4744&3.1358e-7&2.4890&{ 1.8610e-7}&{ 2.6235}\\
$32$&1.7198e-7&1.8480&1.2340e-7&1.8870&8.5040e-8&1.8826&{ 4.5807e-8}&{ 2.0224}\\
$48$&7.8119e-8&1.9463&5.4424e-8&2.0190&3.8185e-8&1.9747&{ 1.6985e-8}&{ 2.4468}\\
$64$&4.4314e-8&1.9707&2.9698e-8&2.1056&2.1401e-8&2.0126&{ 6.6661e-9}&{ 3.2511}\\
\hline
\end{tabular}
\end{table}

\begin{figure}[!ht]
\begin{center}
\begin{minipage}{0.47\textwidth}\centering
\epsfig{figure=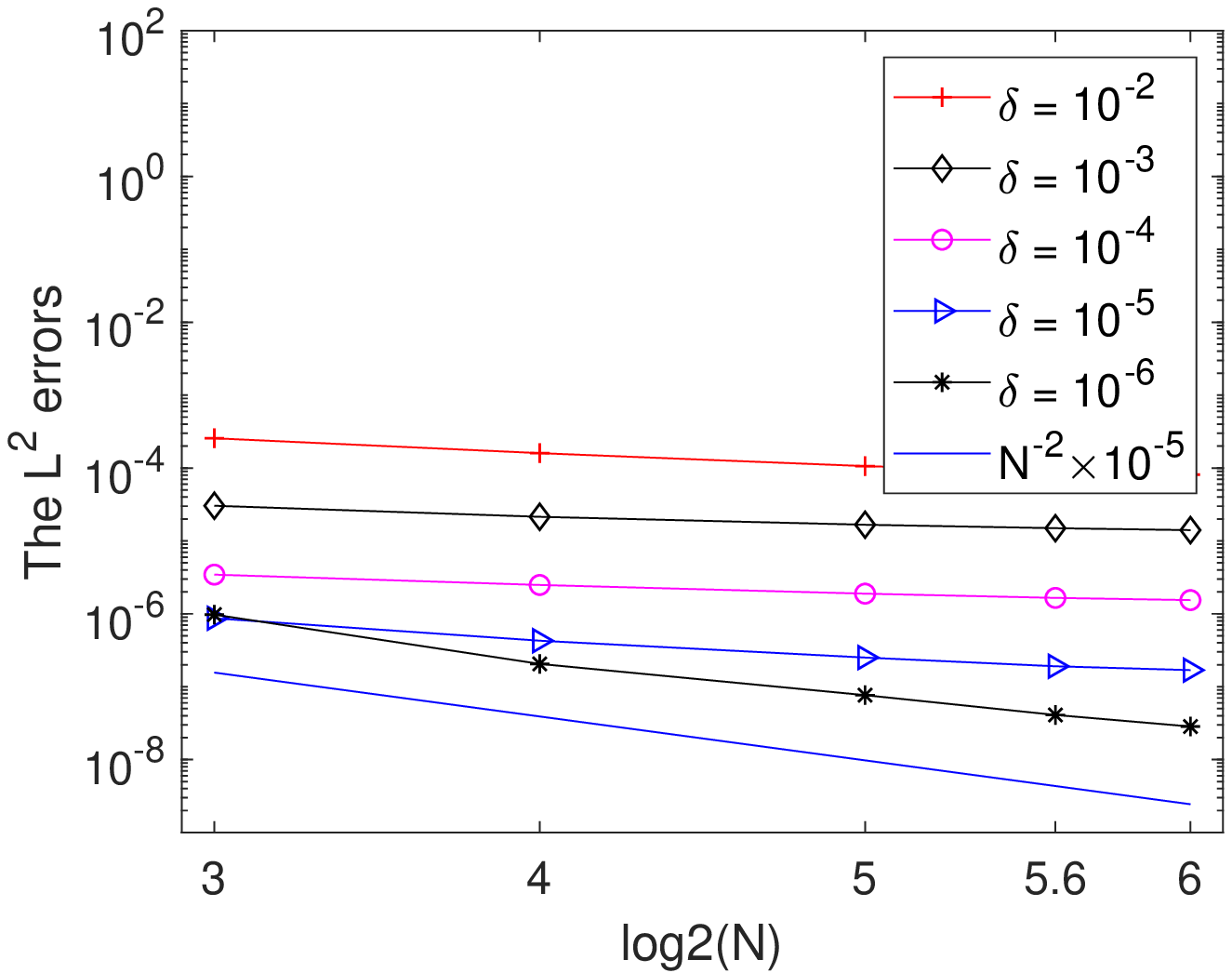,width=7cm} \par {(a) $r=0$.}
\end{minipage}
\begin{minipage}{0.47\textwidth}\centering
\epsfig{figure=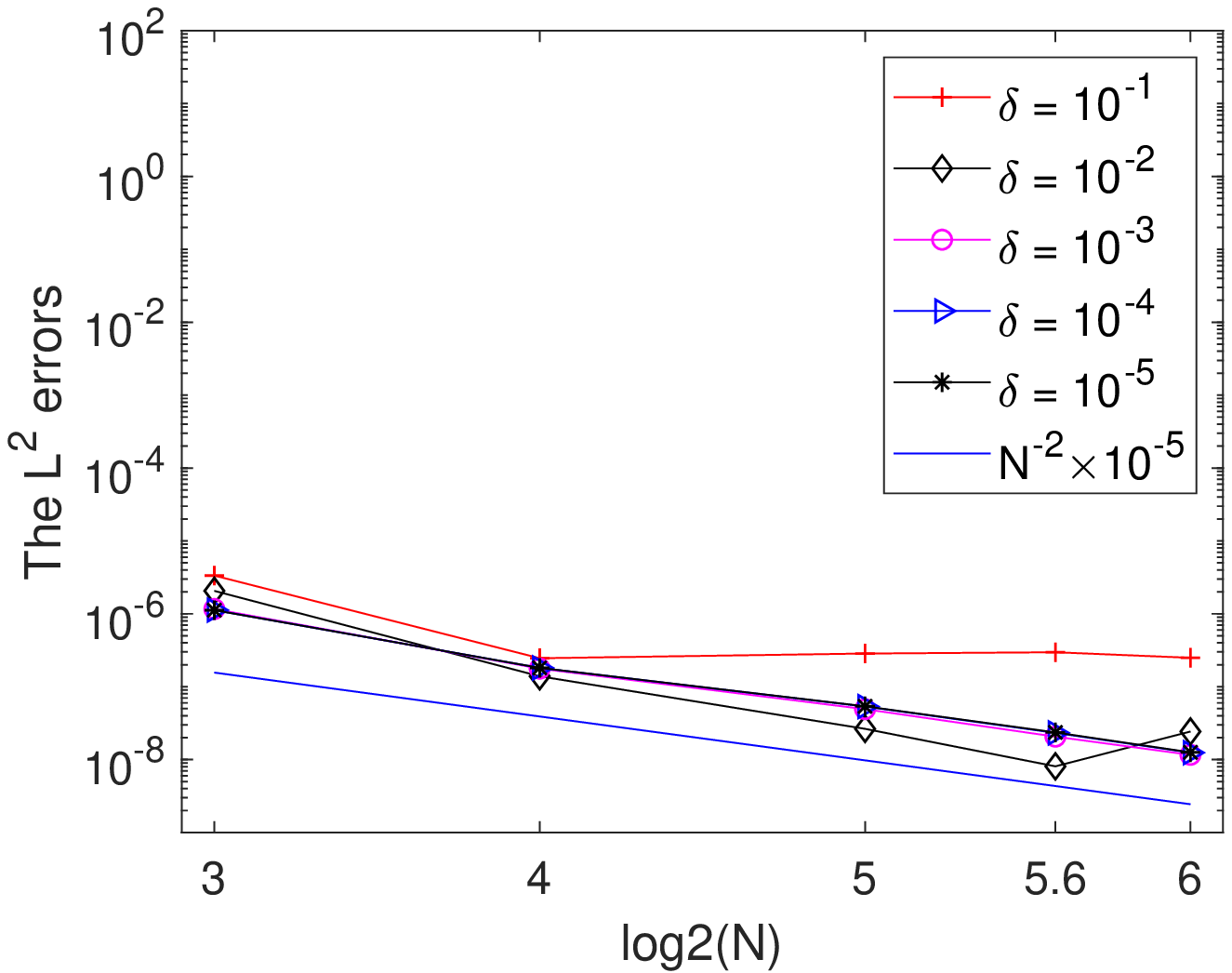,width=7cm} \par {(b) $r=1$.}
\end{minipage}
\end{center}
\caption{The $L^2$ errors at $t=2$ on the circular domain (ii), Example \ref{s5-eg-3}, Case I,
$\kappa=2$, $\gamma_k=k\alpha$, $\alpha=0.5$,  $\tau=2^{-10},m=3$.\label{fig4-1}}
\end{figure}

\begin{figure}[!ht]
\begin{center}
\begin{minipage}{0.47\textwidth}\centering
\epsfig{figure=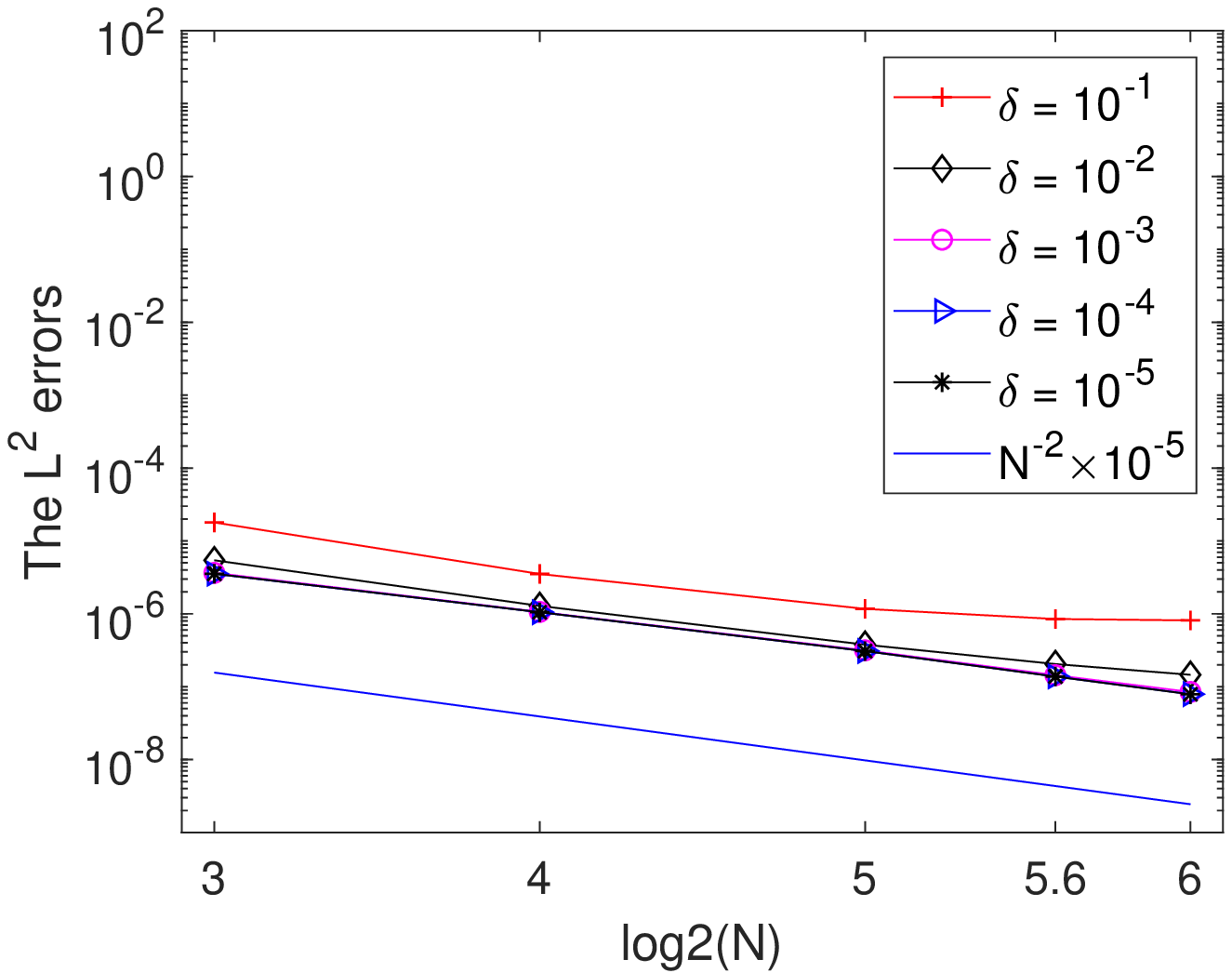,width=7cm} \par {(a) Irregular domain (iii), $r=1$.}
\end{minipage}
\begin{minipage}{0.47\textwidth}\centering
\epsfig{figure=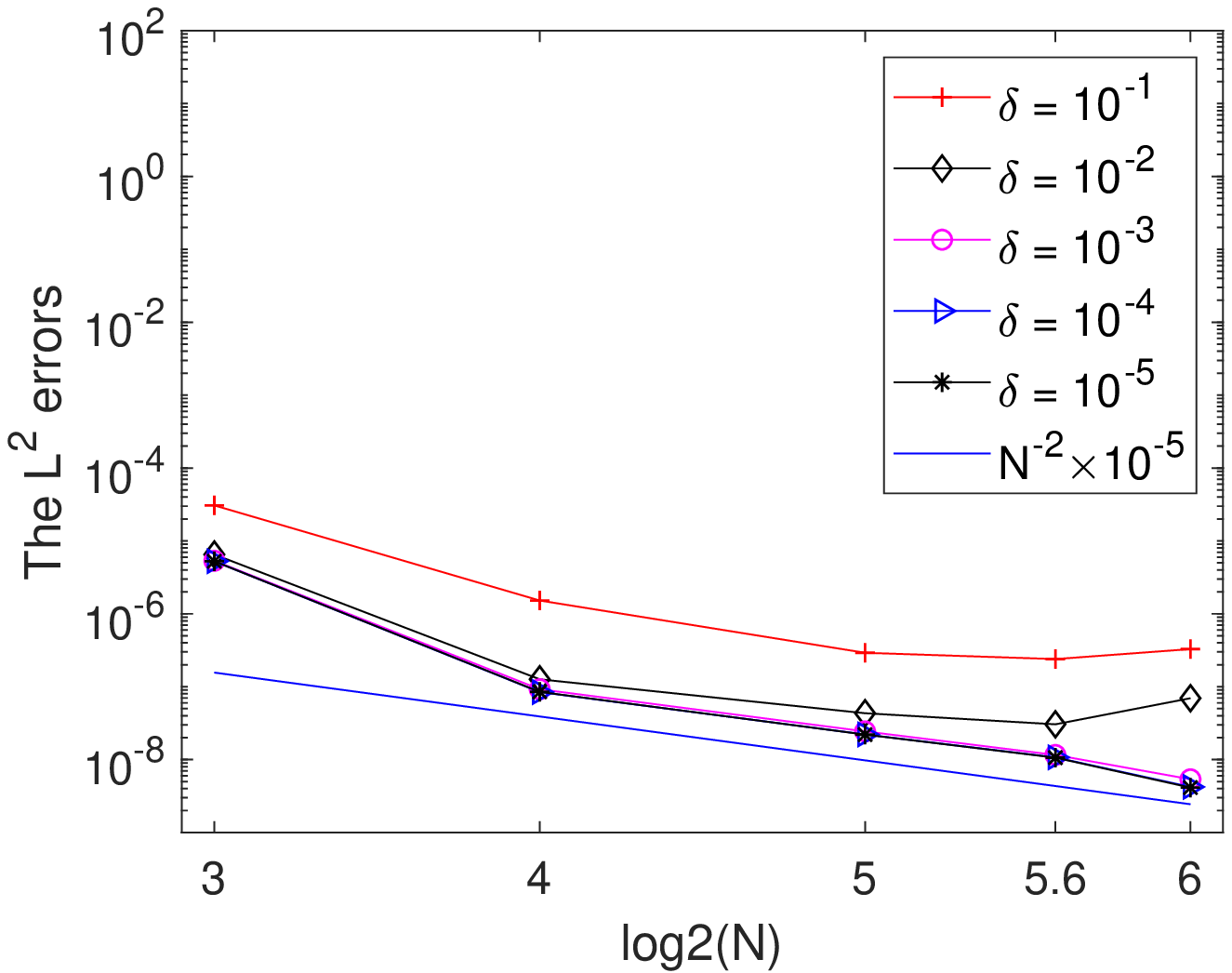,width=7cm} \par {(b) Irregular domain (iv), $r=1$.}
\end{minipage}
\end{center}
\caption{The $L^2$ errors at $t=2$  on the irregular domains (iii) and (iv), Example \ref{s5-eg-3}, Case I, $\kappa=2$, $\gamma_k=k\alpha$, $\alpha=0.5$, $\tau=2^{-10},m=3$.\label{fig4-2}}
\end{figure}

\sjwchange{When the LSC method \eqref{LeastSquare} is applied to
solve \eqref{s5eg1:eq-1} on irregular domains (ii)-(iv), the condition
number of the coefficient matrix in \eqref{LeastSquare} may be large,
but regularization helps to improve the conditioning} while still
yielding accurate numerical solutions. Figure~\ref{fig4-1}(a) shows
the $L^2$ error of the LSC method
\eqref{LeastSquare} at $t=2$ on the circular domain (ii) for $r=0$.
\sjwchange{We observe second-order
accuracy again for $\delta=10^{-6}$.}
\sjwchange{Figure~\ref{fig4-1}(b) shows these errors for $r=1$, where
 we observe more accurate numerical simulations even for larger values
of $\delta$. (Second-order accuracy can be observed even for
$\delta=10^{-2}$.)}  {  For $\delta=10^{-1}$ and $r=1$, the
regularization error plays the dominant role, we observe that the
error is saturated as $N$ increases up to $64$; see
Figure~\ref{fig4-1}(b).  For $\delta=10^{-2}$ and $r=1$, the space
discretization error plays the dominant role when $N\leq 32$. For
$N\in (32,64]$, both the space discretization error and the
regularization error play important roles; \sjwchange{we observe in
Figure~\ref{fig4-1}(b) and Figure~\ref{fig4-2}(b) that these errors
fluctuate slightly with $N$.}  For the irregular domains (iii) and
(iv), we observe similar results as those shown in
Figure~\ref{fig4-1}. Figure~\ref{fig4-2} shows the $L^2$ errors on the
irregular domains (iii) and (iv) for $r=1$.  {  From these
figures, we can see that a relatively large $\delta$ with $r=1$ yields
a well-conditioned system \eqref{LeastSquare} and accurate numerical
solutions, in agreement with the results in Example \ref{eg2-1}.
\sjwchange{In the remainder of this section,
we will use $r=1$ with a relatively large value of $\delta$ in the
computations.} }}

%

\begin{figure}[!ht]
\begin{center}
\begin{minipage}{1\textwidth}\centering
\epsfig{figure=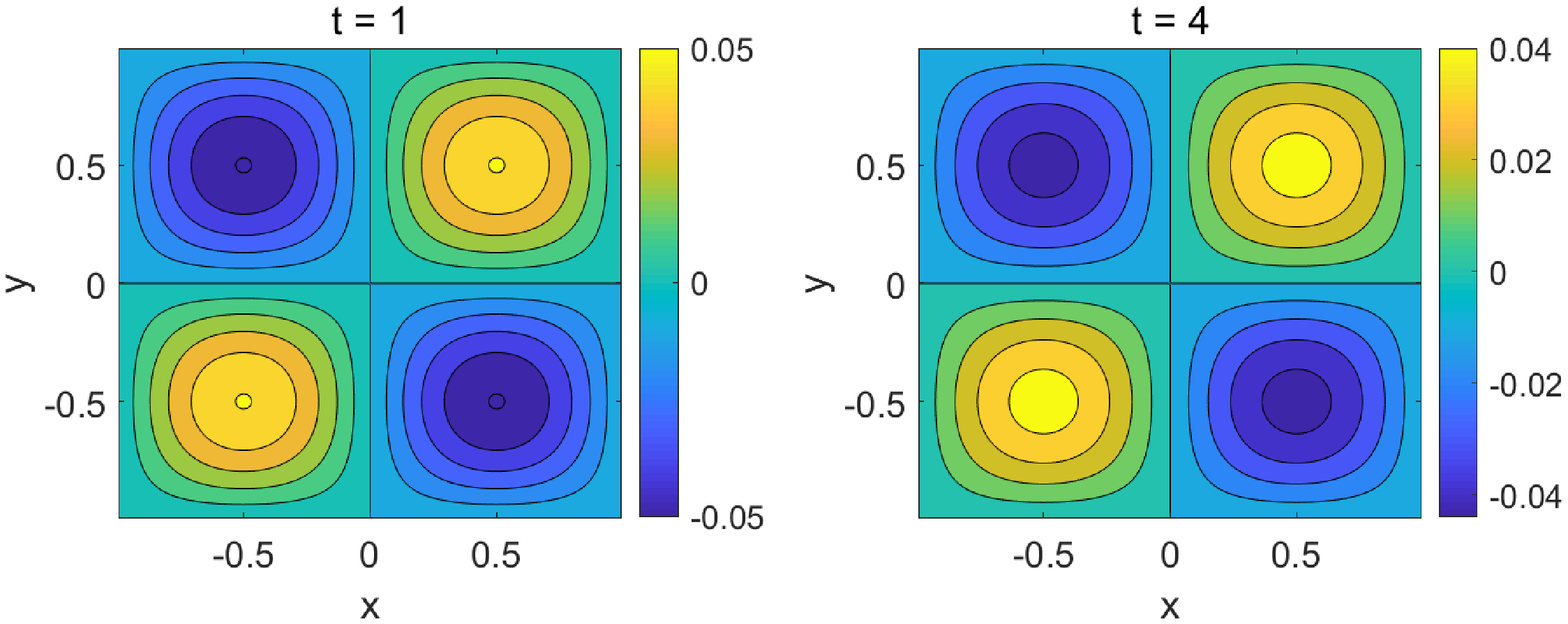,width=14cm}\par {(a) $\alpha=0.2$.}
\end{minipage}
\begin{minipage}{1\textwidth}\centering
\epsfig{figure=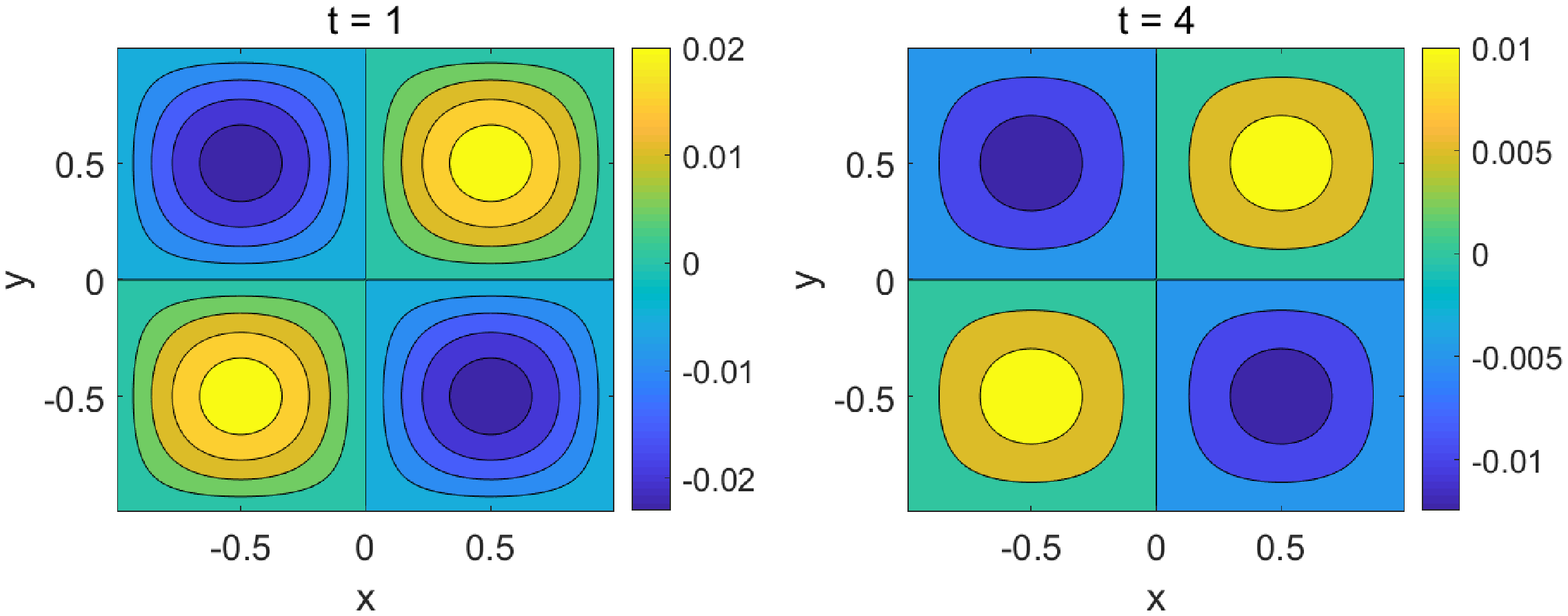,width=14cm}\par {(b) $\alpha=0.8$.}
\end{minipage}
\end{center}
\caption{Numerical solutions on the rectangular domain (i) for Example \ref{s5-eg-3}, Case II,
$\gamma_k=k\alpha$, $\kappa=2,\tau=2^{-7},m=1$, $\mu=1$.\label{eg1fig1-1}}
\end{figure}

\begin{figure}[!ht]
\begin{center}
\begin{minipage}{1\textwidth}\centering
\epsfig{figure=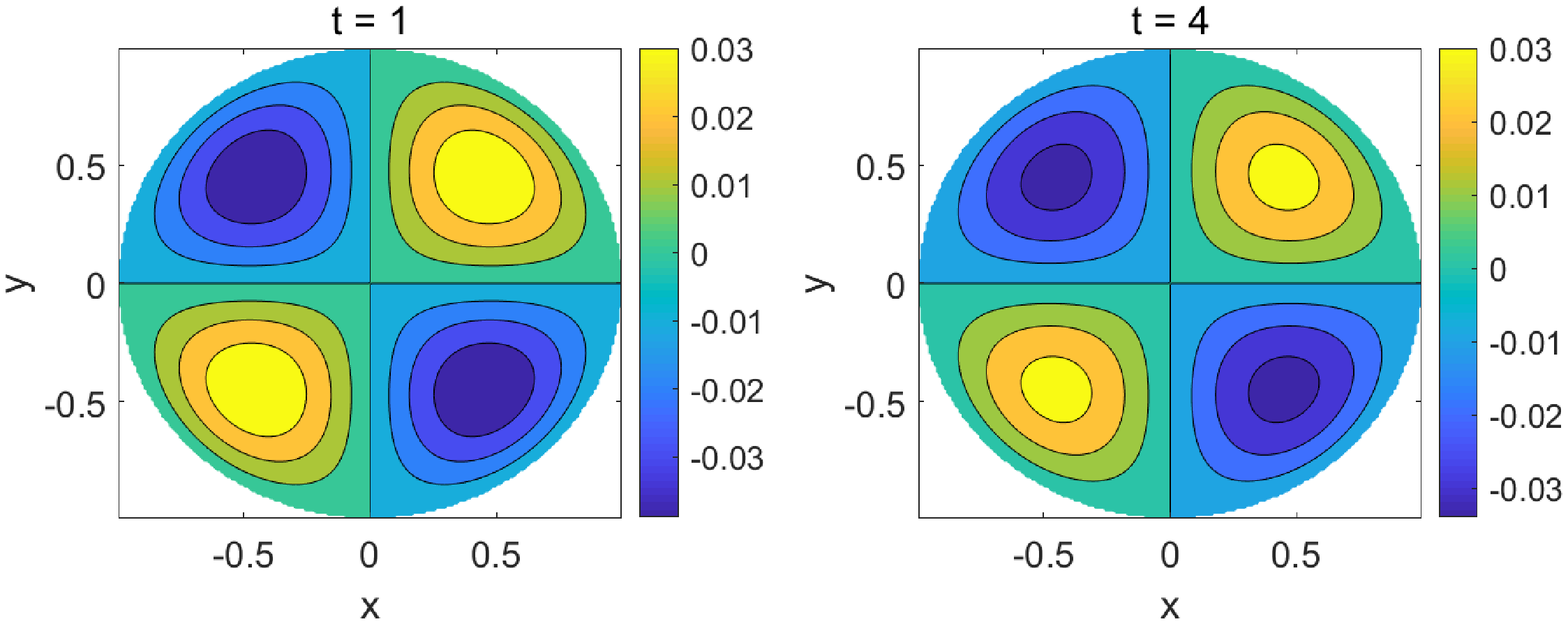,width=14cm} \par {(a) $\alpha=0.2$.}
\end{minipage}
\begin{minipage}{1\textwidth}\centering
\epsfig{figure=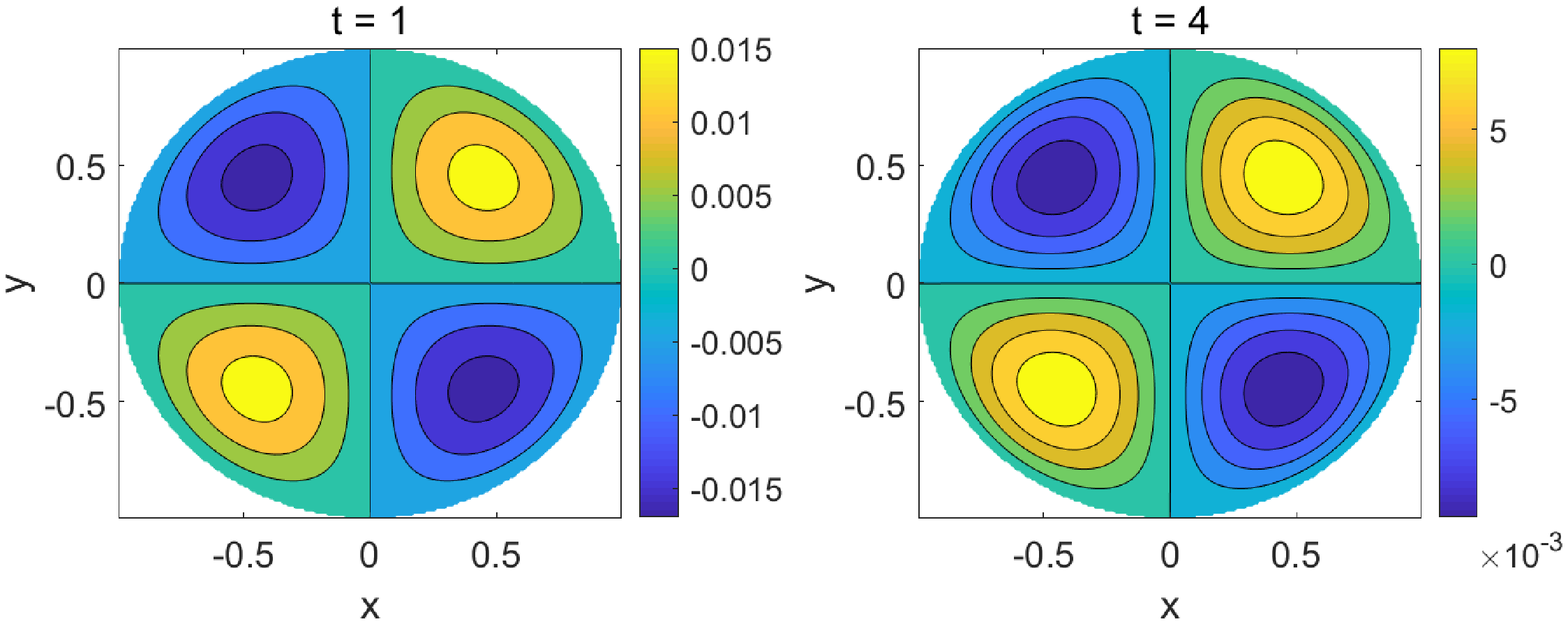,width=14cm} \par {(b) $\alpha=0.8$.}
\end{minipage}
\end{center}
\caption{Numerical solutions on the circular domain (ii) for Example \ref{s5-eg-3}, Case II,
$\gamma_k=k\alpha$,
 $\kappa=2,\tau=2^{-7},m=1$, $\mu=1$.\label{eg1fig2-1}}
\end{figure}

\begin{figure}[!ht]
\begin{center}
\begin{minipage}{1\textwidth}\centering
\epsfig{figure=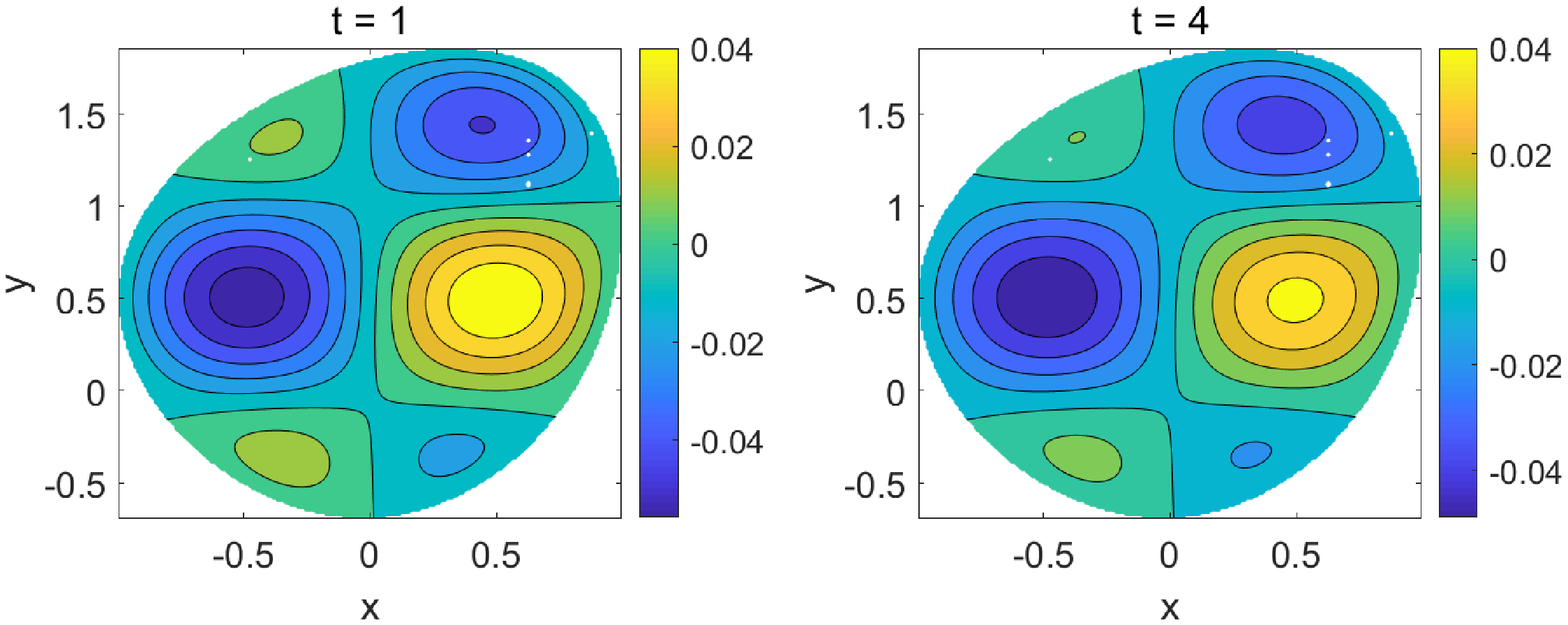,width=14cm} \par {(a) $\alpha=0.2$.}
\end{minipage}
\begin{minipage}{1\textwidth}\centering
\epsfig{figure=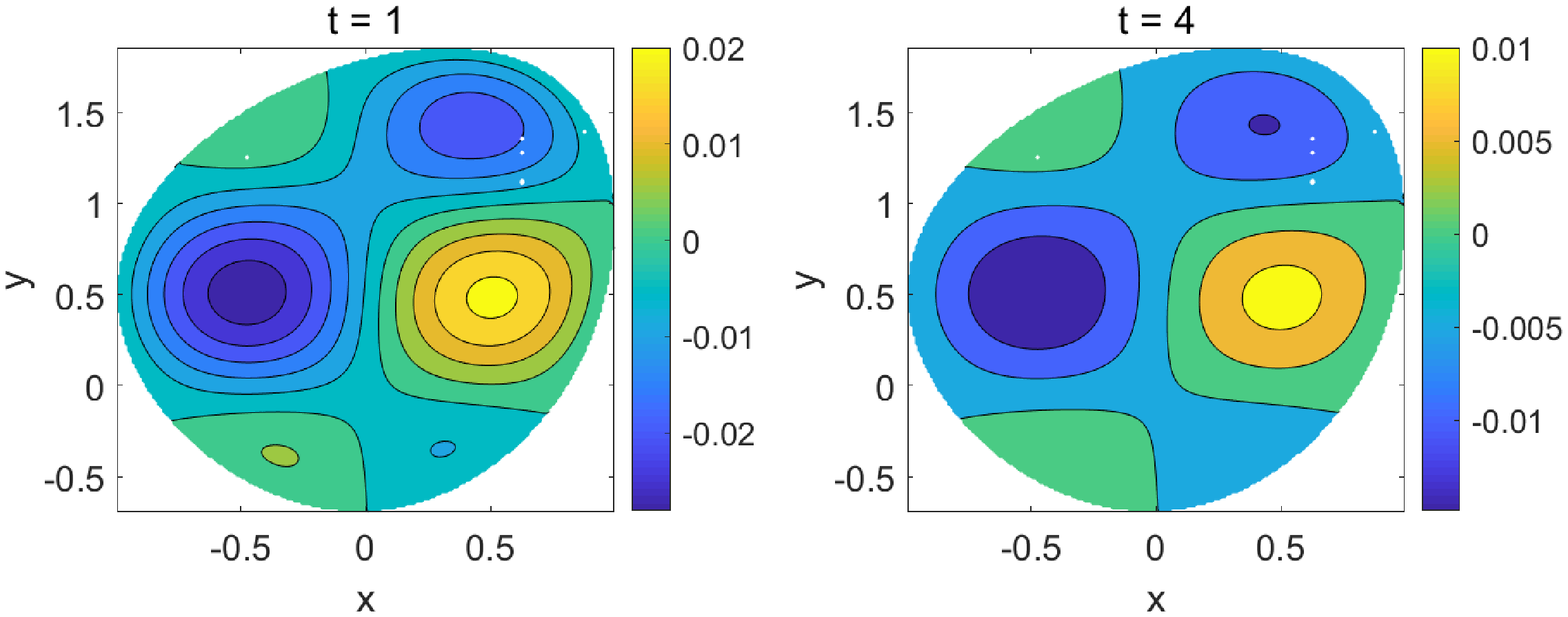,width=14cm} \par {(a) $\alpha=0.8$.}
\end{minipage}
\end{center}
\caption{Numerical solutions on the irregular domain (iii) for Example \ref{s5-eg-3}, Case II,
$\gamma_k=k\alpha$, $\kappa=2,\tau=2^{-7},m=1$, $\mu=1$.\label{eg1fig3-1}}
\end{figure}

\begin{figure}[!ht]
\begin{center}
\begin{minipage}{1\textwidth}\centering
\epsfig{figure=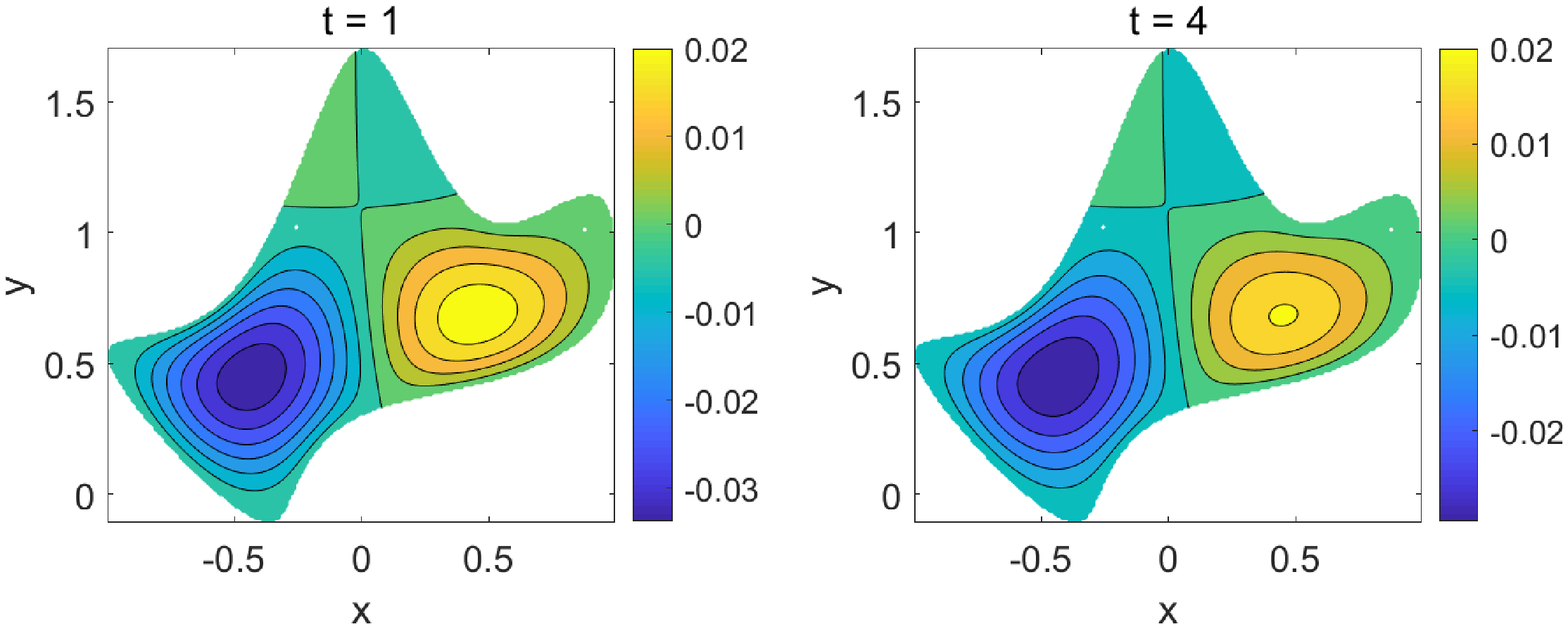,width=14cm} \par {(a) $\alpha=0.2$.}
\end{minipage}
\begin{minipage}{1\textwidth}\centering
\epsfig{figure=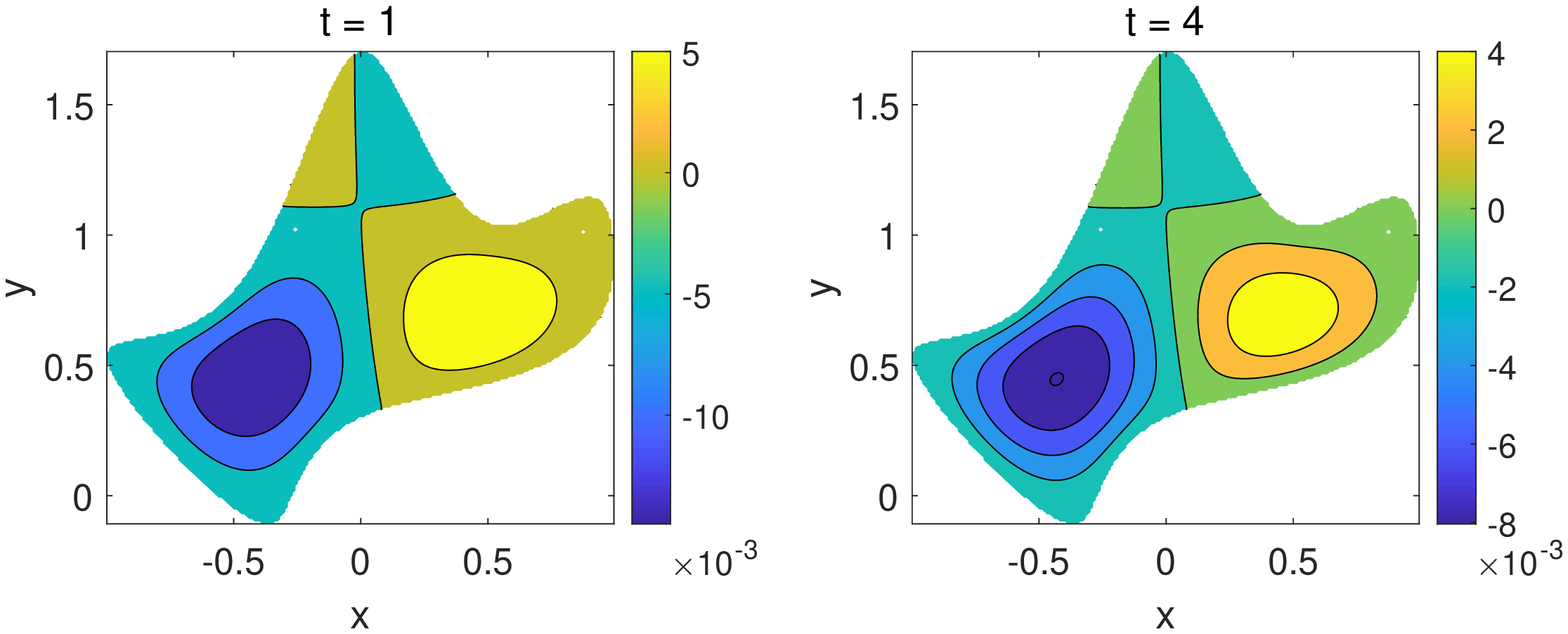,width=14cm} \par {(a) $\alpha=0.8$.}
\end{minipage}
\end{center}
\caption{Numerical solutions on the irregular domain (iv) for Example \ref{s5-eg-3}, Case II,
$\gamma_k=k\alpha$, $\kappa=2,\tau=2^{-7},m=1$, $\mu=1$.\label{eg1fig4-1}}
\end{figure}

\sjwchange{For Case II, we obtain a proxy for the exact solution} by using the LSC
method \eqref{LeastSquare} with $r=1$, $\delta=0.01$, and $N=64$.
Figures~\ref{eg1fig1-1}(a) and (b) show the numerical solution on the
rectangular domain (i) at different times $t=1$ and $t=4$. We observe
that the solution decays as $t$ evolves. Comparing
Figure~\ref{eg1fig1-1}(a) with Figure~\ref{eg1fig1-1}(b), we observe
that the solution decays faster as $\alpha$ becomes larger.  For the
irregular domains defined by (ii), (iii), and (iv), we observe similar
results as those shown in Figure \ref{eg1fig1-1} for the rectangular
domain (i), see Figures~\ref{eg1fig2-1}-\ref{eg1fig4-1}.

\begin{figure}[!h]
\begin{center}
\begin{minipage}{0.47\textwidth}\centering
\epsfig{figure=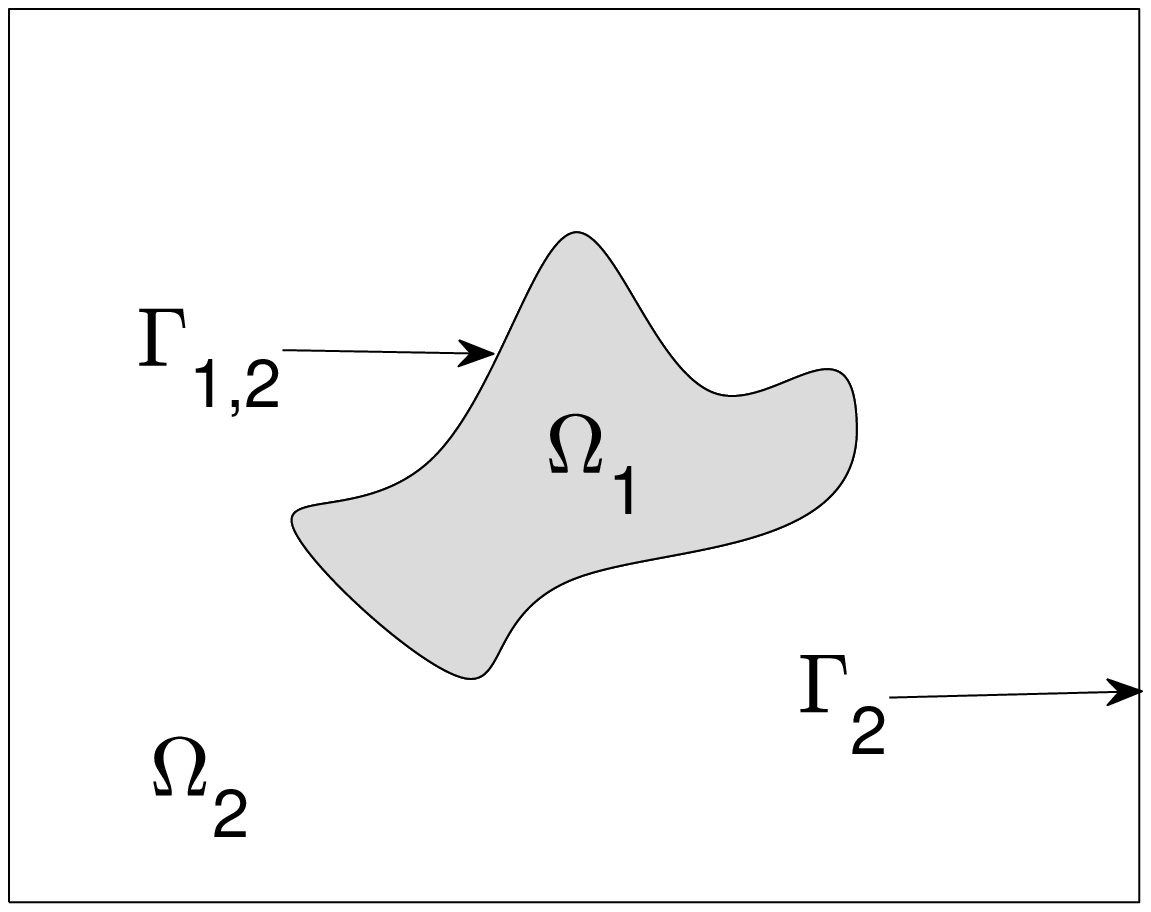,width=6cm}
\end{minipage}
\begin{minipage}{0.47\textwidth}\centering
\epsfig{figure=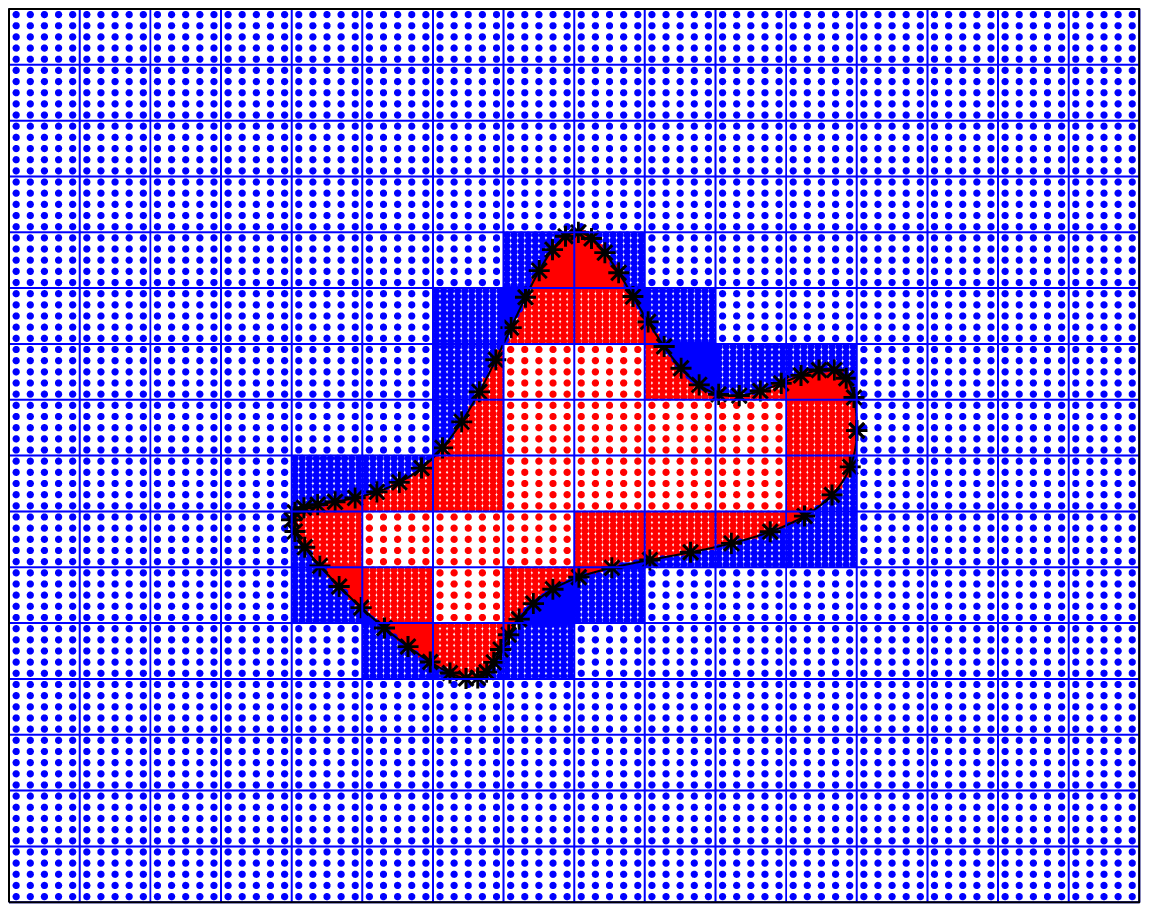,width=6cm}
\end{minipage}
\end{center}
\caption{The  domains $\Omega_1$ and $\Omega_2$ with the shared boundary $\Gamma_{1,2}$ (left) and
their division and distribution of collocation points (right).
\label{fig:domain2}}
\end{figure}

In the following example, we extend the LSC method \eqref{LeastSquare}
to a system of equations with different fractional indices.
\begin{example}\label{s5-eg-4}
Consider the following system of equations
\begin{equation}\label{s5eg1:eq-2}\left\{\begin{aligned}
&{}_{C}D^{\alpha}_{0,t}u(x,y,t) = \mu\Delta u(x,y,t)+f(x,y,t),{\quad}(x,y,t)\in \Omega_1\times(0,T],\\
&{}_{C}D^{\beta}_{0,t}v(x,y,t) = \nu\Delta v(x,y,t)+g(x,y,t),{\quad}(x,y,t)\in \Omega_2\times(0,T],
\end{aligned}\right.\end{equation}
subject to suitable initial  conditions and the following boundary conditions
\begin{equation}\label{s5eg1:eq-3}\left\{\begin{aligned}
&v(x,y,t)=v_b(x,y,t), {\qquad}&&(x,y)\in \Gamma_{2},\\
&u(x,y,t)=v(x,y,t),{\quad}\mathbf{n}\cdot\nabla u = \mathbf{n}\cdot\nabla v,
{\qquad} &&(x,y)\in \Gamma_{1,2}=\px[]\Omega_{1},
\end{aligned}\right.\end{equation}
where, in Figure~\ref{fig:domain2}, $\Omega_1$ is the area inside the
curve $\Gamma_{1,2}$ (the shaded area) and $\Omega_2$ is the area
between the curves $\Gamma_{1,2}$ and $\Gamma_2$.
\end{example}

We first apply the time-stepping method \eqref{s4:eq-1} to each equation of
\eqref{s5eg1:eq-2} to obtain
\begin{equation}\label{s5eg1:eq-4}\left\{\begin{aligned}
&D_{\tau}^{(\alpha,n,m,\gamma)}u= \mu \Delta u^n +f^n + R_1^n,\\
&D_{\tau}^{(\beta,n,\tilde{m},\tilde{\gamma})}v = \nu \Delta v^n +
g^n+ R_2^n,
\end{aligned}\right.\end{equation}
where $D_{\tau}^{(\alpha,n,m,\gamma)}$ is defined by \eqref{s2:Dalf},
and $R_1^n$ and $R_2^n$ are truncation errors in time that depend on
the regularity of $u$ and $v$, respectively.  Omitting the truncation
errors in \eqref{s5eg1:eq-4}, we derive the following semi-discrete
method for
\eqref{s5eg1:eq-2}:
\begin{equation} \label{s5eg1:eq-5}\left\{\begin{aligned}
&\omega_0^{(\alpha)} U^n  - \mu \tau^{\alpha} \Delta U^n = RHS_1^{n-1}
=\tau^{\alpha}D_{\tau}^{(\alpha,n,m,\gamma)}U-\omega_0^{(\alpha)}U^n +\tau^{\alpha}f^n,\\
&\omega_0^{(\beta)} V^n  - \nu \tau^{\beta} \Delta V^n = RHS_2^{n-1}
=\tau^{\beta}D_{\tau}^{(\beta,n,\tilde{m},\tilde{\gamma})}V-\omega_0^{(\beta)}V^n +\tau^{\beta}g^n,
\end{aligned}\right.\end{equation}
subject to the boundary conditions
\begin{equation}\label{s5eg1:eq-7}\left\{\begin{aligned}
&V^n(x,y)=v_b(x,y,t_n), {\qquad}&&(x,y)\in \Gamma_{2},\\
&U^n(x,y)=V^n(x,y),{\quad}\mathbf{n}\cdot\nabla U^n = \mathbf{n}\cdot\nabla V^n,
{\qquad} &&(x,y)\in \Gamma_{1,2}=\px[]\Omega_{1}.
\end{aligned}\right.\end{equation}

The penalized LSC formulation
of \eqref{s5eg1:eq-5}--\eqref{s5eg1:eq-7} is
\begin{equation}\label{s5eg1:eq-8}\begin{aligned}
&\sum_{i=1}^{N^u_{in}}\left[  J^n(\xi^{(in)}_i,\mathbf{d})   \right]^2
 +\sum_{i=1}^{N^v_{in}}\left[  \hat{J}^n(\eta^{(in)}_i,\mathbf{\hat{d}})   \right]^2
 + \lambda^2\sum_{i=1}^{N^v_{b}}\left(R^n_v(\eta^{(b)}_i,\mathbf{\hat{d}})\right)^2\\
&{\qquad}+\lambda^2\sum_{i=1}^{N_{c}}\left[\left(R^n_{u=v}(\varsigma_i,\mathbf{d},\mathbf{\hat{d}})\right)^2
+\left(R^n_{\mathbf{n}\cdot\nabla}(\varsigma_i,\mathbf{d},\mathbf{\hat{d}})\right)^2\right]\\
&{\qquad}+\delta\left(\|\mathbf{d}-\mathbf{d}^{n-1}\|^2
+\|\mathbf{\hat{d}}-\mathbf{\hat{d}}^{n-1}\|^2\right),
\end{aligned}\end{equation}
where \sjwchange{$\delta\geq 0$ is the regularization parameter,
$\lambda$ is the penalty parameter for the boundary conditions,}
$N_{in}^u$ and $N_{in}^v$ are the numbers of collocation points in
$\Omega_1$ and $\Omega_2$, respectively, $N_{b}^v$ is the number of
boundary points on $\Gamma_2$, $N_{c}$ is the number of boundary
points on $\Gamma_{1,2}$, and
\begin{equation} \begin{aligned}
J^n(x,y,\mathbf{d})&=\omega_0^{(\alpha)}U^n_h(x,y)
-\mu\tau^{\alpha} \Delta U^n_h(x,y) -RHS_1^{n-1}(x,y),&& (x,y)\in \Omega_1, \\
\hat{J}^n(x,y,\mathbf{\hat{d}})&= \omega_0^{(\beta)}V^n_h(x,y)
-\nu\tau^{\beta} \Delta V^n_h(x,y)  -RHS_2^{n-1}(x,y),
 &&(x,y)\in \Omega_2,\\
R^n_v(x,y,\mathbf{\hat{d}})&=V^n_h(x,y)-v_b(x,y,t_n),&&{(x,y)\in \Gamma_2},\\
R^n_{u=v}(x,y,\mathbf{d},\mathbf{\hat{d}})
&=U^n_h(x,y)-V^n_h(x,y),&&(x,y)\in \Gamma_{1,2},\\
R^n_{\mathbf{n}\cdot\nabla}(x,y,\mathbf{d},\mathbf{\hat{d}})
&=\mathbf{n}\cdot\nabla U_h^n-\mathbf{n}\cdot\nabla V_h^n,&&(x,y)\in \Gamma_{1,2},\\
\end{aligned}\end{equation}
with
$U^n_h(x,y)=\Phi^T(x,y)\mathbf{d}^n\in\mathcal{M}_{\Omega_1}(\delta_x\times\delta_y),\mathbf{d}^n\in \mathbb{R}^{M_u}$
and
$V^n_h(x,y)=\widehat{\Phi}^T(x,y)\mathbf{\hat{d}}^n\in\mathcal{M}_{\Omega_2}(\delta_x\times\delta_y),\mathbf{\hat{d}}^n\in \mathbb{R}^{M_v}$.
The distribution of collocation points is shown as
Figure~\ref{fig:domain2} (right), which is defined similarly to
$\mathcal{S}^{(p,q)}_{\Omega}$; see \eqref{Spq}.

\begin{itemize}
  \item Case I: \sjwchange{Choose   initial and boundary conditions, and  source terms, to yield the following solution to} \eqref{s5eg1:eq-2}-\eqref{s5eg1:eq-3}:
  \begin{align*}
  u(x,y,t)&=E_{\alpha}(-t^{\alpha})\sin(3x)\sin(3y),{\quad}(x,y)\in \Omega_1,\\
  v(x,y,t)&=E_{\beta}(-t^{\beta})\sin(3x)\sin(3y),{\quad}(x,y)\in \Omega_2,
  \end{align*}
  where $\alpha=\beta$.
{  \item Case II: The initial conditions  are taken as $u(x,y,0)=\sin(2\pi x)\sin(2\pi y)$, $(x,y)\in \bar{\Omega}_1$,   $v(x,y,0)=\sin(2\pi x)\sin(2\pi y),(x,y)\in \bar{\Omega}_2$,   $f=0$ for $(x,y)\in {\Omega}_1$, and $g=0$ for $(x,y)\in {\Omega}_2$.}
\end{itemize}

In this example, \sjwchange{we set $N_x=N_y=N$, $m=\tilde{m}=1$,
$\gamma_1=\alpha$, $\tilde{\gamma}_1=\beta$, $\lambda=10^5$, and
$\delta=0.01$ in the LSC method} \eqref{s5eg1:eq-8}.  The
computational domain in Cases I and II satisfies
$\bar{\Omega}_1 \cup \bar{\Omega}_2=[-1.5,1.5]\times [-1,2]$, where
$\Omega_1$ is defined by Figure~\ref{fig:iregular_domain}(d); see also
(iv) at the beginning of this section.

We show the $L^2$ errors of $u$ and $v$ for Case I at $t=2$ in
Table \ref{eg2:tb1}.  We can see that satisfactory numerical solutions
are obtained, although the convergence rate in space is slightly less
than two, due to \sjwchange{ill conditioning in the least squares
formulation}.

\begin{table}[!ht]
\caption{{The  $L^2$ error of the LSC method \eqref{s5eg1:eq-8}  at $t=2$, Example \ref{s5-eg-4},
Case I,  $\alpha=\beta=0.5$,  and $\tau=2^{-10}$.}}\label{eg2:tb1}
\centering\footnotesize
\begin{tabular}{|c|c|c|c|c|c|c|c|c|c|c|c|c|}
\hline
 $N$ & $L^2$-error $(u)$ & Order& $L^2$-error $(v)$ & Order   \\
 \hline
$8 $&7.7227e-4&      &6.9974e-4&      \\
$16$&7.3497e-5&3.3933&6.2523e-5&3.4844\\
$32$&2.9992e-5&1.2931&2.7286e-5&1.1962\\
$50$&1.3148e-5&2.0338&1.1755e-5&2.0768\\
$64$&8.4209e-6&1.5489&7.5325e-6&1.5471\\
\hline
\end{tabular}
\end{table}

\begin{figure}[!h]
\begin{center}
\begin{minipage}{1\textwidth}\centering
\epsfig{figure=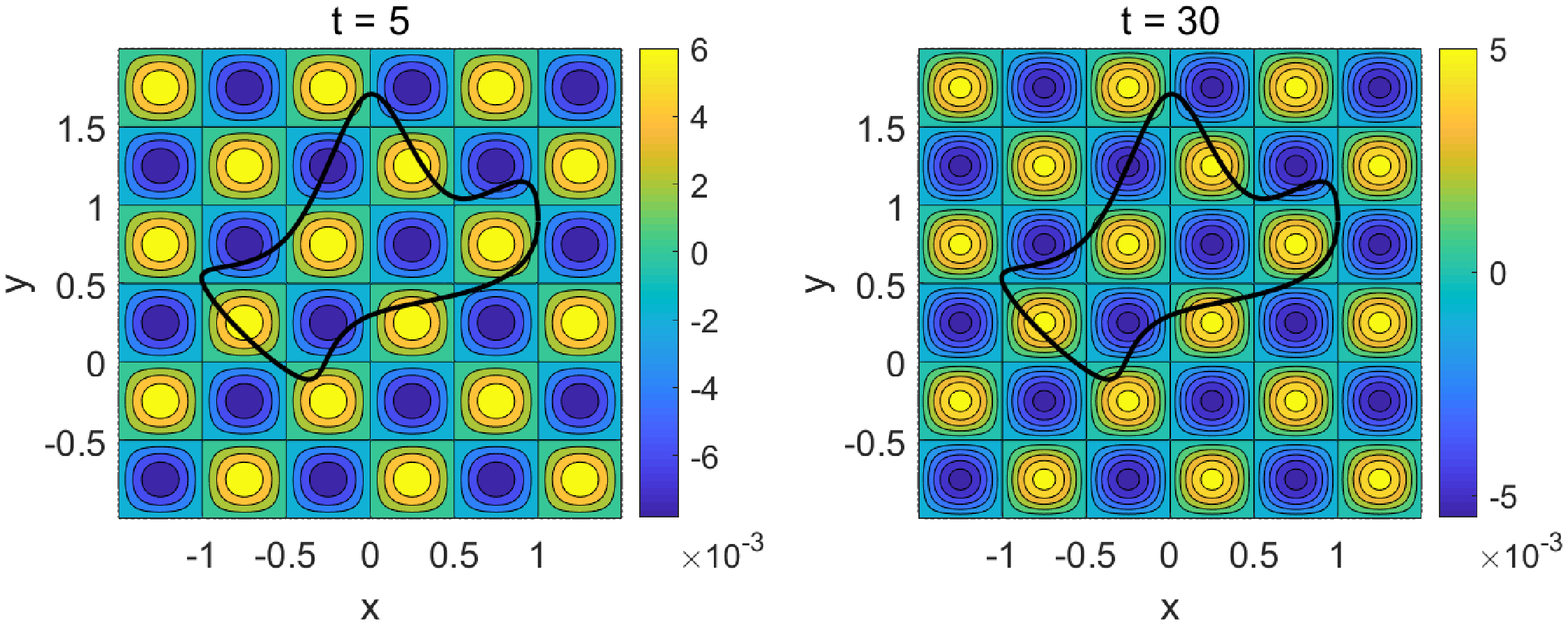,width=13cm}\par {(a) $\alpha=\beta=0.2$.}
\end{minipage}
\begin{minipage}{1\textwidth}\centering
\epsfig{figure=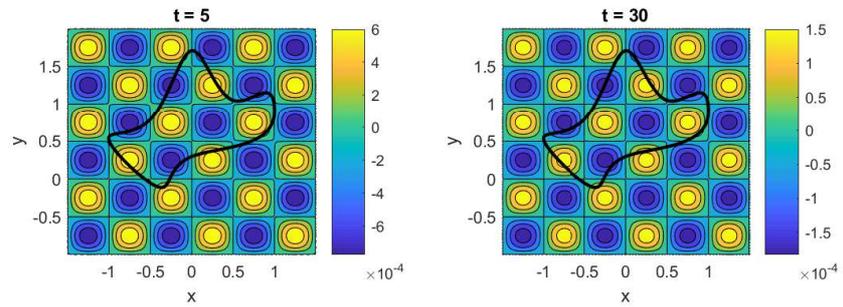,width=13cm}\par {(b) $\alpha=\beta=0.8$.}
\end{minipage}
\begin{minipage}{1\textwidth}\centering
\epsfig{figure=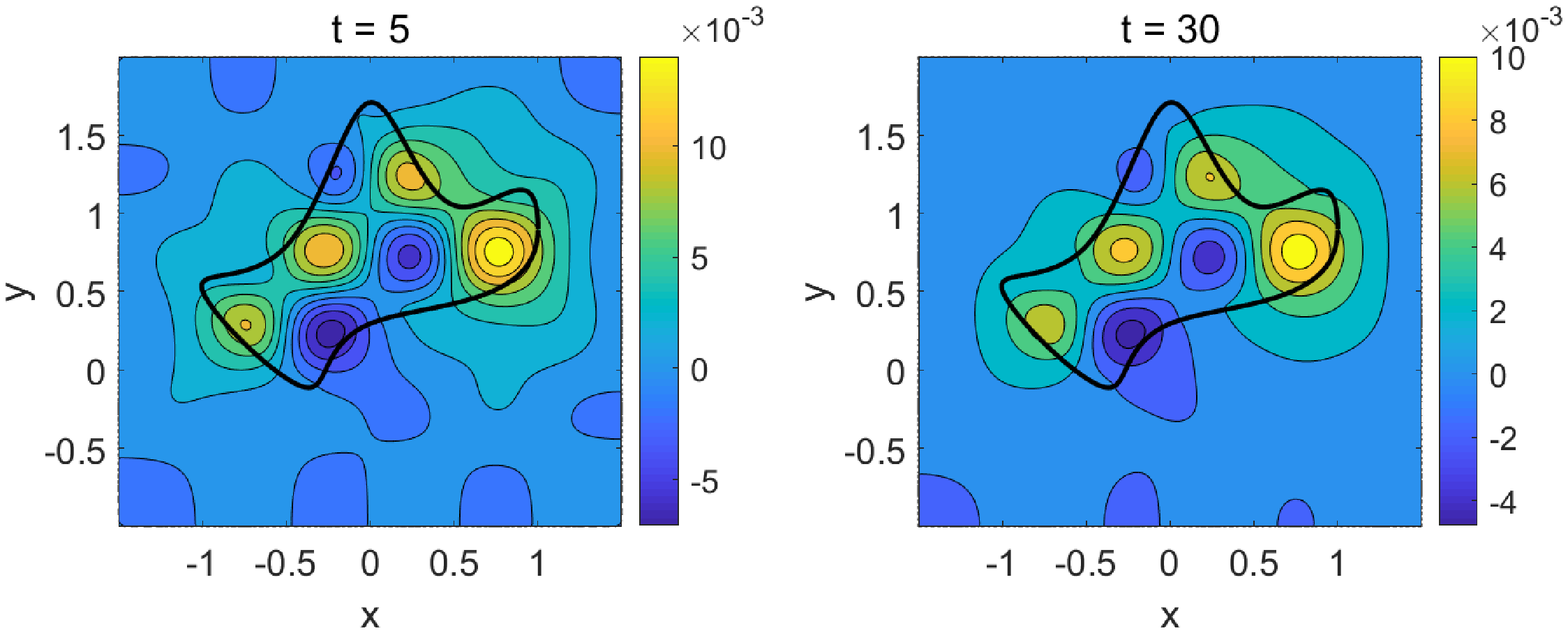,width=13cm}\par {(c) $\alpha=0.2,\beta=0.8$.}
\end{minipage}
\begin{minipage}{1\textwidth}\centering
\epsfig{figure=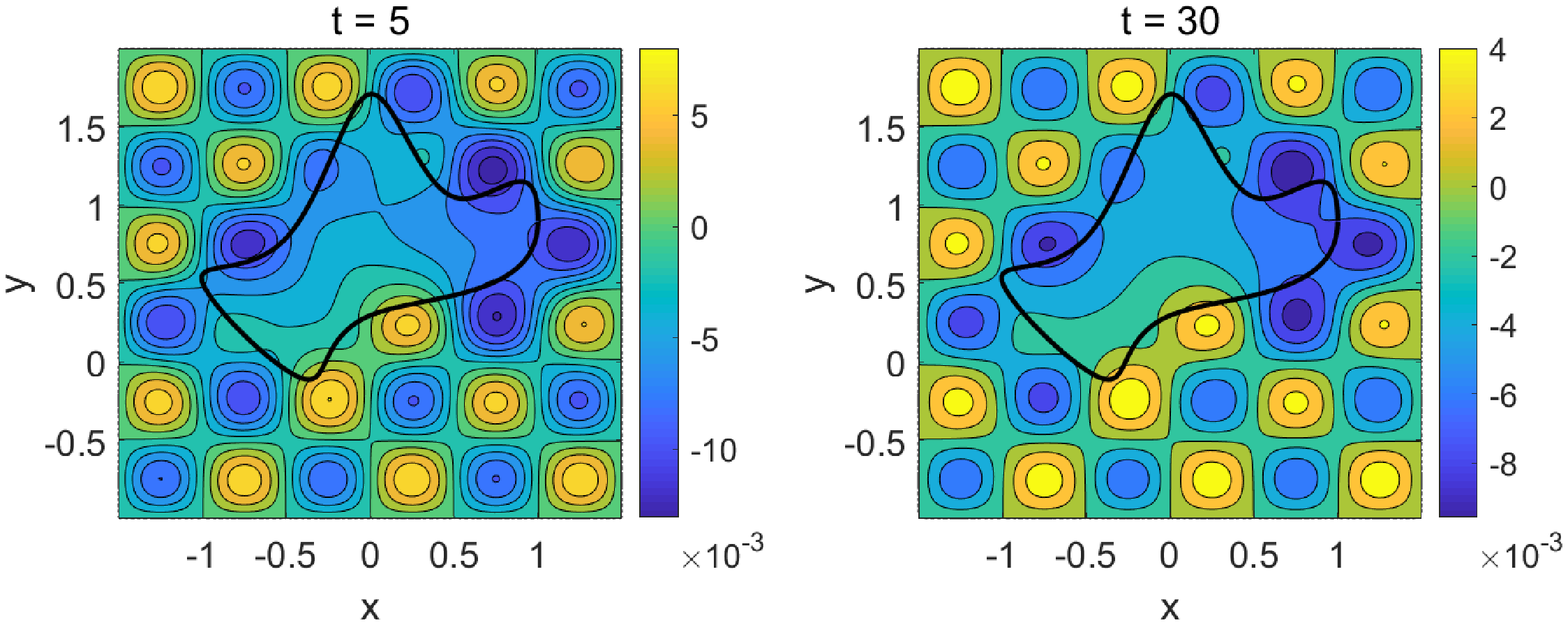,width=13cm}\par {(d) $\alpha=0.8,\beta=0.2$.}
\end{minipage}
\end{center}
\caption{Numerical solutions for Example \ref{s5-eg-4},  Case II, $\mu=\nu=1$.\label{eg2fig1}}
\end{figure}

\begin{figure}[!h]
\begin{center}
\begin{minipage}{1\textwidth}\centering
\epsfig{figure=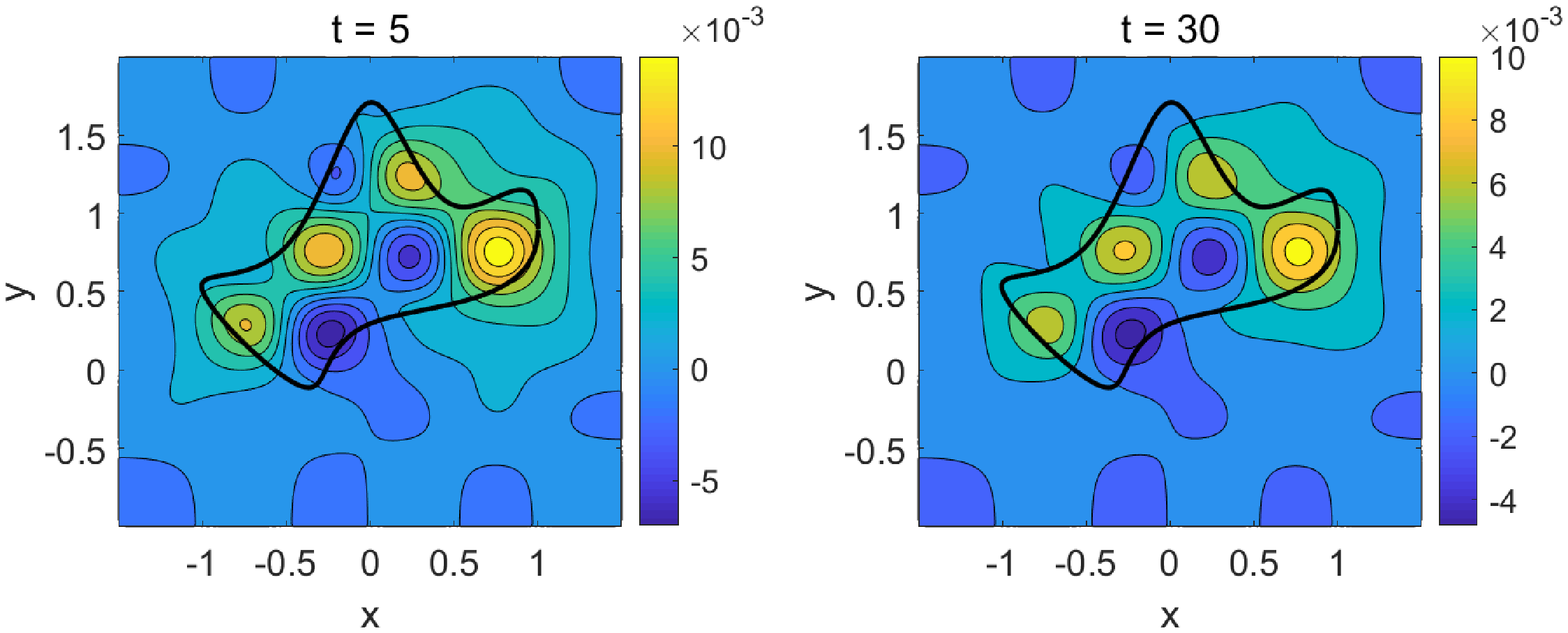,width=13cm}\par {(a) $\alpha=\beta=0.2$.}
\end{minipage}
\begin{minipage}{1\textwidth}\centering
\epsfig{figure=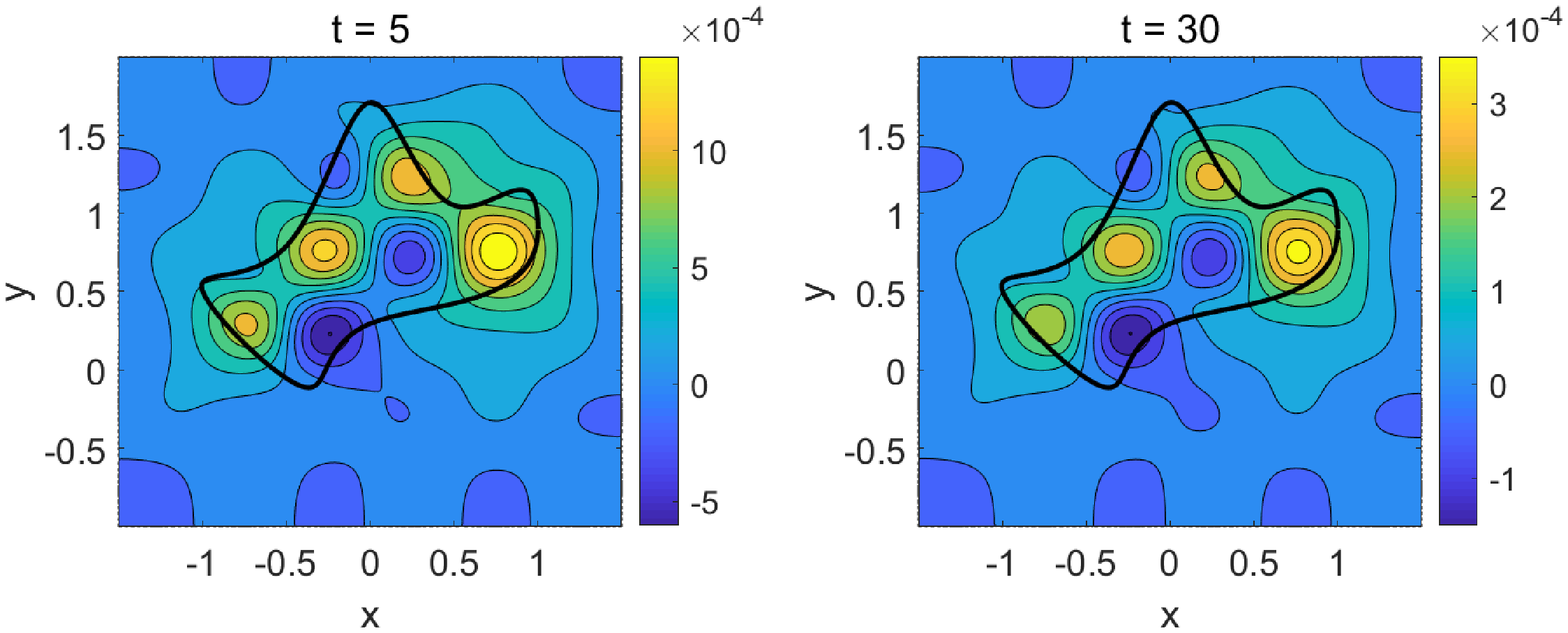,width=13cm}\par {(b) $\alpha=\beta=0.8$.}
\end{minipage}
\end{center}
\caption{Numerical solutions for Example \ref{s5-eg-4},  Case II, $\mu=1,\nu=10$.\label{eg2fig5}}
\end{figure}

\sjwchange{For Case II, we do not have the the analytical solutions, so we
exhibit only the numerical solutions in Figures~\ref{eg2fig1}
and \ref{eg2fig5}. We set $\tau=2^{-7}$, $N=80$, choose oscillating
initial conditions, and set the source terms to zero.}

First, we fix the diffusion coefficients $\mu=\nu=1$ and see how the
fractional orders $\alpha$ and $\beta$ affect the solution
of \eqref{s5eg1:eq-3}.  For $\alpha=\beta=0.2$, both the solutions $u$
(inside the black curve) and $v$ (outside the black curve) decay to
two oscillating patterns as time $t$ increases; see
Figure~\ref{eg2fig1}(a). With $\alpha=\beta$, both $u$ and $v$ decay
faster to the corresponding oscillating patterns as $\alpha$ becomes
larger; see Figure~\ref{eg2fig1}(b). When we choose a smaller
$\alpha=0.2$ and a large $\beta=0.8$, we observe that $u$ decays
slower than $v$; see Figure~\ref{eg2fig1}(c).  For a larger $\alpha$
value of $0.8$ and a smaller $\beta$ value of $0.2$, we observe th
faster decay of $u$ and slower decay of $v$; see
Figure~\ref{eg2fig1}(d).

Second, we choose a larger coefficient $\mu=10$ and a smaller
coefficient $\nu=1$ in the computation.  For $\alpha=\beta=0.2$, we
observe that $u$ decays faster than $v$ because of a larger $\mu$; see
Figure~\ref{eg2fig5}. \sjwchange{Closer} observation reveals that
increasing $\mu$ is somewhat equivalent to increasing $\alpha$, for
fixed $\beta$ and $\nu$; see Figure~\ref{eg2fig1}(c) and
Figure~\ref{eg2fig5}.  Similar behavior is observed for
$\alpha=\beta=0.8$, but $u$ (or $v$) decays faster than $u$ (or $v$)
for the case $\alpha=\beta=0.2$, see Figure~\ref{eg2fig5}(b).

\begin{figure}[!h]
\begin{center}
\begin{minipage}{1\textwidth}\centering
\epsfig{figure=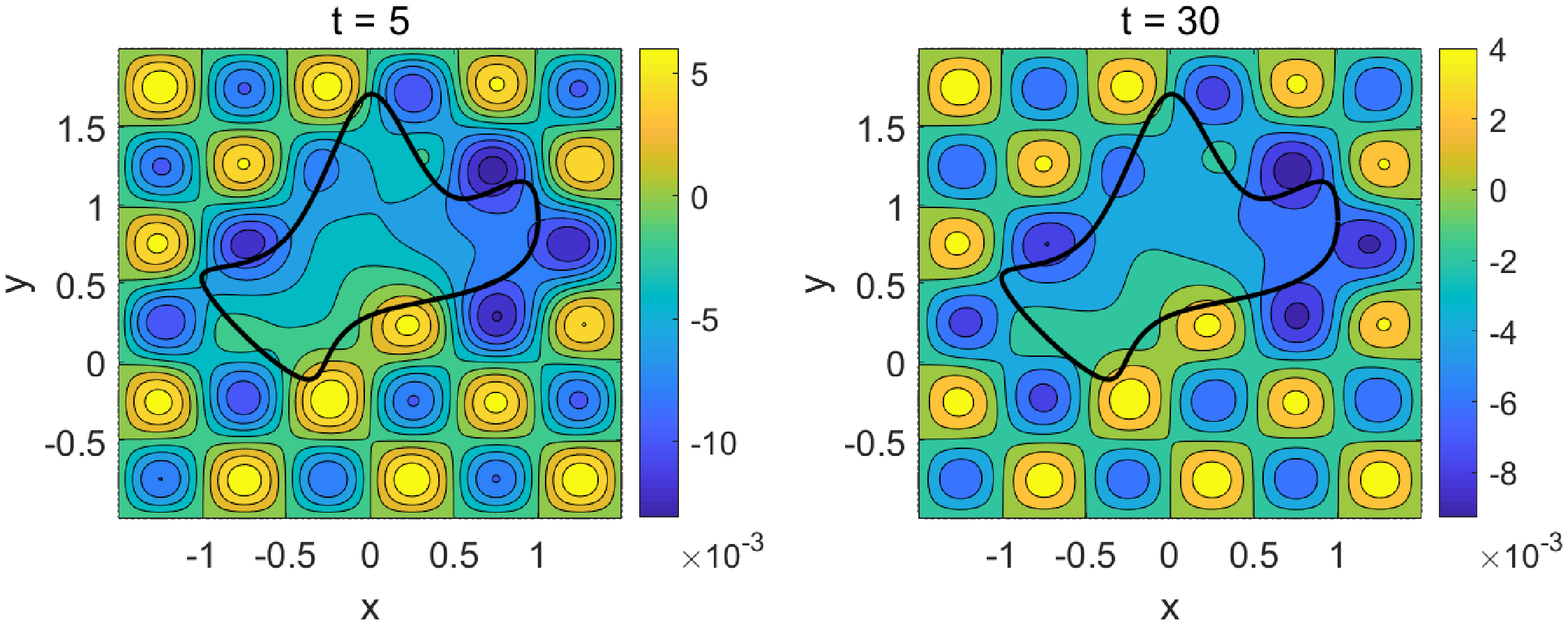,width=13cm}\par {(a) $\alpha=\beta=0.2$.}
\end{minipage}
\begin{minipage}{1\textwidth}\centering
\epsfig{figure=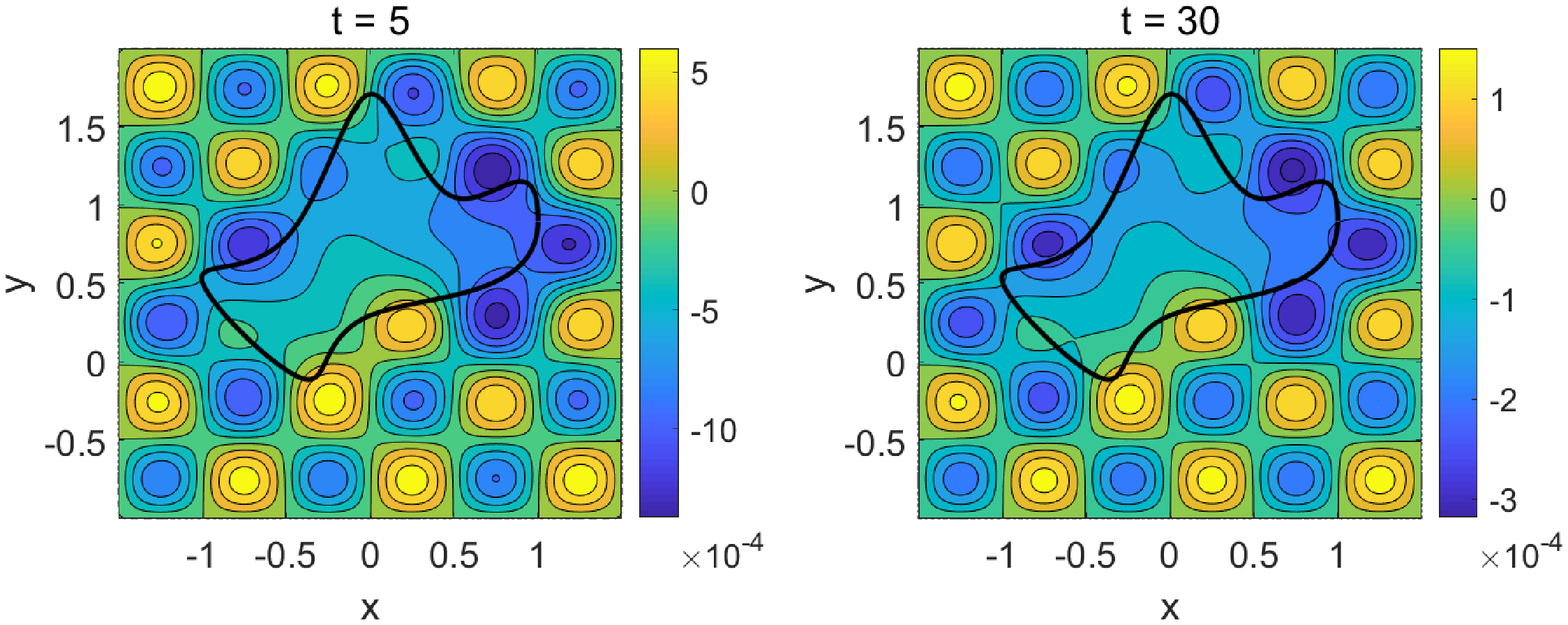,width=13cm}\par {(b) $\alpha=\beta=0.8$.}
\end{minipage}
\end{center}
\caption{Numerical solutions for Example \ref{s5-eg-4},  Case II, $\mu=10,\nu=1$.\label{eg2fig7}}
\end{figure}

Last, we choose a smaller diffusion coefficient $\mu=1$ and a larger
diffusion coefficient $\nu=10$.  We observe faster
decay \sjwchange{for $v$ than for $u$ when $\alpha=\beta$ (see
Figures~\ref{eg2fig7}(a) and (b)), due to the larger value of $\nu$}.
Larger values of $\alpha$ and $\beta$ (with $\alpha=\beta$) lead to
faster decay of both $u$ and $v$.  Compared with
Figure~\ref{eg2fig1}(c), we observe that increasing $\nu$ is somewhat,
but not precisely, equivalent to increasing the fractional order
$\beta$, for fixed $\alpha$ and $\mu$.



%

%
%

%
%

\section{Conclusion and discussion}\label{concl}

\sjwchange{We have  proposed} the regularized LSC method and LSFVM for
solving nonlinear time-dependent PDE systems with reaction terms.  In
the LS method, no grid generalization is required, which makes for
straightforward implementation of LS techniques.  The present method
can be directly extended to three-dimensional problems and can deal
with moving boundary conditions \sjwchange{at low cost, since the mass
and stiffness matrices ($\mathbf{A}_{in}$ and $\mathbf{S}_{in}$
in \eqref{s31:A}, respectively) do not need to be recomputed entirely
with a change of boundary. That is, we just need to add or delete some
rows and columns of these matrices when the boundary
changes. Extensions} to three-dimensional problems and problems with
moving boundary conditions will be studied in future work.

The use of the regularization reduces the condition number of the LS
\sjwchange{method, but it is still somewhat large,} especially when solving the
normal equation (see \eqref{s223-2}) to obtain the numerical
solutions. Fortunately, some methods have been proposed to solve
ill-conditioned linear systems accurately;
see \cite{CarsonHigham18,EldenSimoncini12,Neumaier98,ScottTuma17}.  An
alternative approach is solve the equivalent first-order system
of \eqref{s31:eq-1}, \sjwchange{for which the conditioning is
approximately the square root of the normal equations.}  This approach
will be explored further in our future work, drawing on the related
works in \cite{CaiLMM94,CaiMM97,Rekatsinas2018,ZeiLaiPin95}.
{  In this work, we have considered only the use of tensor
product Hermite cubic spline basis functions.  In the future work, we
will employ other basis functions, such as the tensor product B-spline
basis functions \cite{HughesCB05}. We will also consider the local
refinement to resolve the corner singularity of the
solution \cite{HughesCB05}.  }


\def\cprime{$'$} \def\cprime{$'$} \def\cprime{$'$} \def\cprime{$'$}

\end{document}